\documentclass[final]{siamltex}
\usepackage{amssymb,amscd,amsmath,fancybox,float,epsfig}
\usepackage{afterpage}
\usepackage{appendix}
\usepackage{hyperref}
\usepackage{graphicx}
\usepackage{color}
\usepackage{caption}
\usepackage{subcaption}
\usepackage{times}

\newcommand{\bigO}{{\mathcal O}}         

\newcommand\epsmach{\epsilon_{\operatorname{mach}}}

\newcommand\Span{\operatorname{span}}

\newcommand\cM{\mathcal M}

\newcommand\cJ{\mathcal J}

\newcommand\bone{\boldsymbol 1}

\newcommand\ba{\boldsymbol a}
\newcommand\be{\boldsymbol e}
\newcommand\bg{\boldsymbol g}
\newcommand\bh{\boldsymbol h}
\newcommand\bzero{\boldsymbol 0}
\newcommand\balpha{\boldsymbol \alpha}

\newcommand\btheta{\boldsymbol \theta}

\newcommand\bj{\boldsymbol j}

\newcommand\btau{\boldsymbol \tau}

\newcommand\bm{\boldsymbol m}

\newcommand\bk{{\boldsymbol k}}     
\newcommand\bv{\boldsymbol v}
\newcommand\bn{\boldsymbol n}
\newcommand\bnu{\boldsymbol \nu}
\newcommand\br{\boldsymbol r}

\newcommand\bF{\boldsymbol f}
\newcommand\BF{\boldsymbol F}
\newcommand\bx{\boldsymbol x}
\newcommand\bcc{\boldsymbol c}

\newcommand\cbf{\boldsymbol{ \check{f}}}

\newcommand\by{\boldsymbol y}

\newcommand\hbf{\widehat{\boldsymbol f}}

\newcommand\bl{\boldsymbol l}

\newcommand\bV{\boldsymbol V}

\newcommand\cF{\mathcal F}

\newcommand\Beta{\boldsymbol \beta}

\newcommand\hf{\widehat{f}}
\newcommand\hg{\widehat{g}}

\newcommand\hrho{\widehat{\rho}}

\newcommand\TA{\mathbb A}
\newcommand\hTA{\widehat{\TA}}

\newcommand\cC{\mathcal{C}}

\newcommand\bZ{\boldsymbol {Z}}

\renewcommand\Re{\operatorname{Re}}

\newcommand\bbC{\mathbb C}

\newcommand\bbN{\mathbb N}

\newcommand\bbR{\mathbb R}
\newcommand\bbZ{\mathbb Z}

\newcommand\pa{\partial}

\newcommand\restrictedto{\upharpoonright}

\newcommand\Id{\operatorname{Id}}

\newtheorem{example}{Example}[section]

\newtheorem{remark}{Remark}[section]

\title{Geometry of the Phase Retrieval Problem}
\author{Alexander Barnett\thanks{
    Center for Computational Mathematics (CCM), Flatiron Institute, 
    162 Fifth Avenue, New York, NY 10010. E-mail:
    {abarnett@flatironinstitute.org}, {jmagland@flatironinstitute.org}}
    \and
Charles L. Epstein\thanks{
    Dept.~of Mathematics, University of Pennsylvania,
    209 South 33rd Street, Philadelphia, PA 19104, and CCM, Flatiron Institute. E-mail:
    {cle@math.upenn.edu}.
    } \and
 Leslie Greengard\thanks{Courant Institute,
    New York University, 251 Mercer Street, New York, NY 10012, and
    CCM, Flatiron Institute, 
    162 Fifth Avenue, New York, NY 10010.
    E-mail: {lgreengard@flatironinstitute.org}.\hfill Printed: \today
    } 
\and Jeremy Magland\footnotemark[1]}

\begin{document}

\maketitle

\begin{abstract}
  One of the most powerful approaches to imaging at the nanometer
  length scale is coherent diffraction imaging using X-ray
  sources.  For amorphous (non-crystalline) samples, raw data collected
  in the far-field can be interpreted as the modulus of the
  two-dimensional
  continuous Fourier transform of the unknown object.
  The goal is then to recover the phase through computational means
  by exploiting prior 
  information about the sample (such as its support),
  after which the
  unknown object can be visualized at high resolution.  While many
  algorithms have been proposed for this {\em phase retrieval} problem,
  careful analysis of its well-posedness has received relatively little
  attention.  In this paper, we show that the problem is, in general,
  not well-posed and describe some of the underlying issues that are
  responsible for the ill-posedness. We then show how this analysis can
  be used to develop experimental protocols that lead to better
  conditioned inverse problems.
\end{abstract}



\section{Introduction}
\label{s:intro}

With the increased accessibility of high energy, coherent light
sources, there has been a resurgence of interest in inverse imaging
problems where the measurement can be interpreted as the modulus of
the Fourier transform of an unknown object (see, for example,
\cite{chapman2006,elser2003,elser2007,fienup1987,miao1999,miao2015,Osherovich,bendory17}).
In a typical set-up, an object with electron density $\rho(\bx)$ is
irradiated with a planar beam of coherent x-rays and a measurement of
the scattered wave is made in the far field, i.e.\ the Fraunhofer
regime.  The wave's phase information is not directly measurable, and
the data collected are interpreted as samples of the intensity of the
Fourier transform, $\{|\hrho(\xi_j)|^2:\:j\in \cJ\}$. The
reconstruction problem is then largely reduced to that of
``recovering'' the phases of the complex numbers
$\{\hrho(\xi_j):\:j\in \cJ\}$.  This experimental approach is referred
to as coherent diffraction imaging (CDI).

\begin{remark}
The {\em phase retrieval} problem came to prominence in x-ray crystallography,
where $\rho(\bx)$ is assumed to be a periodic function with some unit cell and
$\{\hrho(\xi_j)\}$ are the coefficients of a discrete Fourier series, defined
only on a regular lattice.  In this context, without additional information,
the phase retrieval problem is obviously ill-posed.  One can assign any value to the phase
of each $\hrho(\xi_j)$ to produce a periodic image.  To circumvent this
problem, for sufficiently small molecules, direct methods that rely on the
non-negativity of the electron density and algebraic relations (Karle-Hauptman
determinants) proved to be very powerful \cite{millane1990}.  For larger
structures, a variety of experimental approaches have been introduced to supply
additional information that permits the reconstruction of the phase information
needed to reconstruct the original crystal \cite{LaddPalmer}.
\end{remark}

In the present paper, we are interested in the setting where
$\rho(\bx)$ is an essentially arbitrary, but compactly supported,
function.  The study of this problem dates back to 1952, when Sayre
noted that for amorphous {\em non-crystalline} objects, one may obtain
values of the intensity $|\hrho(\xi_j)|^2$ on a finer mesh than in the
periodic case \cite{Sayre1952}, since the spectrum is continuous.  In
essence, he proposed that one {\em ``over-sample''} $|\hrho(\xi)|^2$ by
a factor of two in each direction then recover $\rho(\bx)$ as the solution to an
overdetermined, constrained nonlinear least squares problem. The
constraints come from some prior knowledge of $\rho$, such as its
support, whether it is non-negative, etc.  It turns out that Sayre's
conjecture is essentially correct in more than one spatial
dimension. More precisely, a finite approximation to phase retrieval
problem introduced below in Section~\ref{s.dcpp}, with support as the
auxiliary information, has a solution, generically unique up to
``trivial associates;'' see~\cite{barakat1984,hayes1982,hayes1987,bendory17}.
In the continuum case, the trivial associates are obtained by applying
operations to $\rho$ that leave $|\hrho|$ invariant: translations
($\rho(\bx) \rightarrow \rho(\bx-\ba)$ for some $\ba$) and inversion
($\rho(\bx) \rightarrow \rho(-\bx)$).

The most common additional information used in phase retrieval is a
support constraint: that $\rho(\bx)$ is nonzero only within some
closed and bounded region $D$ in ${\mathbb R}^d$. It is with reference
to an estimate for the support of $\rho(\bx)$ that one speaks about
over-sampling its magnitude Fourier transform.  If $R\supset D$ is the
smallest rectangle covering $D$, then the Fourier data
$\{|\hrho(\bk_j)|\}$ must be sampled on a grid fine enough to
represent a periodic function with fundamental cell a rectangle whose
side lengths are twice those of $R$. In the phase retrieval
literature, this is called ``double oversampling.'' In fact some
degree of oversampling is clearly needed so that the sampled data 
contains adequate information about the support of $\rho(\bx)$ for
reconstruction to be possible. Double oversampling is the assumption
required for the proof of the basic uniqueness theorem (Hayes's
theorem) in coherent diffraction imaging; it also allows the autocorrelation function of $\rho,$
\begin{equation}
  \rho\star\rho(\bx)=\int\rho(\bx+\by)\rho(\by)d\by,
\end{equation}
to be reconstructed from the magnitude Fourier data without aliasing artifacts.

The earliest practical method for solving the phase retrieval problem
is a variant of the alternating projection algorithm due to Saxton and
Gerchberg \cite{gerchberg1972}.  The basic idea, which was first
introduced in the context of Banach spaces by von Neumann as a method
to find the intersections of convex sets, is quite general. In phase
retrieval we let $A$ denote the collection of images $\rho$ with the
given magnitude Fourier data and let $B$ denote the set of images
which satisfy the support constraint.  Given a function $f(\bx)$,
projection onto $A$ corresponds to computing its Fourier transform
$\hf$, keeping the phase information from $\hf$ and replacing the
modulus with the measured data $|\hrho(\xi)|$.  We denote this
operator by $P_A$.  Projection onto $B$ corresponds to multiplying
$f(\bx)$ by the characteristic function of $D$. We denote this
operator by $P_B$.  Alternating projection can then be written as the
following iteration:
\begin{equation}
\rho_{k+1} = P_B\circ P_A(\rho_k), 
\label{apiter}
\end{equation}
with some initial guess $\rho_0$.

This algorithm has a long history when $A$ and $B$ are convex sets,
which we do not 
seek to review here (see, for example, \cite{bauschke1996}). It has
also received a lot of study in the non-convex setting 
\cite{andersson2013,bauschke2002,borwein}.
Unfortunately, alternating projection often
converges to fixed points unconnected to the reconstruction problem at
hand. To overcome this, Fienup proposed a new class of so-called
hybrid input-output (HIO) algorithms 
\cite{bauschke2002,fienup1982,fienup1987}, which
were placed into the larger framework of \emph{difference-maps} by
Elser and collaborators, see~\cite{elser2003,elser2007}. 
The fixed points
of these algorithms all specify correctly reconstructed objects.
We will describe this method in detail below in sections
\ref{sec4} and \ref{sec5}.
Note that methods from continuous optimization have also been
applied to this problem; see \cite{Osherovich}.

Despite the enormous effort that has gone into finding robust algorithms,
the state of the art is generally unsatisfactory and reconstructions
are typically not very accurate.  That is to say, the
phase retrieval problem with a support constraint has all the hallmarks 
of an ill-posed problem. Our
main purpose in this paper is to describe recent work aimed at understanding
what aspects of the phase retrieval problem render it ill-posed, and how this
knowledge can be used to modify the experimental protocols to obtain better
conditioned inverse problems. 

\begin{remark}
The {\em phase retrieval problem}
sometimes refers to the more general setting where $\rho(\bx)$ is
unknown and measurements $M(k)$ are  of the form
\[ M(k) = | \langle \rho, a_k \rangle |,    \] 
for some set of querying functions $a_k$. Here, 
$\langle \rho, a \rangle$ denotes the inner product of the two
functions.
If the phase information
were available, then solving for $\rho$ would correspond 
to a linear least squares problem. 
Without the phase information, the problem is non-convex.
When the map from $\rho$ to $M(k)$ is invertible, a variety of 
optimization methods have been developed based, for example, on
semidefinite relaxation or gradient descent 
\cite{candes2015a,candes2015b}.
Unfortunately, the phase retrieval problem of interest in x-ray
scattering does not satisfy the necessary hypotheses for these methods to 
apply, namely that the forward map is injective and the solution is unique. 
The recent paper \cite{alaifari2017} contains a detailed analysis of phase retrieval
in the invertible case and an interesting discussion of stability in that context.
\end{remark}

\subsection{The Discrete Classical Phase Retrieval Problem}\label{s.dcpp}
For the sake of simplicity, we analyze a finite-dimensional analogue
of the phase retrieval problem described above, which, in the limit of
infinitely many samples, converges to the continuum problem.
We assume that $\rho$  is \emph{real valued}, and imagine that the
unknowns are the samples $f_{\bj}=\rho\left(\frac{\bj}{N}\right),$
where $\bj\in J$ are points in a finite cubical integer lattice,
$J=\{0,1,\dots,2N-1\}^d \subset \bbZ^d$.
Here $d$ is the ambient dimension ($d=2$ in our examples, but
3D phase retrieval is also possible \cite{chapman2006}).
The vector $\bF := (f_{\bj}:
\bj\in J)$ denotes an image, which can be viewed as a uniform
pixelization of a density function $\rho$ lying in $[0,2)^d$.  We use
  the notation $\bbR^J$ to denote the set of all possible such
  images. Using $J$ as the index set is somewhat non-standard in the
  engineering literature.

The measured data values are
modeled as $a_\bk := |\hat f_{\bk}|$,
where
\begin{equation}
  \hat f_{\bk}=\sum_{\bj\in J}f_{\bj}\exp\left(\frac{2\pi i\bj\cdot\bk}{2N}\right)~,
  \qquad \bk\in J~,
  \label{dft}
\end{equation}
is the usual $d$-dimensional discrete Fourier transform (DFT)
taking the $(2N)^d$ pixel values to $(2N)^d$ frequency data.
We call $a_\bk$ \emph{magnitude DFT} data, and denote the data vector
by $\ba := (a_\bk: \bk\in J)$. We define the measurement map,
$\cM:\bbR^J\to\bbR_+^J,$ by setting
\begin{equation}
  \cM(\bF):=(|\hf_{\bk}|:\:\bk\in J).
\end{equation}
\begin{remark}
  One may connect the above discrete model to the continuous case as follows.
  From \eqref{dft}, the indices $\bk$ are $2N$-periodic in each
  dimension.  Let $\tilde\bk$ be the periodic folding of $\bk$ into
  the origin-centered cube $\{-N,-N+1,\dots,N-1\}^d$.  Then define the
  spatial frequencies $\xi_\bk := \pi \tilde{\bk}$ for $\bk\in J$,
  which lie in the cube $[-\pi N,\pi N]^d$.  The sum in \eqref{dft} can be
  interpreted as approximate samples of $(2N)^d\hat\rho(\xi_\bk)$,
  where the Fourier integral over $[0,2)^d$ has been approximated by a
    $2N$-point trapezoid quadrature (in each dimension)
    at the nodes $\bj/N$.  Thus for
    a continuous function $\rho$ the sequence $(2N)^{-d} \hat f_{\bk}$
    tends to the exact Fourier transform as $N\to\infty$.
An alternative (but less physically realistic) interpretation is:
  \eqref{dft} gives point samples of the {\em exact}
  Fourier transform of a ``sum of point masses'' scatterer model
  $\rho(\bx) := \sum_{\bj\in J} f_{\bj} \delta(\bx-\bj/N)$.
\end{remark}

The advantage of a discrete model over the continuous one is that it
admits, given a support condition to be presented shortly, an exact
solution that is generically unique up to trivial associates. In
contrast, for the continuous $\rho$ problem, given any finite
collection of samples of $|\hrho(\bk)|$, there is an infinite
dimensional space of functions with these Fourier coefficients, which
also satisfy the support constraint. 

\begin{definition}
  Given a magnitude DFT data vector $\ba = (a_\bk: \bk\in J)$,
  the {\em magnitude torus}, denoted by $\TA_{\ba}$, is the 
  collection of images $\bF$ in $\bbR^J$ with this magnitude DFT data,
  i.e.
  \begin{equation}
    \TA_{\ba}=\{\bF\in\bbR^J:\:|\hat f_\bk| = a_\bk \text{  for all }\bk\in J\}.
  \end{equation}
\end{definition}
\begin{remark}
Note that $\TA_{\ba}$ is either empty (if $\ba$ does not obey the
inversion symmetry demanded by \eqref{dft} for a real image), or is a
real torus, or union of tori, of dimension equal to approximately half
of the cardinality of $J$.  The reason for the approximate nature, and
the possible existence of multiple connected components (which are all
tori), is the fact that some data is forced to obey various
symmetries; for example $\hat f_{\mathbf{0}}$ is always real, whereas
most DFT data is generically complex. Vanishing DFT coefficients
also lower the dimension of the torus.
\end{remark}

The prior information about the image is encoded as a second set
$B\subset \bbR^J$. The \emph{discrete classical phase retrieval problem} is
then the problem of finding points in the intersection $\TA_{\ba}\cap B$;
see Fig.~\ref{tangentfig}.  For an image $\bF$ we denote its
true support by
    \begin{equation}
      S_{\bF}:=\{\bj\in J:\:f_{\bj}\neq 0\}.
      \label{Sf}
    \end{equation}
If $S_{\bF} \subset S\subset J,$ then we say that
$S$ is an estimate for the support of $\bF$, and let
\begin{equation}
  B_S:=\{\bF\in\bbR^J:\: f_{\bj}=0\text{ for }\bj\notin S\}.
\end{equation}
This is clearly a linear subspace of $\bbR^J$.

To define the operations in the finite, discrete case that generate the
set of trivial associates of an image $\bF$ we need to extend the
image to be $2N$-periodic. That is, the integer indices are defined mod
$2N$ in each dimension. With this understood, the image $f_{\bj}$ is
defined for all $\bj\in\bbN^d,$ with its restriction to $\bj\in J$
representing a single period. We call such images $J$-periodic.  This
is consistent with the formula for the DFT, which defines $\hf_{\bk}$
for all $\bk\in\bbN^d,$ and the inverse formula, which defines
$f_{\bj}$ for all $\bj\in\bbN^d,$ as  $J$-periodic images.

With this periodic extension, the
operations that generate the trivial associates are:
\begin{enumerate}
\item If $\bv\in J$, then the {\bf translate} of $\bF$ by $\bv$ is defined
  by its components
  \begin{equation}\label{eqn2}
    f^{(\bv)}_{\bj} := f_{\bj-\bv}~, \qquad \bj\in J~.
  \end{equation}

\item The {\bf inversion} $\check{f}$ is defined by its components
  \begin{equation}\label{eqn3}
    \check{f}_{\bj} := f_{-\bj}~, \qquad \bj\in J~.
  \end{equation}

\end{enumerate}
\noindent
Each image in the set of all $2(2N)^d$ trivial 
associates has the same magnitude DFT data as $\bF$.

We now state a well-known uniqueness theorem due to Hayes
\cite{barakat1984,hayes1982,hayes1987}, concerning the support
constraint.  If $S$ is contained in a rectangular subset of $J$, with
side lengths at most half the corresponding side-lengths of $J$,
i.e.\ at most $N$, and $S_{\bF} \subset S$, then the set
$\TA_{\ba}\cap B_S$ is finite, and generically consists of trivial
associates of a single point in this set.

\begin{definition} \label{ssdef}
  An image $\bF$ with support $S_{\bF}\subset S$, for a set $S$
as above, is said to have \emph{small support}.
\end{definition}

Note that small support corresponds to the density function
$\rho$ having support lying within a $d$-dimensional cube of side length 1.


In some experimental situations one has a constraint
$\bF\in B_+$, where
\begin{equation}
  B_+:=\{\bF:\: f_{\bj}\geq 0,\text{ for all }\bj\in J\}.
  \label{Bplus}
\end{equation}
 This auxiliary condition alone does not uniquely specify a set of trivial
 associates, or even a finite set. It is easy to show, however, 
that the squared magnitude data
 $(|\hat f_{\bk}|^2:\:\bk\in J)$ are the DFT coefficients of the
 autocorrelation image
\begin{equation}
  [\bF\star\bF]_{\bj}=\sum_{\bl\in J}f_{\bl}f_{\bj+\bl}.
\end{equation}

In~\cite{BEGM} we prove that if the support of $\bF\star\bF$ is
sufficiently small, then the set of non-negative images in $\TA_{\ba}$
is finite, and generically consists of trivial associates of a single
element. This follows because, for $\bF\in B_+$, if the support of
$\bF\star\bF$ is sufficiently small, as a subset of $J,$ then it
provides a non-trivial upper bound on the support of $\bF$, for which
we can conclude uniqueness, up to trivial associates.  As all points
$\bF\in\TA_{\ba}$ have the same autocorrelation image, without a
non-negativity constraint, the support of the autocorrelation image
does \emph{not}, in general, provide a bound on the support of the
image itself.  We  refer to the problem of finding the
intersection $\TA_{\ba}\cap B_+$ as {\em phase retrieval with
  non-negativity constraints}.

\begin{definition} \label{adequatedef}
We say that the set $B,$ defined by some auxiliary conditions, 
is {\em adequate} if the intersection $\TA_{\ba}\cap
 B$ is a finite set for any $\ba$ in the range of $\cM$.
(Thus, if the support of $\bF\star\bF$ is
sufficiently small, then non-negativity of the image is adequate
data for phase retrieval.)
\end{definition}

\subsection{Well-posedness of the Discrete Phase Retrieval Problem}

 The concept of well-posedness for an inverse problem comprises two
 distinct questions: the first is the uniqueness of the solution and
 the second concerns the continuity properties of the local inverse
 map near a solution.  For the phase retrieval problem, uniqueness
 should be understood as uniqueness up to trivial associates.

 Neither aspect of well-posedness has been analyzed in  detail for
 the phase retrieval problem.  From the proof of Hayes' uniqueness
 theorem, it follows that there are pairs of images, $\bF_1, \bF_2$,
 which are not trivial associates, with identical magnitude DFT data
 and support for which $\|\bF_1-\bF_2\|\approx 1$. Yet, since uniqueness is
 generic, we can find images, $\bF'_1, \bF'_2$, as near to $\bF_1,
 \bF_2$ as we like for which the phase retrieval problem does have a
 unique solution. Moreover we can assume that the supports of $\bF'_1,
 \bF'_2$ are the same as those of $\bF_1, \bF_2$. From this
 observation it is clear that, at any finite precision, this problem
 does not always have a unique solution, even up to trivial associates. A
 similar observation was made by Fienup and Seldin
 in~\cite{FienupSeldin:90}.

A primary concern here is that of understanding
what makes it difficult to find points in $\TA_{\ba}\cap B,$ even if
it is assumed that the set $B$ is selected so that this intersection
consists of finitely many points.  In coherent diffraction imaging, 
the cardinality of the index
set $J$ is in the hundreds of thousands, millions, or even
billions (depending on whether a two-dimensional or three-dimensional
object is being imaged and at what resolution). 
High dimensionality certainly complicates the problem at
hand, but it is not the root cause of its difficulty. Rather, it
is the geometry \emph{near to} points in $\TA_{\ba}\cap B$ that
renders this problem so difficult.      This sort of local geometry is usually discussed in
terms of the relationship of the fibers of the tangent bundles to
$\TA_{\ba}$ and $B$ at points of intersection. For the remainder of
this discussion we focus on the support condition case $B=B_S$,
which is itself a linear subspace.
\begin{definition}
For a point $\bF\in\TA_{\ba},$ the {\em fiber of the tangent bundle}
to $\TA_{\ba}$ at $\bF$, denoted by $T_{\bF}\TA_{\ba}$, is the affine
subspace of $\bbR^J$ through $\bF$ that is the best linear
approximation to $\TA_{\ba}$ near to $\bF$.
\end{definition}
\begin{remark}
  An embedded submanifold of a Euclidean space has a ``best''
  approximating affine subspace if it is at least $\cC^1.$ A
  magnitude torus is a product of round circles and is therefore a
  real analytic subspace of $\bbR^J.$
\end{remark}

\begin{definition}
Let $\bF\in \TA_{\ba}\cap B_S$.
The intersection
of $\TA_{\ba}$ with $B_S$ is {\em transversal} at $\bF$ if
\begin{equation}
  T_{\bF}\TA_{\ba}\cap B_S= \{\bF\},
\end{equation}
that is, the affine space
$T_{\bF}\TA_{\ba}$ intersects the linear subspace $B_S$ only at $\bF$.
\end{definition}

 If the intersection $\bF$ is transversal, then the geometry of
 $\TA_{\ba}\cup B_S$ near to $\bF$ is accurately modeled by a
 neighborhood of $\bF$ in $T_{\bF}\TA_{\ba}\cup B_S.$ In this case,
 the conditioning of the problem of finding a point $\bF\in
 \TA_{\ba}\cap B_S$ is determined by the angles between
 $T_{\bF}\TA_{\ba}$ and $B_S.$ If there are positive dimensional
 subspaces of $T_{\bF}\TA_{\ba}$ and $B_S$ that make a very small
 angle with one another, then the condition number, though finite, will be very large.

 For high dimensional non-linear submanifolds of $\bbR^J$ one does not
 generally expect to have an explicit description of the fibers of the
 tangent bundle. It is a remarkable feature of magnitude tori that such
 a description is accessible.  Using this description we can show that
  the intersections between $\TA_{\ba}$ and $B_S$ are {\em typically
   not transversal}.  If the intersection at $\bF$ is not transversal,
 then $T_{\bF}\TA_{\ba}\cap B_S$ is a positive dimensional affine
 subspace.  In this case the local geometry of $T_{\bF}\TA_{\ba}\cup
 B_S$ near to $\bF$ does not resemble that of $\TA_{\ba}\cup B_S;$ 
 the linearized problem does not have a unique solution. Formally
 speaking, the condition number of the non-linear problem is infinite.
 Moreover, in this case, linear analysis fails to adequately describe
 the behavior of algorithms for finding intersection points.

Above we defined $\cM:\bbR^J\to \bbR_+^J$ as the ``forward operator''
(in the language of inverse problems), i.e.\ the measurement map from
an image to its corresponding DFT magnitude data $\ba:=(|\hat
f_{\bk}|:\:\bk\in J)$.  The inverse image $\cM^{-1}(\ba)$ of a point
$\ba\in\bbR_+^J$ is simply its magnitude torus $\TA_{\ba}$.  Suppose
now that $\bF\in\TA_{\ba}$ has small support contained in the set
$S$. Then, by Hayes' theorem, the set $\cM^{-1}(\ba)\cap B_S$ is
finite, and non-empty. This remains true if $\bF$ is replaced by a
nearby point $\bF'\in B_S$.  Therefore, the map
$\cM\restrictedto_{B_S}$ has a local inverse, defined on the manifold
of consistent data $\cM(B_S)$, near to $\cM(\bF)$. We denote this
local inverse by $\cM^{-1}_{\bF,S}.$ As noted above, the other issue
that arises in a discussion of well-posedness is the continuity of
this local inverse.  In~\cite{BEGM} we prove the following theorem.
\begin{theorem}
  The local inverse $\cM^{-1}_{\bF,S}$ satisfies a Lipschitz estimate if 
  and only if the intersection of $\TA_{\ba}$ with $B_S$ at $\bF$ 
  is transversal.
\end{theorem}

\begin{figure}
\centering
\includegraphics[width=1.0\linewidth]{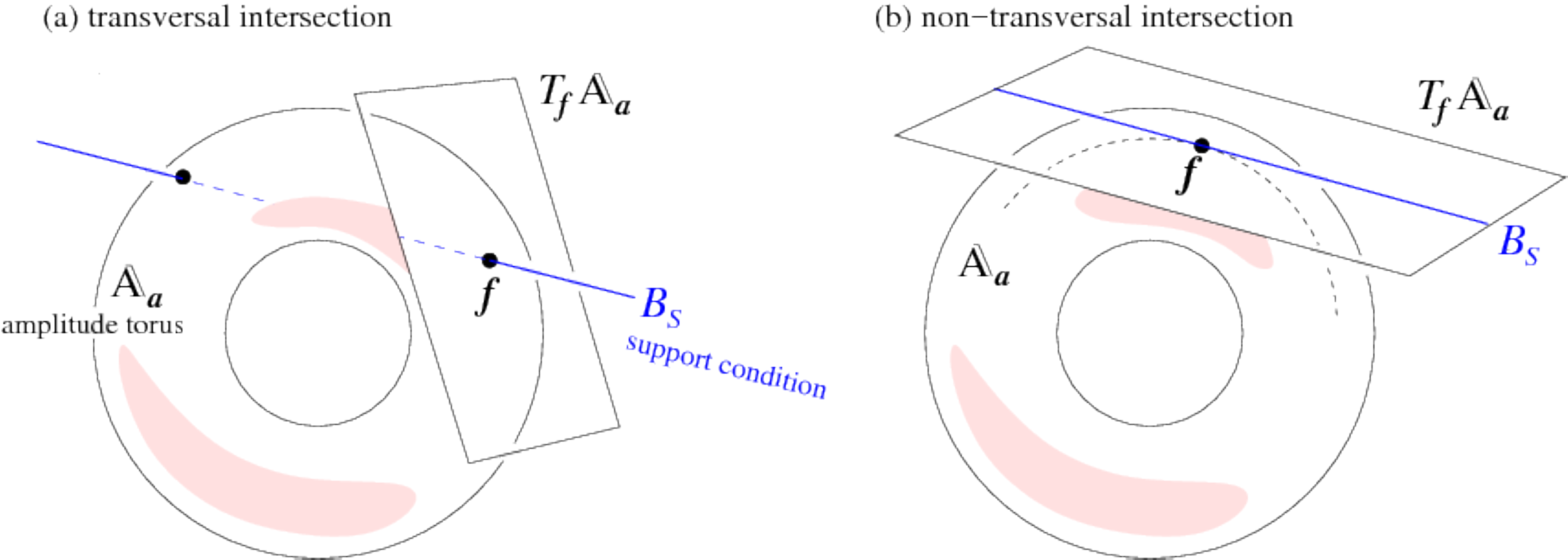}
\caption{ Illustration of two types of intersection in image space
  between the magnitude torus (where $\ba = \cM(\bF)$) and the
  constraint $B_s$.  (a) Transversal case. The angle between the
  tangent space and $B_S$ is positive. Also visible (left-most black
  dot) is a trivial associate of the image $\bF$.  Finding $\bF$ given
  the data $\ba$ is (locally) well-conditioned.  (b) Non-transversal
  case. The angle between the fiber of the tangent bundle at $\bF$ and
  $B_S$ is zero (since we sketch in $\mathbb{R}^3$ we are forced to
  show $B_S$ lying within $T_{\bF}\TA_{\ba}$; in general this is only
  true for a subspace of $B_S$). In the case shown, the distance from
  the torus grows quadratically with distance from $\bF$ for points in
  $B_S$.  High-order contact of this type is much more problematic in
  high dimensions, where the dimension of the fibers of the tangent bundle can be
  large.}
\label{tangentfig}
\end{figure}

By itself, 
the failure to have a local Lipschitz inverse leads to a kind of ill-conditioning.
When the intersection is non-transversal, the local inverse is, 
at best, H\"older continuous of order $\alpha <1$,
which implies an infinite condition number.
In fact the number of 
accurate digits possible in the reconstructed image
cannot exceed $\alpha d$ when the data is available with a
relative precision of $d$ digits.
For phase retrieval, it is often the case that $\alpha\leq \frac 1 2$
(see Fig. \ref{tangentfig}(b) for an illustration with $\alpha=\frac 1 2$). 
More critical, however, is that non-transversality
stalls the convergence
for standard reconstruction algorithms,
even in the ideal case of noise-free data,
as discussed in section \ref{sec4}.
  
%
%
\subsection{Contents of the Paper}
An outline of the paper follows: in sections \ref{sec2} and
\ref{sec3}, we discuss the geometry of phase retrieval with support
constraints.  We turn to the practical consequences of our analysis in
sections \ref{sec4} and \ref{sec5}, and extend the analysis to the
case of non-negativity constraints in section \ref{sec5.2}.  In
section \ref{sec7}, we consider the possibility of alternate
experimental protocols that yield better conditioned inverse problems.
We draw heavily here on results from the text~\cite{BEGM}, which
contains, among other things, complete proofs of the main theorems
used (as well as more detailed numerical experiments).

\section{The Tangent Bundle to $\TA_{\ba}$} \label{sec2}
In this section we let $\bF$ denote an image with small support.
The question of transversality of the intersection at
$\bF\in\TA_{\ba}\cap B$ concerns the relationship between the fiber
of the tangent bundle to $\TA_{\ba}$ at $\bF$
and a linear approximation to the set $B$. For any $\bF\in\TA_{\ba}$
we let $T_{\bF}\TA_{\ba}$ denote the fiber of the tangent bundle to
$\TA_{\ba}$ at $\bF$.  The fiber of normal bundle at $\bF,$
$N_{\bF}\TA_{\ba},$ is the affine subspace through $\bF$ orthogonal to
$T_{\bF}\TA_{\ba}.$ We have defined these fibers as affine subspaces
of the ambient space $\bbR^J,$ and let $T^0_{\bF}\TA_{\ba},
N^0_{\bF}\TA_{\ba},$ denote the linear subspaces of $\bbR^J$ so that
\begin{equation}
  T_{\bF}\TA_{\ba}= \bF+T^0_{\bF}\TA_{\ba}\quad \text{ and }\quad
  N_{\bF}\TA_{\ba}=\bF+ N^0_{\bF}\TA_{\ba}.
\end{equation}

Recall that for $B=B_S,$ where $S$ is an estimate for the support of $\bF$,
$B_S$ is a linear subspace and the intersection is transversal if and
only if $T_{\bF}\TA_{\ba}\cap B_S=\{\bF\}$. 

\begin{remark}
If $B=B_+$ (see \eqref{Bplus}) then the
intersection lies on $\pa B_+,$ which is not a smooth submanifold of
$\bbR^J,$ but rather a stratified space. This renders the concept of
transversality more subtle to define. As the $\pa B_+$ is ``piecewise
linear'' in that it is locally a union of orthants in linear spaces of
various dimensions, it again makes sense to say that the intersection
is transversal provided that $T_{\bF}\TA_{\ba}\cap\pa
B_+=\{\bF\}.$ Since it is conceptually much simpler (and more general
in its applicability), most of our
discussion of transversality uses a support constraint as auxiliary
information. In Section~\ref{sec5.2} we briefly discuss the transversality of
the intersection with $\pa B_+.$
\end{remark}

\subsection{The Tangent Bundle in the DFT Representation}
The key to analyzing these intersections is to have an
explicit, readily computable description of the fibers of the tangent
bundle to $\TA_{\ba}.$ In this section we give two such descriptions.
The DFT \eqref{dft}, which we denote by $\cF$,
maps the torus $\TA_{\ba}$ onto a torus in $\bbC^J$
defined by
\begin{equation}
\begin{split}
  \hTA_{\ba}&:= \cF \TA_{\ba} =
  \{\hbf:\: |\hat f_{\bj}|=\ba_{\bj}, \; \text{ for }\bj\in J\}\\
&=\{(e^{i\theta_{\bj}}a_{\bj}:\:\bj\in J)\text{ for all }:\: \btheta\in \bbR^J\}.
\end{split}
\end{equation}
Taking $\theta_{\bj}$-derivatives
gives a very simple description of the tangent
bundle: for $\hbf\in\hTA_{\ba}$,
\begin{equation}
  T^0_{\hbf}\hTA_{\ba}=\Span_{\bbR}\{i\frac{\hat f_{\bj}}{|\hat f_{\bj}|}\be^{\bj}:\bj\in J\},
\end{equation}
where the standard basis vector
$\be^{\bj}\in\bbR^J$ has a $1$ in the $\bj$th location and is
otherwise zero. It is the ``real-span'' because 
$\hTA_{\ba}$ is a real submanifold of $\bbC^J.$

As the images, $\bF,$ we consider are real, this is reflected in a
symmetry of $\hbf$: for each index $\bj\in J$ there is a conjugate
index $\bj':=2(N-1)\bone-\bj,$ where $\bone=(1,\dots,1),$ for which
\begin{equation}
  \hat f_{\bj'}=\overline{\hat f_{\bj}}~;
\end{equation}
note that $(\bj')'=\bj.$  In fact, the fiber of the tangent bundle is the span of a smaller
set of vectors:
\begin{equation}
  T^0_{\hbf}\hTA_{\ba}=\Span_{\bbR}\left\{i\left[\frac{\hat f_{\bj}}{|\hat f_{\bj}|}\be^{\bj}-
\frac{\overline{\hat f_{\bj}}}{|\hat f_{\bj}|}\be^{\bj'}\right]:\bj\in J\right\}.
\end{equation}
The fiber of the normal bundle has a similar
description:
\begin{equation}
  N^0_{\hbf}\hTA_{\ba}=\Span_{\bbR}\left\{\left[\frac{\hat f_{\bj}}{|\hat f_{\bj}|}\be^{\bj}+
\frac{\overline{\hat f_{\bj}}}{|\hat f_{\bj}|}\be^{\bj'}\right]:\bj\in J\right\}.
\end{equation}

Up to a scale factor, the DFT is a unitary map, and therefore
$T_{\bF}\TA_{\ba}=\cF^{-1}[T_{\hbf}\hTA_{\ba}],$ though this is not a
very explicit, or useful description.

\subsection{The Tangent Bundle in the Image Representation}
We now give a second description,
in the image domain, of bases for
the tangent and normal bundles to $\TA_{\ba},$ whose elements share
many properties with that of the image itself. Recall that, for
$\bv\in J,$ the translate, $\bF^{(\bv)},$ of $\bF$ by $\bv$ is defined
in~\eqref{eqn2}.  Introduce the following images, which are the difference and sum of
translates of the image by $\bv$ and $-\bv,$
\begin{equation}
  \btau^{\bv}:= \bF^{(\bv)}-\bF^{(-\bv)}~,
  \qquad
    \bnu^{\bv}:= \bF^{(\bv)}+\bF^{(-\bv)}~.
\end{equation}
\begin{theorem}\label{thm1}
  Suppose that $\ba\in\bbR^J$ is the DFT magnitude data of a real
  image, $\bF.$ For any point $\bF\in \TA_{\ba}$, we have that
  \begin{equation}
      T^0_{\bF}\TA_{\ba}=\Span_{\bbR}\{\btau^{\bv}:\, \bv\in J\}~,\qquad
N^0_{\bF}\TA_{\ba}=\Span_{\bbR}\{\bnu^{\bv}:\, \bv\in J\}.
  \end{equation}
\end{theorem}
\begin{remark} The images $\{\bF^{(\bv)}\}$ are of course just  trivial
  associates of $\bF$, whose existence makes the solution of the phase
  retrieval problem non-unique. In fact the distances between the
  trivial associates are fairly large; an effective algorithm
  defined by a map with strong contraction properties would not have
  problems on this account. The theorem describes a far more insidious
  effect of the existence of trivial associates,
  as explained in the next paragraph:
  {\em it often renders the
    intersections of $\TA_{\ba}$ and $B_S$ non-transversal.}
  As we shall see in
  the next section, this adversely affects the continuity properties
  of the inverse map, which is entirely algorithm-independent. In
  Section~\ref{sec4} we see that it also vastly diminishes the contraction
  properties of the maps used to define phase retrieval algorithms,
  which inevitably leads to stagnation and even poorer reconstructions
  than would be expected from the results of Section~\ref{sec3}.
\end{remark}

Given an image, $\bF$, there is a subset $J_t$ of $J$ so that
$\{\btau^{\bv}:\bv\in J_t\},$ is a basis for the vector space
$T^0_{\bF}\TA_{\ba}.$ If $S$ is a realistic estimate for the support
of $\bF$, then there is usually a non-empty subset $J_{it}\subset J_t$
such that the tangent vectors $\{\btau^{\bv}:\:\bv\in J_{it}\}$ also
have support in $S$. In this case
\begin{equation}
  T^0_{\bF}\TA_{\ba}\cap B_S \;\supset\; \Span_{\bbR}\{\btau^{\bv}:\:\bv\in J_{it}\},
\end{equation}
which implies that the intersection at $\bF$ is \emph{not}
transversal.

For a subset $W\subset J$ the $p$-pixel neighborhood, $W_p,$ of $W$ is
defined to be
\begin{equation}\label{eqn2.10.002}
  W_p\overset{d}{=}\{\bj\in J:\:
\exists \bk\in W\text{ with } \|\bk-\bj\|_{\infty}\leq p\}.
\end{equation}
Once again this distance should be understood in the $J$-periodic sense.
A simple combinatorial argument shows that
\begin{equation}\label{eqn2.11.002}
  \dim T_{\bF}\TA_{\ba}\cap B_{S_{\bF,p}}\geq 2p(p+1),
\end{equation}
showing that a looser support constraint leads to a greater failure of
transversality.

\subsection{The Convolution Property of $T_{\bF}\TA_{\ba}$}\label{sec2.3.004}
In this section we examine a surprising property of the tangent bundle
to magnitude torus $\TA_{\ba}$ defined by a image $\bh$ that is a
convolution of two other images, that is
\begin{equation}
  \bh=\bF\ast \bg.
\end{equation}
In this section it is important to recall that we regard images in
$\bbR^{J}$ as periodic, i.e. as elements of $\bbR^{\bbN^d}$ with
indices in $J\subset\bbN^d$ representing a single period.

This discussion requires some additional notation. For
$\bF\in\bbR^J,$ we let
\begin{equation}
  \ba_{\bF}\overset{d}{=}(|\hf_{\bk}|:\:\bk\in J),
\end{equation}
so that $\TA_{\ba_{\bF}}$ is the magnitude torus defined by $\bF,$ and we
let
\begin{equation}
  \btau^{(\bv)}_{\bF}\overset{d}{=}\bF^{(\bv)}-\bF^{(-\bv)}.
\end{equation}
Suppose that $J_t\subset J$ is chosen so that $\{\btau^{\bv}_{\bF}:\: \bv\in
J_t\}$ is a basis for $T_{\bF}\TA_{\ba_{\bF}}.$ If $\balpha\in\bbR^{J_t},$
then we let
\begin{equation}
  \btau^{\balpha}_{\bF}\overset{d}{=}\sum_{\bv\in
    J_t}\alpha_{\bv}\btau^{(\bv)}_{\bF}\in T_{\bF}\TA_{\ba_{\bF}}.
\end{equation}

Recall that (periodic) discrete convolution is defined by
\begin{equation}
  [\bF\ast\bg]_{\bj}=\sum_{\bk\in J} f_{\bj-\bk}g_{\bk}.
\end{equation}
The support of a convolution satisfies
\begin{equation}\label{eqn2.15.002}
  S_{\bF\ast\bg}\subset S_{\bF}+S_{\bg},
\end{equation}
where we recall that if $X,Y\subset J,$ then $X+Y\overset{d}{=}\{\bj+\bk:\:\bj\in
X,\, \bk\in Y\}\mod J.$ The DFT coefficients of a convolution satisfy
\begin{equation}\label{2.16.002}
  \widehat{\bF\ast\bg}_{\bj}=\hf_{\bj}\hg_{\bj}.
\end{equation}

It is an elementary computation to show that
\begin{equation}\label{eqn2.15.01}
  \btau^{(\bv)}_{\bF\ast\bg}=\btau^{(\bv)}_{\bF}\ast\bg=\bF\ast\btau^{(\bv)}_{\bg}.
\end{equation}
For generic images, $\bF$ and $\bg,$ a single index set $J_t$ can be
used so that $\{\btau^{\bv}_{\bF}:\: \bv\in J_t\}$ is a basis for
$T_{\bF}\TA_{\ba_{\bF}},$ and $\{\btau^{\bv}_{\bg}:\: \bv\in J_t\}$ is a
basis for $T_{\bg}\TA_{\ba_{\bg}}.$ Letting $\balpha\in \bbR^{J_t}$, it
follows from~\eqref{eqn2.15.01} that
\begin{equation}\label{eqn2.20.003}
  \btau^{\balpha}_{\bF\ast \bg}=\btau^{\balpha}_{\bF}\ast\bg=\bF\ast\btau^{\balpha}_{\bg}.
\end{equation}
More succinctly we can write:
\begin{equation}
  T_{\bF\ast\bg}\TA_{\ba_{\bF\ast\bg}}=T_{\bF}\TA_{\ba_{\bF}}\ast\bg=
  \bF\ast T_{\bg}\TA_{\ba_{\bg}}.
\end{equation}
If we let 
\begin{equation}
  S_{\btau^{\balpha}_{\bF}}=\{\bj\in J:\:\btau^{\balpha}_{\bF,\bj}\neq 0\},
\end{equation}
then~\eqref{eqn2.15.002} and~\eqref{eqn2.20.003} imply that
\begin{equation}
  S_{\btau^{\balpha}_{\bF\ast \bg}}\subset S_{\btau^{\balpha}_{\bF}}+S_{\bg}.
\end{equation}

Suppose that $\bg$ is an image for which there exists an $\balpha\in J_t$ so that
\begin{equation}\label{eqn2.21.002}
  S_{\btau_{\bg}^{\balpha}}\subset S_{\bg}.
\end{equation}
The vector field $\bF\ast\btau_{\bg}^{\balpha}\in
T_{\bF\ast\bg}\TA_{\ba_{\bF\ast\bg}}$ then satisfies
\begin{equation}\label{eqn2.22.002}
  S_{\bF\ast\btau_{\bg}^{\balpha}}\subset S_{\bg}+S_{\bF}.
\end{equation}
On its face the condition in~\eqref{eqn2.21.002} seems very unlikely to hold for
any $\balpha,$ as it is equivalent to the system of
linear equations for $\balpha\in J_t$:
\begin{equation}\label{eqn2.23.002}
  \btau_{\bg,\bj}^{\balpha}=0\text{ for }\bj\in S^c_{\bg}.
\end{equation}
Since $\dim T_{\bg}\TA_{\ba_{\bg}}=|J_t|\approx|J|/2,$ and $|S_{\bg}^c|= |J|-|S_g|\geq
\frac{3|J|}{4},$ for an image with small support, the equations
in~\eqref{eqn2.23.002} appear to be overdetermined.

It turns out that if $\bg$ is inversion symmetric, that is
\begin{equation}
  \bg_{\bj}=\bg_{-\bj},
\end{equation}
then the basis vectors for $T_{\bg}\TA_{\ba_{\bg}}$ satisfy the
equations
\begin{equation}
  \btau^{(\bv)}_{\bg,\bj}=- \btau^{(\bv)}_{\bg,-\bj}\text{ for all
  }\bv\in J_t,\,\bj\in J.
\end{equation}
For the case of an inversion symmetric image, half of the equations
in~\eqref{eqn2.23.002} imply the other half, from which the following
theorem follows easily.
\begin{theorem}\label{thm2}
  Let $\bg\in\bbR^J$ be inversion symmetric, then the solution space
  to the equations in~\eqref{eqn2.23.002} has dimension at least
  \begin{equation}
   |J_t|- \frac{|J|-|S_{\bg}|}{2}\approx \frac{|S_{\bg}|}{2}.
  \end{equation}
\end{theorem}
\begin{remark} This theorem also holds for images that
  are anti-symmetric, i.e. $\bg_{\bj}=-\bg_{-\bj},$ and images that are inversion
  symmetric (or anti-symmetric) with respect to any point in $J.$ 
\end{remark}

An inversion symmetric image has real DFT coefficients, and therefore
that phase retrieval problem for such an image reduces to the much
easier sign retrieval problem. But suppose that $\bF$ is an arbitrary
image with small  support, and $\bg$ is inversion symmetric
so that $\bF\ast\bg$ also has small support. The theorem along
with~\eqref{eqn2.22.002} imply that
\begin{equation}
  \dim T_{\bF\ast\bg}\TA_{\ba_{\bF\ast\bg}}\cap
  B_{S_{\bF}+S_{\bg}}\geq \frac{|S_{\bg}|}{2}.
\end{equation}
For non-negative images $S_{\bF\ast\bg}=S_{\bF}+S_{\bg};$ in any case, the sum
is a reasonable estimate for $S_{\bF\ast\bg}.$ The failure of
transversality, with $S=S_{\bF\ast\bg}$  is \emph{inherited} by
$\bF\ast\bg$ from $\bg,$ even though $\bF\ast\bg$ has no obvious symmetries.

\subsection{Examples of the Failure of Transversality for Convolutions}\label{sec2.4.004}
The analysis in the previous section shows that for images that are
convolutions with inversion symmetric images, the failure of
transversality occurs even if we use the exact support of the image to
define the support constraint. In order to diminish the effects of
noise, it is a very common practice to multiply measured data by a
smooth cut-off function, such as a Gaussian. Since our measurements are
in the DFT domain, equation~\eqref{2.16.002} shows that this is
equivalent to convolving the unknown image with a Gaussian.  Our
analysis suggests that this makes the problem of recovering the phase
much more difficult. In this section we present the results of
numerical experiments that demonstrate this phenomenon.

\begin{figure}[h]
\begin{centering}
 \centerline{\includegraphics[width=13cm]{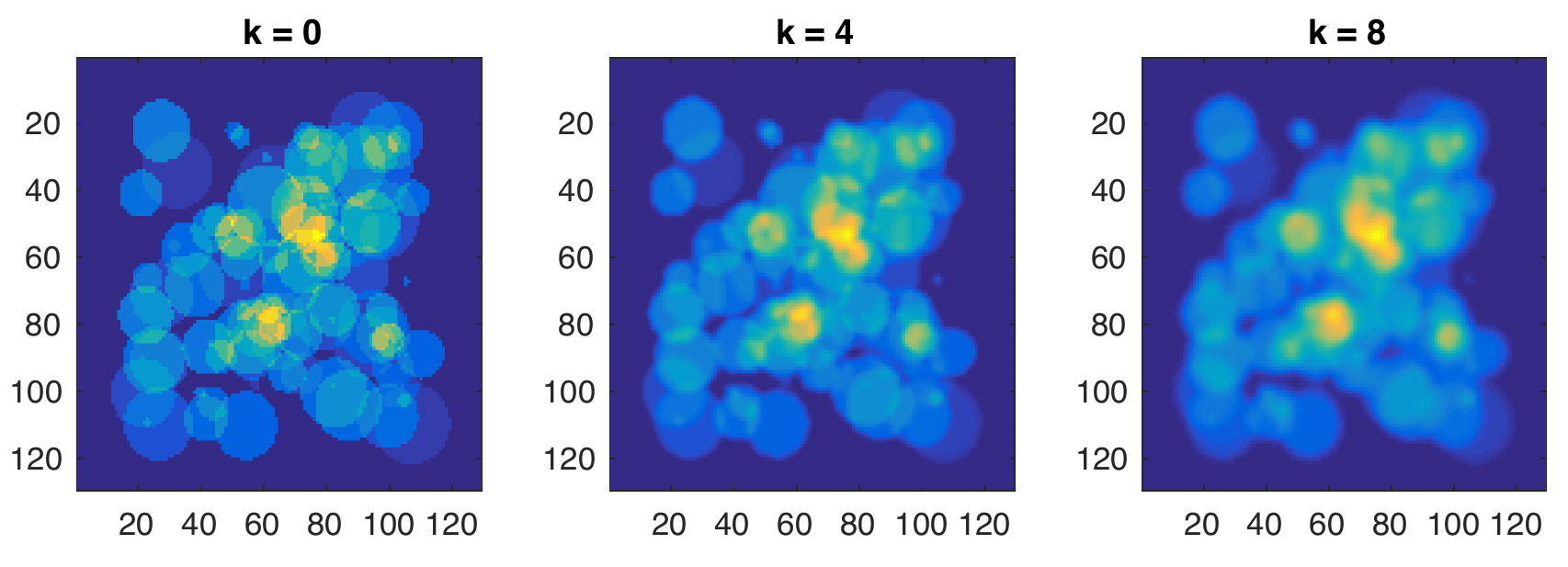}}
        \caption{Images similar to those used in numerical experiments
          below of various smoothness levels: $k=0, 4, 8.$  The parameter $k$ scales
             the width of the Gaussian in \eqref{smoother} used to smooth the 
             image. These are $256\times 256$ images ($N=128$), with the object
             contained in a $128\times 128$ square.}\label{fig3}
\end{centering}
\end{figure}

For our numerical experiments, we use images defined by a sum
of radial functions,
\begin{equation}
  \rho(\bx)=\sum_{i=1}^IR_i(\bx)~.
\end{equation}
In the simplest case, each function $R_i$ is a scaled characteristic function
of a disc,
\begin{equation}
  R^0_i(\bx)=\alpha_i\chi_{[0,r_i]}(\|\bx-\bcc_i\|)~,
\end{equation}
where the intensities $\alpha_i>0$, radii $r_i>0$, and centers $\bcc_i$
are set randomly; see left side of Fig.~\ref{fig3}(a).
In this case the function $\rho$ is piecewise constant; in the imaging
literature one would say that $\rho$ represents a \emph{hard} object.
The discrete image $\bF$ is then generated from point samples of $\rho$
on a regular grid.
We also generate smoother images by
(discrete) convolution of this $\bF$ with the discretely sampled Gaussian,
\begin{equation}
  G_k(\bx) = c_k e^{-\frac{16N^2\|\bx\|^2}{(k+1)^2}}~, \qquad k>0,
\label{smoother}
\end{equation}
where $k$ controls the smoothness,
and where $c_k$ is chosen to make $\|G_k\|_1=1$. 
(The unsmoothed case we denote by $k=0$.)
In Figure~\ref{fig3} we show examples of such images with
smoothness levels $0$ (unsmoothed), $4$, and $8$. 

Since the smoothing in \eqref{smoother} results in full support,
we instead set a threshold $\epsilon\approx\epsmach$
appropriate for finite-precision arithmetic,
and define the support to be
\begin{equation}
  S_{\bF}\overset{d}{=}\{\bj\in J:\:|f_{\bj}|\geq \epsilon\}~.
\end{equation}
The $p$-pixel neighborhoods, $S_{\bF,p}$ of $S_{\bF}$ are defined as
in~\eqref{eqn2.10.002}.  To see how $\dim T_{\bF}\TA_{\ba}\cap
B_{S_{\bF,p}}$ depends on the degree of smoothing, $k$, and the size
of the support padding, $p$, we generate 5 random samples,
$\{\bF_i^{k}:\:i=1,\dots,5\},$ for each smoothness level
$k=0,1,2,3,4.$
Since the intersection dimension computation requires a dense SVD
of a matrix of size $\bigO(N^2\times N^2)$, the study is
limited to small images; we choose the image size so that the
double oversampled image is $64\times 64$. For these images
$  \dim T_{\bF_i^k}\TA_{\ba}=2046\text{ if }k=0$  and $1984\text{ if }k>0.$  
The dimension decreases when $k>0$ as the symmetries of $G_k$
forces certain DFT coefficients to vanish.

\begin{figure}[h]
    \centering
    \begin{subfigure}[H]{.4\textwidth}
        \centering
        \includegraphics[height= 4cm]{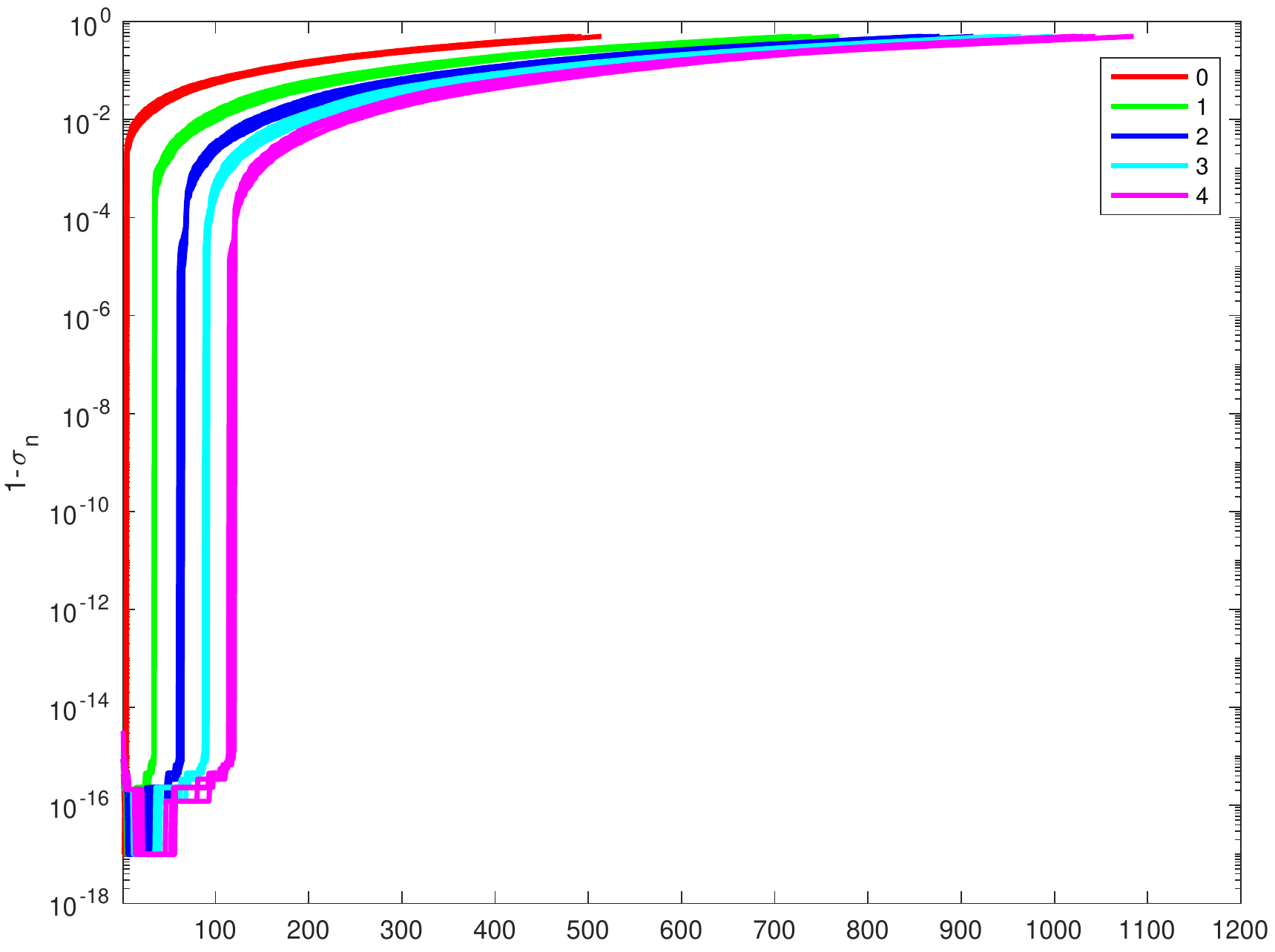}
        \caption{1-pixel neighborhoods}
    \end{subfigure}\qquad
\begin{subfigure}[H]{.4\textwidth}
        \centering
        \includegraphics[height= 4cm]{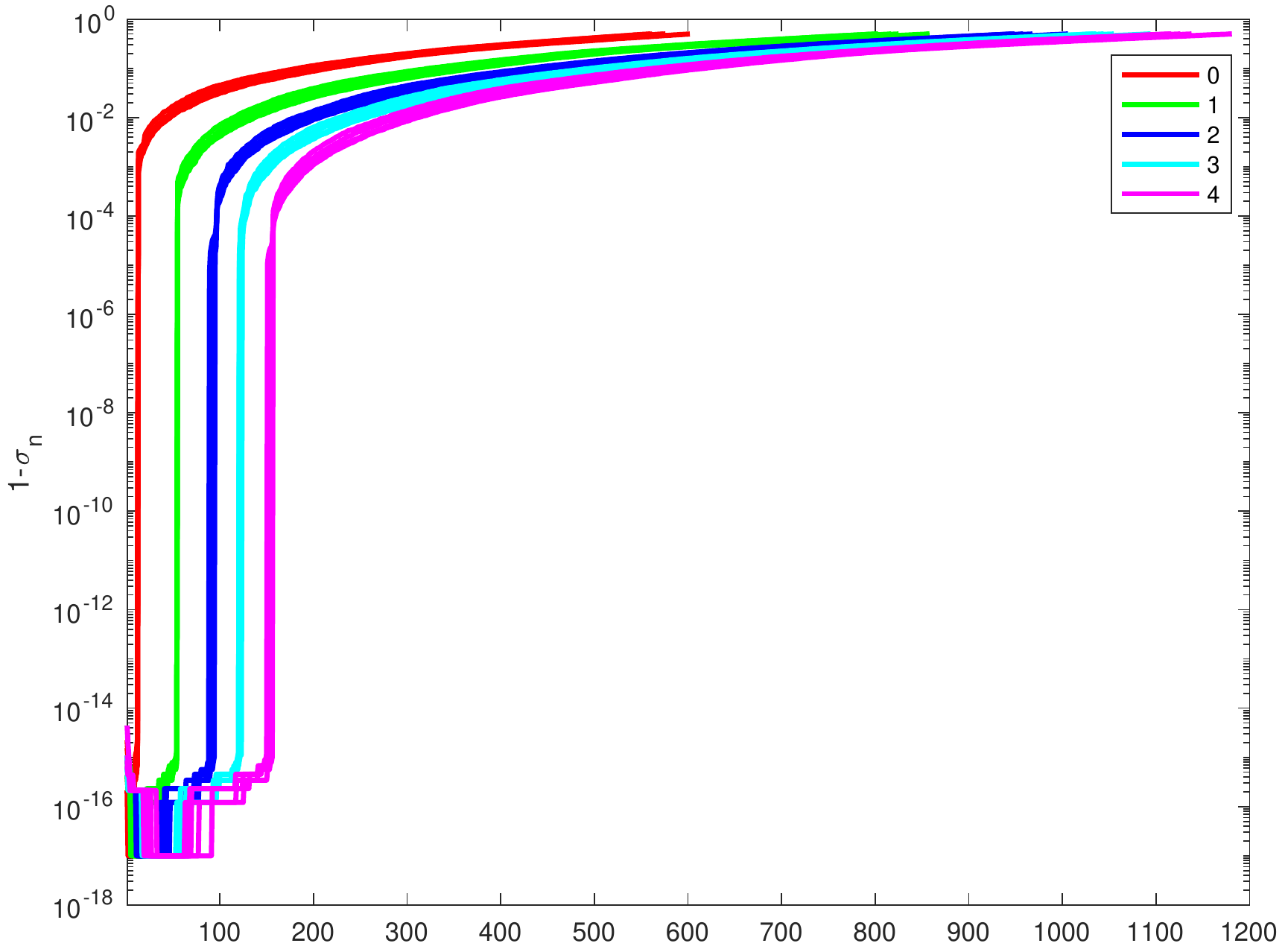}
        \caption{2-pixel neighborhoods}
    \end{subfigure}
\begin{subfigure}[H]{.4\textwidth}
        \centering
        \includegraphics[height= 4cm]{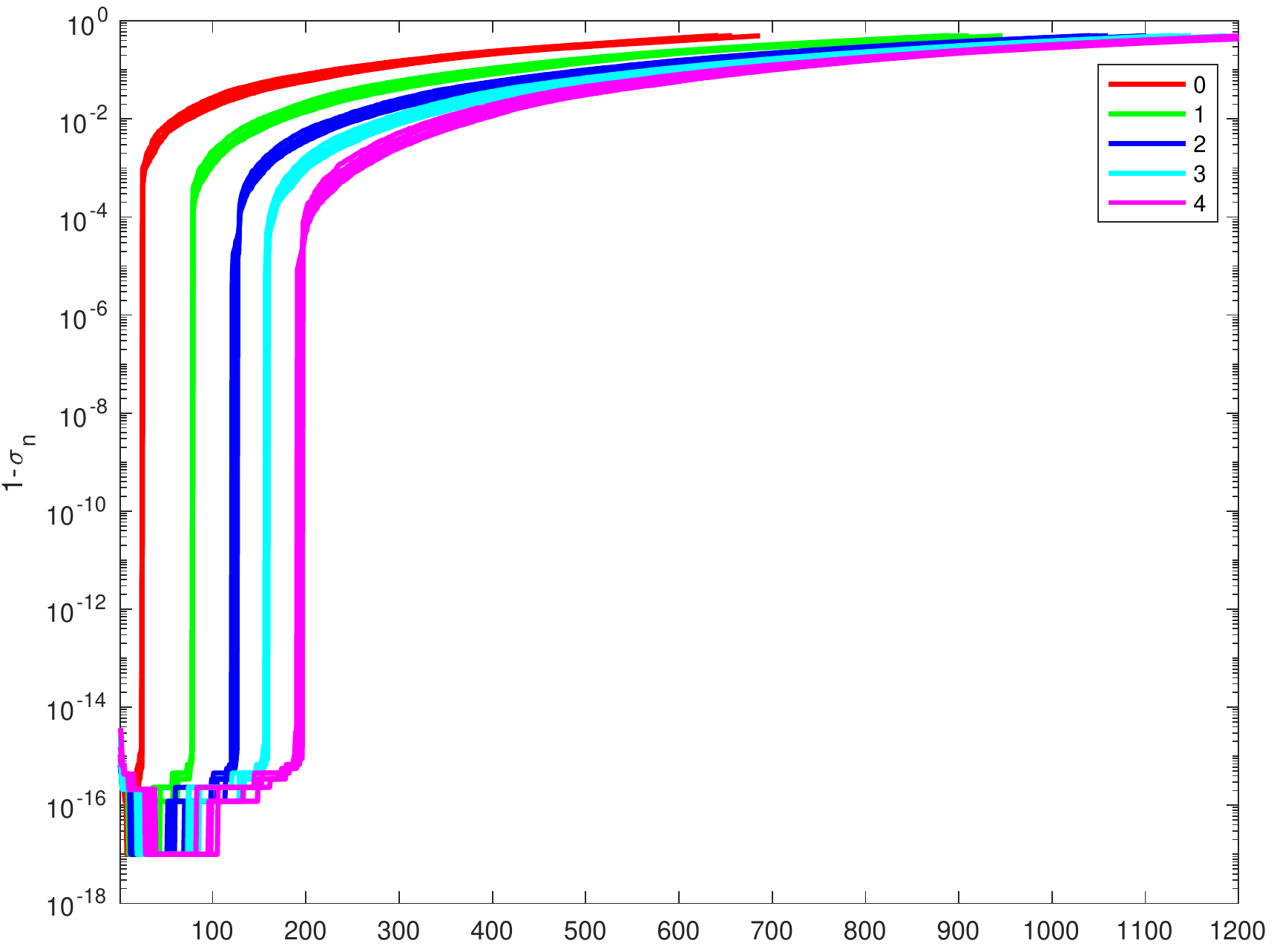}
        \caption{3-pixel neighborhoods}
    \end{subfigure}\qquad
\begin{subfigure}[H]{.4\textwidth}
        \centering
        \includegraphics[height= 4cm]{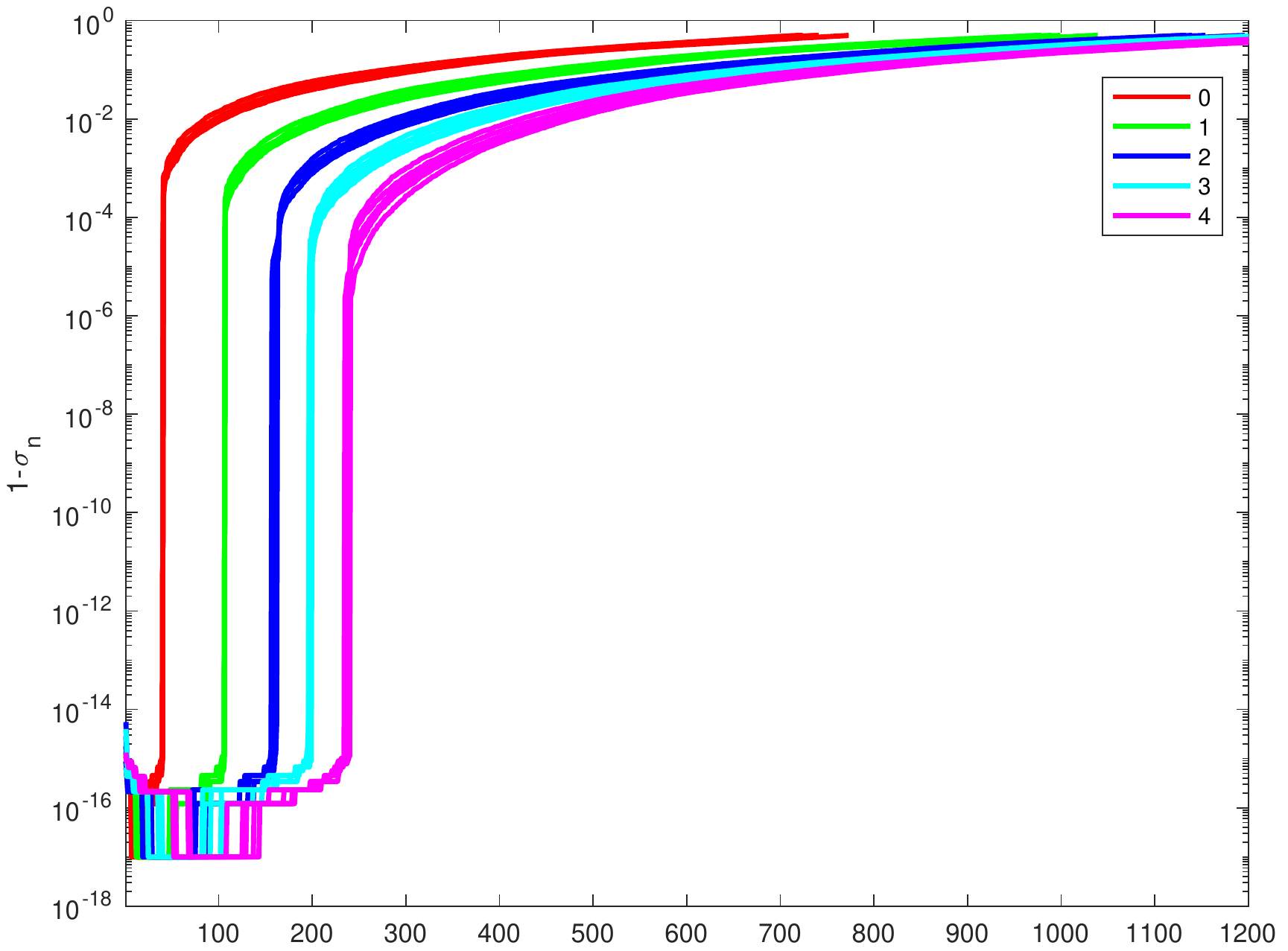}
        \caption{4-pixel neighborhoods}
    \end{subfigure}
    \caption{Singular values of $V^tU$ from images with varying degrees of smoothness
      $k=0,1,2,3,4$ (corresponding to the colors: red, green, blue,
      cyan, magenta). Plots produced by 5 random $64\times
      64$-examples are shown for $p=1,2,3,4.$ The plots show
      $\log_{10}(1-\sigma_n),$ where $\{\sigma_n\}$ are the singular
      values of $H.$ }\label{fig4}
\end{figure}

For each sample image, we numerically compute $U$, an orthonormal
basis for $T^0_{\bF_i^{k}}\TA_{\ba}$, and $V$, an orthonormal basis
for $B_{S_{\bF_i^k,p}}$, and then compute the SVD of $H=V^tU$.  In
exact arithmetic, $\dim T_{\bF_i^{k}}\TA_{\ba}\cap B_{S_{\bF_i^k,p}}$
is the number of singular value, $\{\sigma_n(H)\},$ equal to 1.  Since we
work in finite-precision arithmetic, such singular values are only
approximately 1.  Figure~\ref{fig4} contains plots of
$\log_{10}(1-\sigma_n).$ Note that there is a sharp transition from
singular values within $10^{-15}$ of 1 to smaller ones, which is
remarkably consistent across the samples.  As the support of the
Gaussian, $G_k,$ (at machine precision) grows with $k,$ this is
essentially as predicted by Theorem~\ref{thm2}. 

We summarize these dimension measurements in Table~\ref{tab1}.  The
dimensions shown correspond to the number of singular values greater
than $1-10^{-15}$. The $k=0$ row is precisely $2p(p+1)$ as predicted
in~\eqref{eqn2.11.002}.  The table also has a $p=0$ column where we
have used the exact support, $S_{\bF_i^{k}},$ to define the support
constraint.  As predicted by Theorem~\ref{thm2}, for $k>0$ the
dimensions of these intersections are non-zero.  As $k$ increases, the
support of $G_k$ grows, and this produces a larger and larger
intersection between $T_{\bF_i^{k}}\TA_{\ba}$ and $B_{S_{\bF_i^k}}.$
As $p$ increases the dimension of these intersections grow beyond
$|S_{G_k}|/2:$ the larger $p$ is, the more ways there are to obtain
tangent vectors of the form $\btau_{G_k}^{\balpha}\ast\bF^0$ so that
$S_{\btau_{G_k}^{\balpha}}+S_{\bF^k_i}\subset S_{\bF^k_i,p}.$ In
examples below we show that even a very low dimensional failure of
transversality can lead to stagnation in standard reconstruction
algorithms.

\begin{table}[h]
\begin{center}
    \begin{tabular}{| c || c | c | c |c|c|}
    \hline
    smthns\textbackslash supp &$p=0$& $p=1$ & $p=2$ & $p=3$& $p=4$ \\ \hline\hline
    $k=0$&0 &4 &12 &24 &40 \\ \hline
$k=1$&18 &34 &54 &78 &106 \\ \hline
 $k=2$&38 &64 &92 &124 &160 \\ \hline
 $k=3$&61 &88 &120 &156 &196 \\ \hline
 $k=4$&85 &119 &155 &195 &239 \\ \hline
    \end{tabular}
\end{center}
\caption{Table showing the typical dimensions of
  $T_{\bF^k_i}{\TA_{\ba}}\cap {B_{S_{\bF^k_i,p}}}$ for $p=1,2,3,4$ and
  varying degrees of smoothness. }\label{tab1}
\end{table}

\begin{figure}[h]
    \centering
    \begin{subfigure}[H]{.4\textwidth}
        \centering
        \includegraphics[height= 4cm]{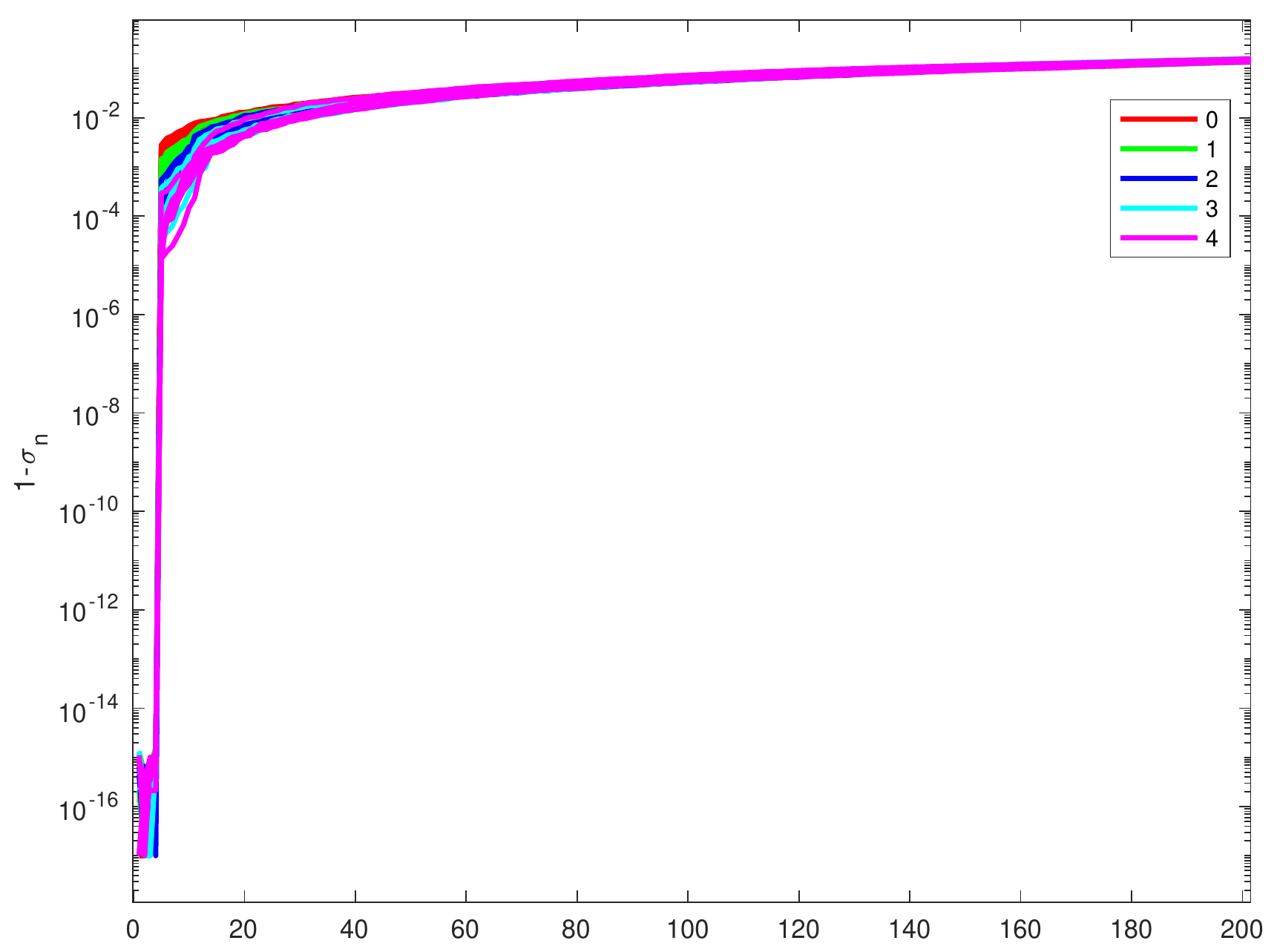}
        \caption{1-pixel neighborhoods}
    \end{subfigure}\qquad
\begin{subfigure}[H]{.4\textwidth}
        \centering
        \includegraphics[height= 4cm]{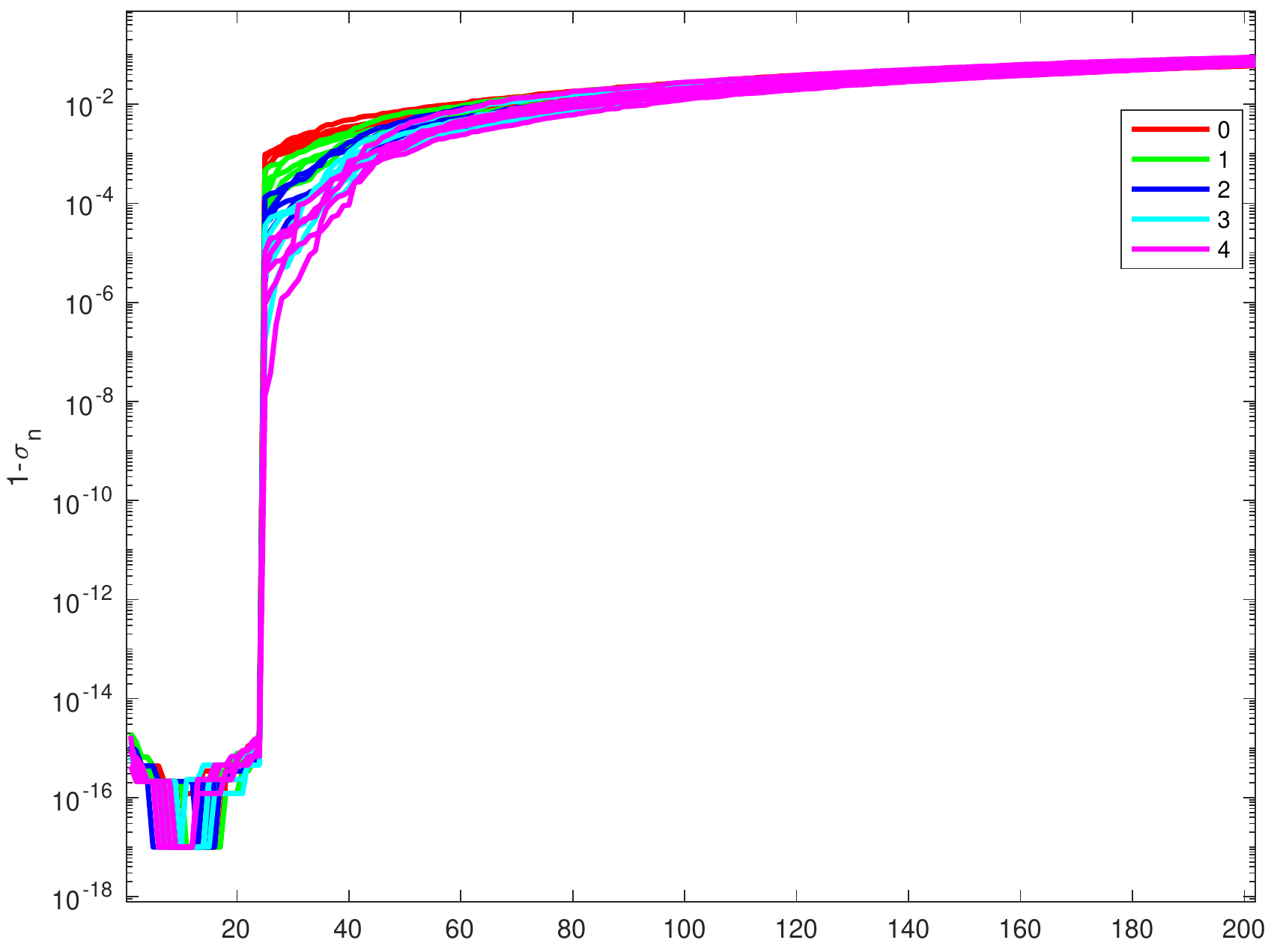}
        \caption{3-pixel neighborhoods}
    \end{subfigure}
    \caption{The first 200 singular values of $H=V^tU$ from images
      with varying degrees of smoothness $k=0,1,2,3,4$ (corresponding
      to the colors: red, green, blue, cyan, magenta). Plots produced
      by 5 random $64\times 64$-examples are shown for $p=1,3.$ The
      plots show $(\log_{10}(1-\sigma_n):\: n=1,\dots,200),$ where
      $\{\sigma_n\}$ are the singular values of $H.$ }\label{fig4.1}
\end{figure}

This experiment can be repeated using images defined as samples, $\bF^k,$ of
functions of the form
\begin{equation}\label{eqn2.35.004}
  \rho(\bx)=\sum_{i=1}^IR^k_i(\bx),
\end{equation}
where
\begin{equation}\label{eqn2.36.004}
  R^k_i(\bx)=\alpha_i\chi_{[0,r_i]}(\|\bx-\bcc_i\|)\|\bx-\bcc_i\|^k.
\end{equation}
These functions increase in smoothness with $k,$ but are not defined
as convolutions.  For each $k\in\{0,1,2,3,4\}$ we produce 5 random
choices of $\rho.$ In this case we find that the dimensions of the
intersections are independent of $k,$ satisfying $ \dim
T_{\bF^k_i}\TA_{\ba_{\bF^k_i}}\cap B_{S_{\bF^k_i,p}}=2p(p+1);$ in this
case $ \dim T_{\bF^k_i}\TA_{\ba_{\bF^k_i}}=2046$ for all $k.$ In
Figure~\ref{fig4.1} we show plots of $(\log(1-\sigma_n):\:
n=1,\dots,200),$ for each choice $\bF^k_i,\, i=1,\dots,5,$ with
$p=1,3.$ From these plots we see that $\dim
T_{\bF^k_i}\TA_{\ba_{\bF^k_i}}\cap B_{S_{\bF^k_i,p}}$ does not depend
on $k,$ or the choice of example, and that the number of directions in
which $T_{\bF^k_i}\TA_{\ba_{\bF^k_i}}$ meets $B_{S_{\bF^k_i,p}}$ at
very small angles does increase with $k.$ These small angles also
dramatically stall the convergence of standard algorithms for finding
these intersections.

\section{Transversality, Well-Posedness and Microlocal Non-uniqueness}\label{sec3}
We turn now to an analysis of the effects of a non-transversal
intersection on the computational difficulty of the phase retrieval
problem. While the results in this section are algorithm independent,
they have direct implications about the loss of solution accuracy
given finite precision data and computations.
The results in this section are related to and, in part,
inspired by those in~\cite{cahill2016}
and~\cite{alaifari2017}.

\subsection{Transversality and Well-Posedness}

Recall from the introduction
that $\cM:\bbR^J\to \bbR_+^J,$ denotes the measurement map
$\cM(\bF)=(|\hat f_{\bj}|:\:\bj\in J)$.
First note that this map is Lipschitz continuous in the $2$-norm,
\begin{equation}
  \|\cM(\bF)-\cM(\bg)\|_2 \;\le\; C_{\cM}\|\bF-\bg\|_2~,
  \label{lip}
\end{equation}
with $C_{\cM} = \sqrt{|J|} = (2N)^{d/2}$,
which follows from the Plancherel theorem for the DFT \eqref{dft}
and the triangle inequality.
The inverse image of a point
$\ba\in\bbR_+^J$ is the magnitude torus $\TA_{\ba}$.
Suppose that
$\bF\in\TA_{\ba}$ has small support contained in the set $S$, then
the set $\cM^{-1}(\ba)\cap B_S$ is finite, and non-empty,
and the local inverse $\cM^{-1}_{\bF,S}$ is defined on a neighborhood
of $\ba=\cM(\bF)$ in the set of consistent data $\cM(B_S)$.

As is typical in the field of inverse problems, in order for the
problem to be locally well posed at $\bF$ it is necessary for this local
inverse to be a Lipschitz map. That is, there must be a neighborhood
$U\subset B_S$ of $\bzero$ and a constant $C>0$, such that for
$\delta\bF\in U$, we have the estimate
\begin{equation}\label{eqn19}
  C\|\delta\bF\|_2\leq \|\cM(\bF)-\cM(\bF+\delta\bF)\|_2,
\end{equation}
so that if $\ba=\cM(\bF),$ and $\ba+\delta\ba=\cM(\bF+\delta\bF),$ then
\begin{equation}
  \|\cM^{-1}_{\bF,S}(\ba)-\cM^{-1}_{\bF,S}(\ba+\delta\ba)\|_2\leq
  \frac{1}{C}\|\delta\ba\|_2.
\end{equation}

Now assume that $\TA_{\ba}\cap B_S$ is non-transversal at $\bF,$ and
let $\btau\in T_{\bF}\TA_{\ba}\cap B_S$ be a unit vector, then the
definition of the tangent bundle, and \eqref{lip}, imply that there
is a constant $c$, dependent on $\bF$, so that
\begin{equation}\label{eqn21}
  \|\cM(\bF)-\cM(\bF+t\btau)\|\leq ct^2~,
  \qquad \mbox{ for all sufficiently small real $t$~.}
\end{equation}
In this case the best general bound one can hope for is that 
\begin{equation}
  \|\cM^{-1}_{\bF,S}(\ba)-\cM^{-1}_{\bF,S}(\ba+\delta\ba)\|_2 \; \leq \;
  \frac{1}{C}\sqrt{\|\delta\ba\|_2}~.
\end{equation}
This implies that the local inverse is, at best, H\"older
continuous of order $\frac 12$, and therefore has an unbounded condition number.  As noted
in the introduction, it also implies that if the measurements have $d$
significant digits, then, generally, it will be impossible to
reconstruct an image with more than $\frac{d}{2}$ significant digits.

In~\cite{BEGM}, we prove the following result:
\begin{theorem}\label{transversality_thm}
  For an image $\bF\in B_S$, let $\TA_{\ba}$ denote the magnitude torus
  defined by $\ba=\cM(\bF)$. Suppose that $a_{\bj}\neq 0$ for all
  $\bj\in J$.  There are positive constants, $\eta$, $C$, so that, if
  $\delta\bF\in B_S$ and $\|\delta\bF\|_2<\eta$, then
  \begin{equation}
   C\|\delta\bF\|_2 \;\leq\; \|\cM(\bF)-\cM(\bF+\delta\bF)\|_2
  \end{equation}
  if and only if $T_{\bF}\TA_{\ba}\cap B_S=\{\bF\}$.
\end{theorem}

The theorem says that, if a support condition is the auxiliary
information that is available, and the DFT data is generic
(non-vanishing), then the phase retrieval problem can only be
well-conditioned near to $\bF\in \TA_{\ba}\cap B_S$ if this
intersection is transversal.  This statement is intrinsic to the phase
retrieval problem, i.e.\ is algorithm independent. The results in
section~\ref{sec4} indicate that, with a realistic support condition,
these intersections are very rarely transversal (see Table~\ref{tab1}
and Figure~\ref{fig4.1} above). As our numerical experiments below
show, this failure of transversality can also dramatically harm the
convergence properties of standard algorithms.
%

\subsection{$\epsilon$-Non-Uniqueness} \label{nonuniquesec}
As discussed in the introduction, the solution to the phase retrieval problem
with support condition is not always unique up to trivial
associates. From the discussion in the previous section we already know that
the conditioning of the phase retrieval problem depends subtly on the unknown
image, and the precise nature of the auxiliary information. In this section we
explore various ways in which this problem can fail to have a unique solution
to a given precision $\epsilon>0.$ Suppose that there are two images
$\bF_1,\bF_2,$ and a subset $S\subset J,$ adequate for generic uniqueness, such
that
\begin{enumerate}
\item The norms $\|\bF_1\|_{2}=\|\bF_2\|_2,$ but the minimum distance between trivial
  associates of $\bF_1$ and $\bF_2$ is much larger than $\epsilon\|\bF_1\|_{2}.$
\item The sets $\{\bj:\epsilon<|f_{i\bj}|\}\subset S,$ for $i=1,2.$
\item $\|\cM(\bF_1)-\cM(\bF_2)\|_{2}<\epsilon.$
\end{enumerate}
then we say that the solution to the phase retrieval problem defined by the
data $(\cM(\bF_1),S)$ is $\epsilon$-non-unique. In the remainder of this
section we describe two distinct mechanisms
that lead to $\epsilon$-non-uniqueness.

\subsubsection{Consequences of Genuine Non-Uniqueness}
The fact that a discrete image, with sufficiently small support, is
generically determined by the magnitude DFT data is a consequence of
the classical theorem that polynomials in two or more variables are
generically irreducible over the complex numbers. If $(f_{\bj}:\bj\in
J)$ is the image, then its $\bZ$-transform is
\begin{equation}
  \BF(\bZ)=\sum_{\bj\in J}f_{\bj}\bZ^{-\bj},
\end{equation}
where $\bZ^{-\bj}=z_1^{-j_1}\cdots z_d^{-j_d}.$ There is a minimal integer
vector $\bm$ so that $\bZ^{\bm}\BF(\bZ)$ is a polynomial.

Suppose that $\bF$ is an image whose $\bZ$-transform, $\BF(\bZ)$ is reducible,
in the sense that there are polynomials, $\BF_1,\BF_2$ in $\bZ$ such that 
\begin{equation}
  \BF(\bZ)=\bZ^{\bn}\BF_1(\bZ)\BF_2(\bZ),
\end{equation}
for some integer vector $\bn.$ If $\bF_1$ and $\bF_2$ are images with
$\bZ$-transforms $\BF_1,$ $\BF_2,$ (up to a factor of $\bZ^{\bm_i}$ for some
$\bm_i$) then, up to a translation, $\bF=\bF_1\ast\bF_2,$ where $\ast$ denotes
discrete convolution. If no trivial associate of either $\bF_1$ or $\bF_2$ is
inversion symmetric, then the image $\bF'=\bF_1\ast\check{\bF}_2$ is not
a trivial associate of $\bF$ and, typically, the minimum distance between the
trivial associates of $\bF$ and $\bF'$ is large. If $\bF_1$ and $\bF_2$ are
non-negative, then so are $\bF$ and $\bF',$ and the smallest rectangles
containing each image coincide.

Suppose that $\bF$ and $\bF'$ both have small support contained in a
set $S,$ which is small enough to generically imply uniqueness, up to
trivial associates, and let $0<\epsilon$ be chosen with $\epsilon\ll
\|\bF-\bF'\|_2.$ Because uniqueness is generic we can modify these two
images to obtain generic images $\bg$ and $\bg',$ so that the norms
satisfy the estimates $\|\bg-\bF\|_2<\epsilon/2$,
$\|\bg'-\bF'\|_2<\epsilon/2$, and $S_{\bg}=S_{\bF}$,
$S_{\bg'}=S_{\bF'}$. The data $(\cM(\bg),S)$ defines both a phase
retrieval problem with a unique solution, up to trivial associates,
and an $\epsilon$-non-unique problem.  That is, we can construct images
$\bg$ and $\bg',$ with support in $S,$ and nearly
identical magnitude DFT data:
$$\|\cM(\bg)-\cM(\bg')\|\leq \epsilon,$$
but satisfying
$$\|\bg-\bg'\|_2\geq \|\bF-\bF'\|_2-\epsilon>\!\!>0.$$

\begin{figure}[h]
\begin{centering}
    \begin{subfigure}[t]{.8\textwidth}
 \centerline{\includegraphics[width=11cm]{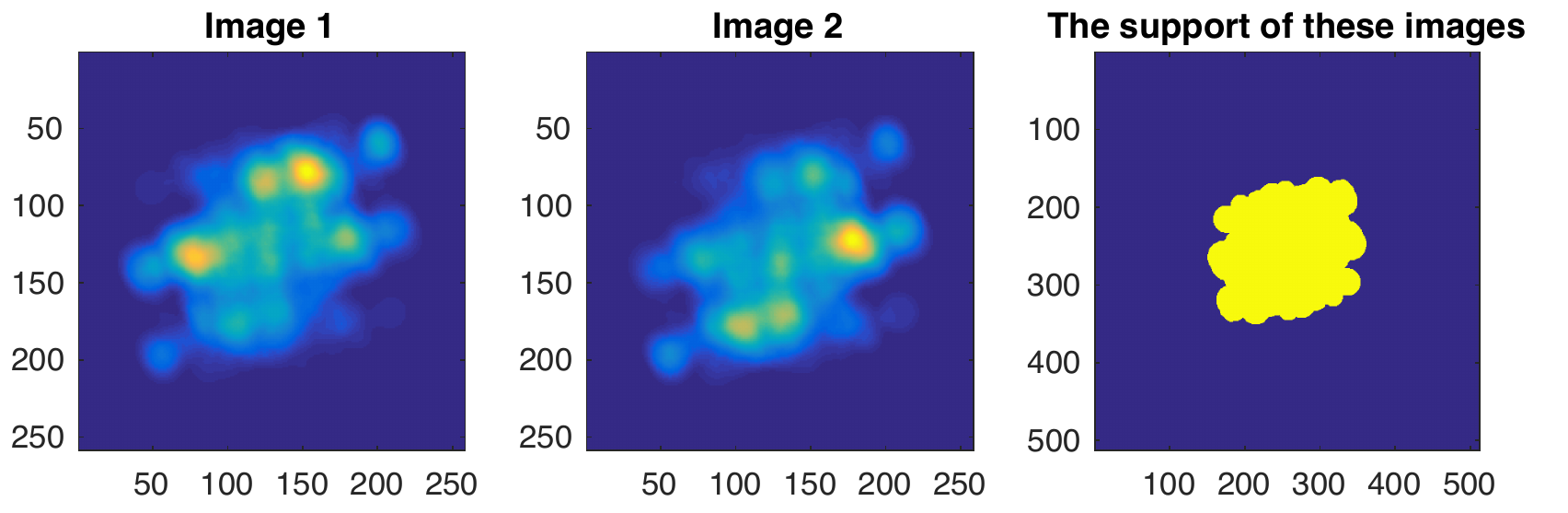}}
        \caption{A pair of images with essentially identical magnitude-DFT data
          that are not trivial associates, and their common support. For
          clarity the left and center images show the central $256\times 256$
portion of the original $512\times 512$ image.}
    \end{subfigure}\qquad
\begin{subfigure}[t]{.8\textwidth}
          \centerline{\includegraphics[width=11cm]{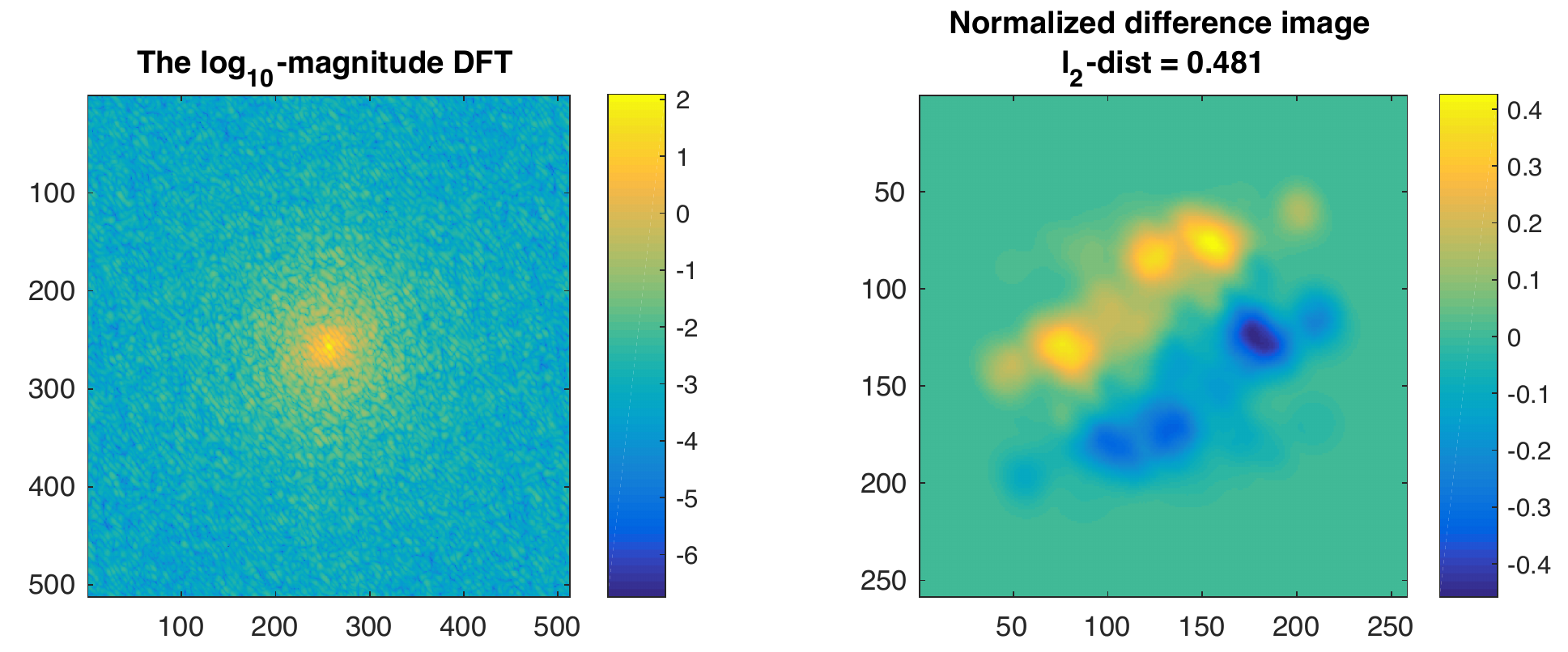}}
        \caption{$\log_{10}$-Magnitude-DFT data for the images above,  and the
          difference between the images themselves.
          The image on the right shows the central $256\times 256$ portion of the original $512\times 512$ image.}
    \end{subfigure}
           \caption{ An illustration of true non-uniqueness in the phase retrieval
             problem.}\label{fig2.1}
\end{centering}
\end{figure}

We conclude this section with an example of a pair of non-negative images,
$\bF,\bF'$, with exactly the same support and magnitude-DFT data such that
$\|\bF-\bF'\|_{2}\approx .48\|\bF\|_2.$ The minimum distance between trivial
associates is about $.18\|\bF\|_2,$ but the closest trivial associates have
rather different supports.  These images are obtained as described above with
$\bF_1$ and $\bF_2$ non-negative images whose supports are inversion symmetric,
but the images themselves are not. The left and middle images in
Figure~\ref{fig2.1}(a) show the central $256\times 256$ portion of $\bF$ and $\bF',$
and the right image shows their common support.  The left image in
Figure~\ref{fig2.1}(b) is the $\log_{10}$-magnitude-DFT of both images,
the right image is the central $256\times 256$ portion of the difference of the
two images.

What is striking about this example is how perfectly ordinary the
images and their magnitude-DFT data look. The only criterion that we
know of (in the continuum model) to exclude this phenomenon is that it
cannot occur in an image with jump discontinuities, because the
convolution of two bounded measurable functions is continuous. For the
discrete model it is difficult to make this statement precise, as
there are images $\bF=\bF_1\ast \bF_2,$ where, say, $\bF_1$ is a ``sum
of $\delta$-functions,'' which provide counterexamples.  Note that,
for discrete images, a sum of $\delta$-functions is modeled by an
image with support in a set of isolated pixels.  In fact such examples
can be found in~\cite{FienupSeldin:90}. A complete (asymptotic)
analysis of this problem might require an analysis, as $q$ tends to
infinity, of the density of the subset of reducible polynomials of
degree $q$ within the set of all polynomials of this degree. For
results in this direction see~\cite{KaltofenMay2003}

\subsubsection{Microlocal Non-Uniqueness}\label{sec3.2}

There is a second mechanism that leads to $\epsilon$-non-uniqueness,
which we call \emph{microlocal non-uniqueness}. We now explain the
mechanism underlying this phenomenon.  For the construction, let $S$
be a set so that images supported in $S$ have small support.   Let
$\bF\in\TA_{\ba}$ be an image that can be decomposed as a sum,
\begin{equation}
  \bF=\bF_1+\cdots+\bF_k,
\end{equation}
where the components have the following properties:
\begin{enumerate}
\item For each $1\leq l\leq k$ we have $\{\bj:\: \epsilon<|f_{l\bj}|\}\subset S.$
\item Each pair $1\leq l\neq m\leq k$ has distinct spectral $\epsilon$-support,
  \begin{equation}
  \{\bj:\: \epsilon<|\hat f_{l\bj}|\}\cap \{\bj:\: \epsilon<|\hat f_{m\bj}|\}=\emptyset.
  \end{equation}
\end{enumerate}

From the second condition, it follows that, for $\bV=\{\bv_l\in J:1\leq
  l\leq m\},$ and $\Beta=(\beta_1,\dots,\beta_k),$ a binary string, the set of images
  \begin{equation}\label{eqn26}
    \bF^{\bV,\Beta}=\sum_{l=1}^k(-1)^{\beta_l}\bF_l^{(\bv_l)},
  \end{equation}
all have the same DFT magnitude data to precision $k\epsilon.$ Indeed,
we are also free to replace some of the $\bF_l$ with their
inversions $\cbf_l$. For a realistic estimate $S$ of the support,
images $\bF^{\bV,\Beta}$ defined by a collection of small
translations $\bV$ also have their supports within $S$ to precision
$k\epsilon$. In this way a large collection of images, which are not
trivial associates, can be constructed that belong to the intersection
$\TA_{\ba}\cap B_S$ up to a fixed, very small, error. The data for any one of
these images is $m\epsilon$-non-unique, for some fixed $m$. 

\begin{figure}[t]   
    \centering
\begin{subfigure}[t]{.4\textwidth}
     \centerline{\includegraphics[width= 5.5cm, height=5.5cm]{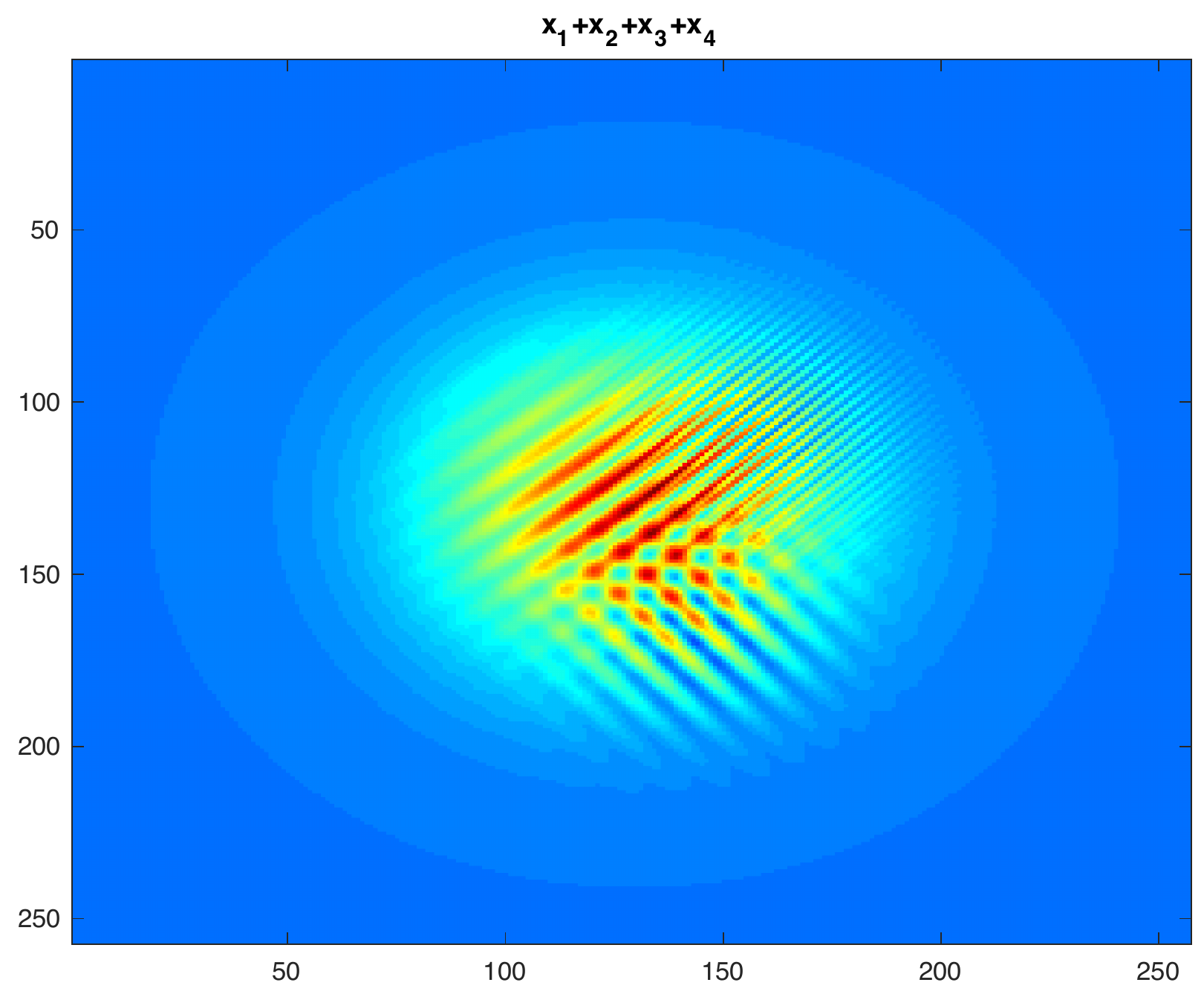}}
\caption{A sum of 4 component Gaussians.}
\end{subfigure}\qquad
\begin{subfigure}[t]{.4\textwidth}
     \centerline{\includegraphics[width= 5.5cm, height=5.5cm]{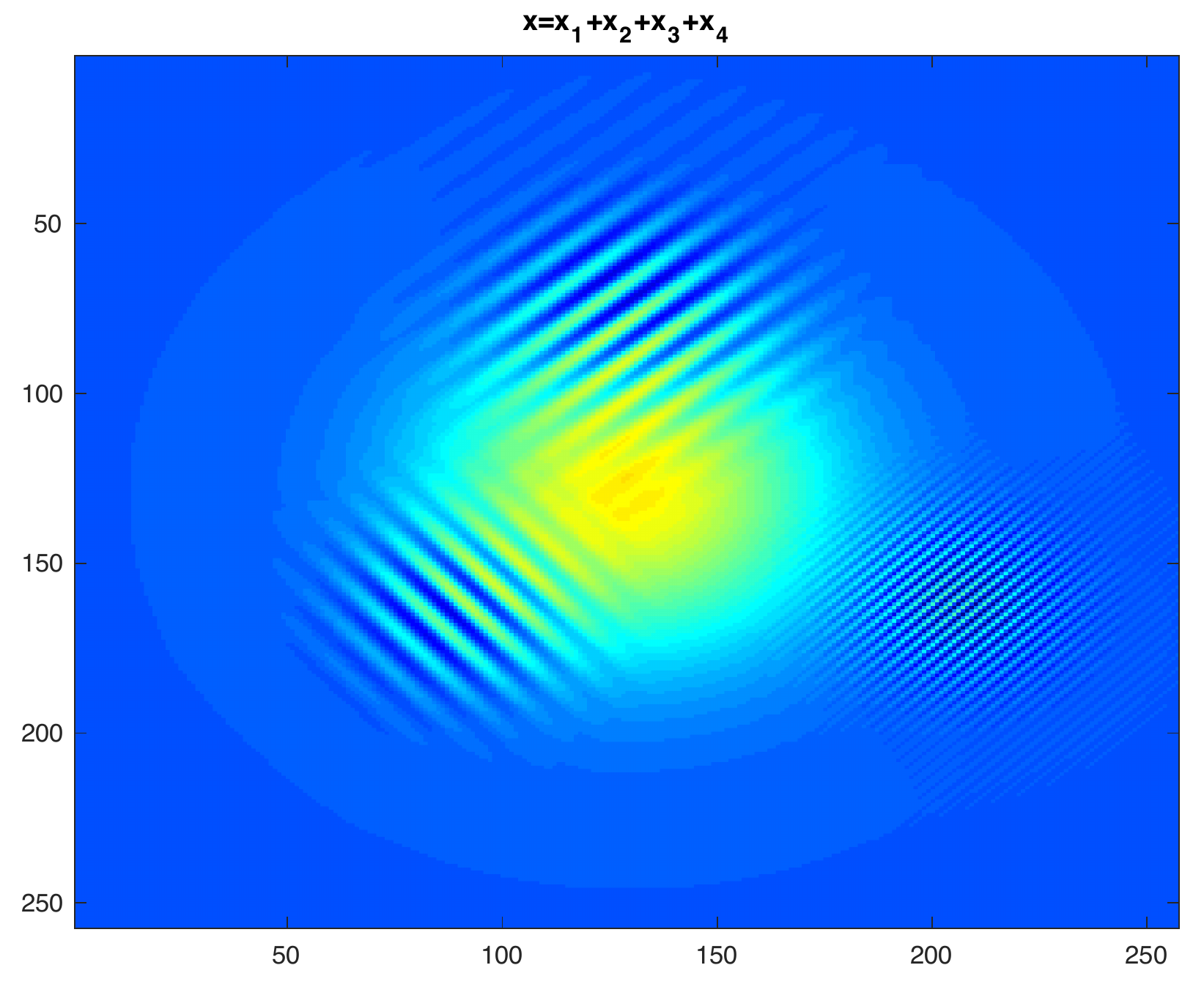}}
\caption{A different sum of 4 component Gaussians.}
\end{subfigure}
\begin{subfigure}[t]{.4\textwidth}
     \centerline{\includegraphics[height= 4.5cm]{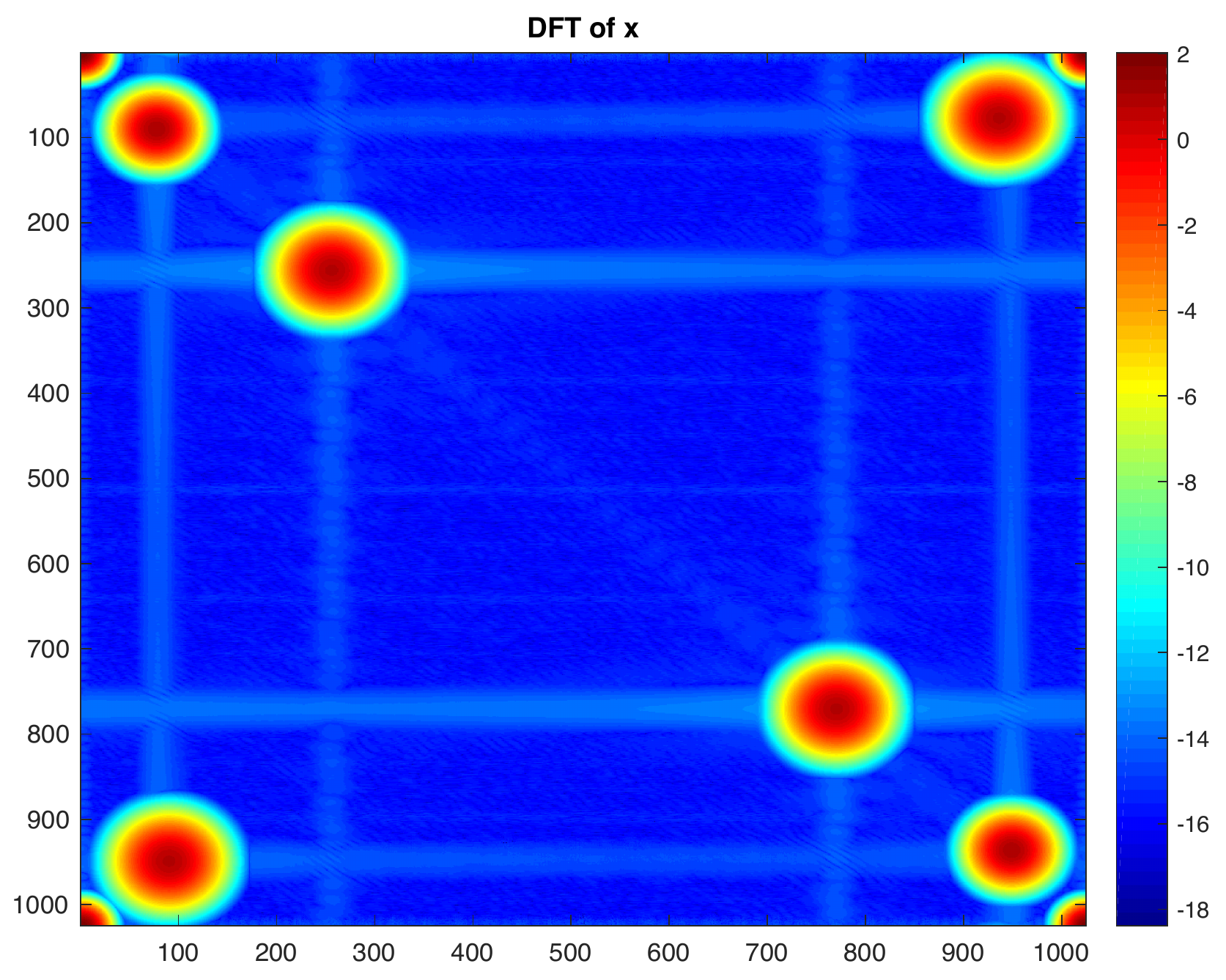}}
\caption{The $\log_{10}$-DFT magnitude data for the object in (a).}
\end{subfigure}\qquad
\begin{subfigure}[t]{.4\textwidth}
     \centerline{\includegraphics[height= 4.5cm]{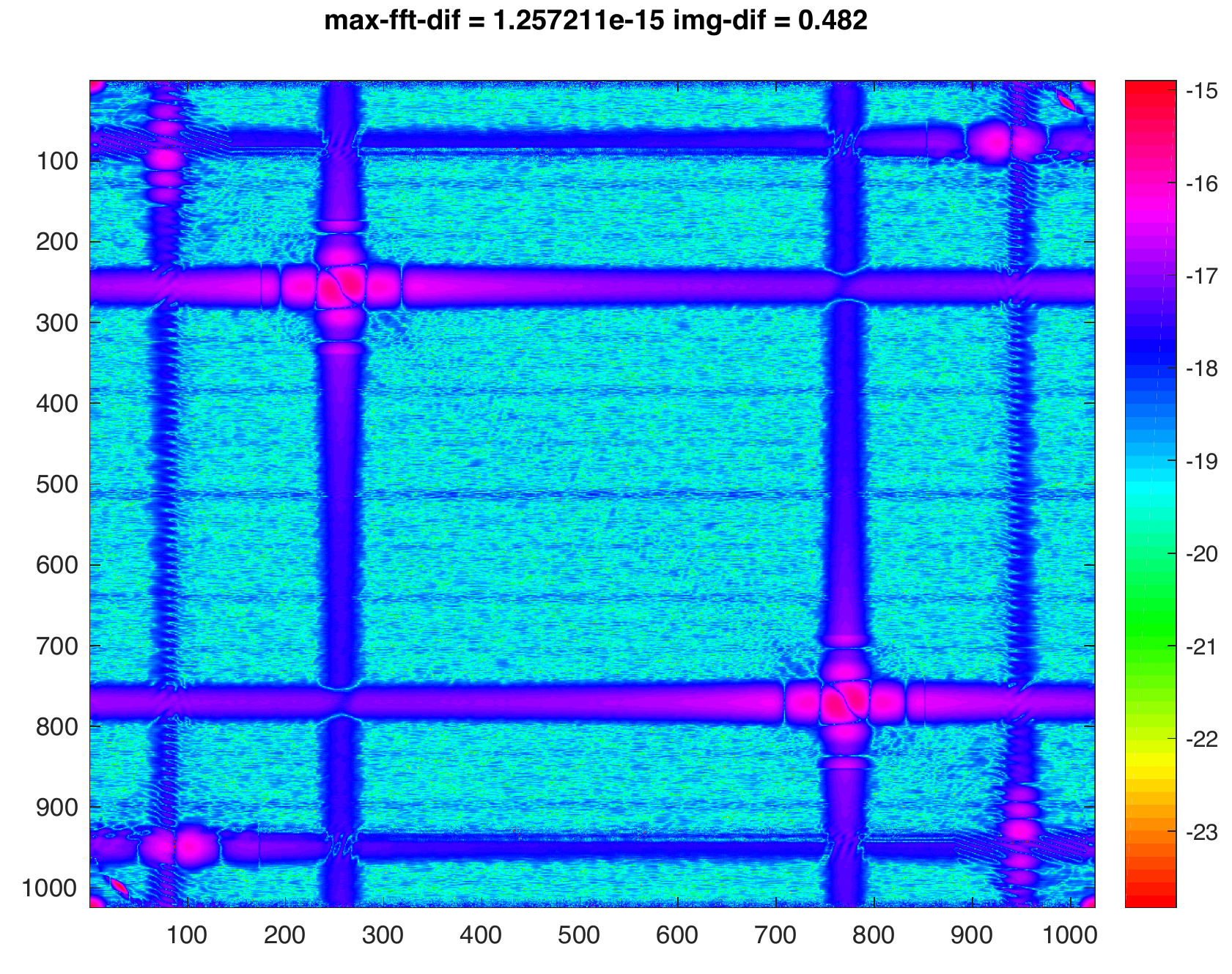}}
\caption{The $\log_{10}$ of the magnitude of the differences of
the DFT magnitude data for the objects in [a] and [b].}
\end{subfigure}
\caption{The top row shows two different objects $\bF_a$ and $\bF_b$ defined
in \eqref{illposedimgs}.
More precisely, we plot the central $256\times 256$ portion of the
$1024\times 1024$ array used in constructing these examples.
(c) is a plot of
of the $\log_{10}$ of the magnitude DFT data for the object in (a), while (d) is
the $\log_{10}$ of the difference of the magnitude DFT data for the objects in
(a) and (b). Note that the maximum difference is about $10^{-15}.$}
\label{fig2}
\end{figure}

 We close this section with an example of a pair images, whose
 difference is $O(1)$, but with identical support, and magnitude DFT
 data to precision $\epsilon=10^{-12}.$ In this example, which is
 shown in Figure~\ref{fig2}, $k=4.$

\begin{example}\label{illposed}
  In $d=2$, let
  the four images $\bF_1$ through $\bF_4$ be defined by their components
\begin{equation}\label{eqn28}
f_{i\bj}=e^{-\sigma_i^2|\bj-\bl_i|^2}\cos(\langle\bk_i,\bj-\bl_i\rangle),\quad i=1,2,3,4,
\end{equation}
where
\begin{equation}
  \bk_1=\bzero,\quad \bk_2=(70, 60),\quad \bk_3=(-60,70),\quad \bk_4=(200,200),
\end{equation}
and $N=512$, so $J=\{0,1,\dots,1023\}^2$.
Then we construct
\begin{align}
\begin{split}
\bF_a &= \bF_{1}+\bF_{2}+\bF_{3}+\bF_{4}~,\\
\bF_b &= \bF_{1}+\bF^{(\bv_2)}_{2}+\bF^{(\bv_3)}_{3}+\bF^{(\bv_4)}_{4},
\label{illposedimgs}
\end{split}
\end{align}
where the translation vectors are given by 
$$\bv_2 = (-8,0), \quad \bv_3 = (0,-8), \quad \bv_4 = (8,8).$$

Figure~\ref{fig2}(a,b) shows a plot of $\bF_a$ and $\bF_b$.
The support sets are defined as
$S=\{\bj:|f_{1\bj}+f^{(\bv_2)}_{2\bj}+f^{(\bv_3)}_{3\bj}+f^{(\bv_4)}_{4\bj}|
>10^{-12}\}.$
In both cases this is a disk of diameter $475$ pixels, thus the support is
small.
Yet the magnitude DFT data of these images are equal to precision
$10^{-15}$, thus phase retrieval is
incapable of distinguishing $\bF_a$ from $\bF_b$,
even if the data is measured to, say, 12 digits of accuracy.
\end{example}

\begin{remark}
The reader may note that our construction is somewhat pathological,
since the DFT of the image consists of well-separated Gaussian
``islands" of non-zero data, which leads to easier detection of this
sort of $\epsilon$-non-uniqueness.  The examples in Section
\ref{nonuniquesec} are less pathological and this form of
non-uniqueness is more difficult to detect. It remains an open problem
to describe the class of images for which $\epsilon$-uniqueness can be
proven, even for very small values of $\epsilon.$
\end{remark}

\section{Algorithms for Phase Retrieval}\label{sec4}

 We now see what the results of the previous sections imply about the
 behavior of standard algorithms used for phase retrieval. These
 algorithms are defined by iterating maps, which are, in turn, built
 from ``closest point maps.'' If $W$ is a subset of $\bbR^J,$ then
 $P_W(\bF)$ is defined to be the point in $W$ closest to $\bF$ with
 respect to the Euclidean distance. If $W$ is a linear subspace then
 $P_W$ is the orthogonal projection. If $W$ is convex then $P_W$ is
 defined and continuous everywhere, whereas for a non-convex set,
 these maps are defined, and continuous, on the complement of a
 positive codimensional subset.

 \begin{figure}  
  \includegraphics[width=2.5in]{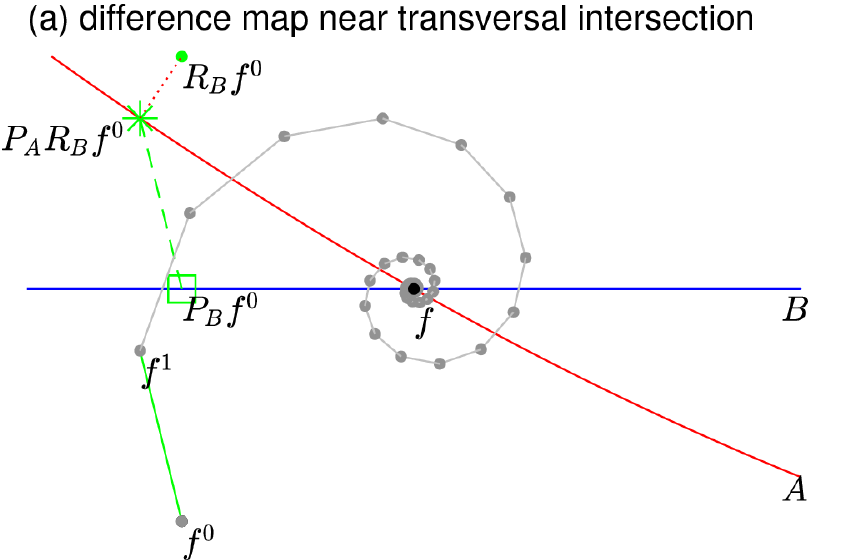}
  \includegraphics[width=2.9in]{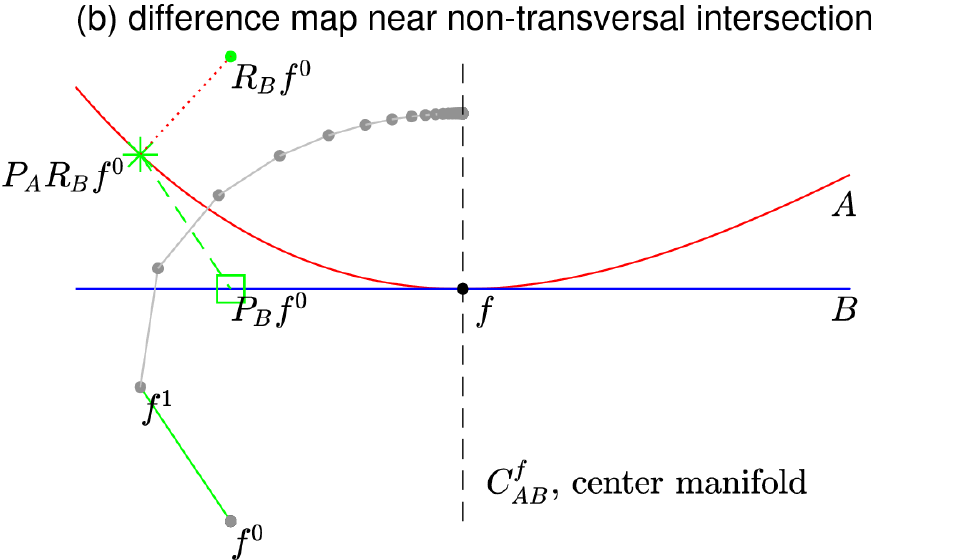}
  \caption{The hybrid iterative map $\bF^{n+1} = D_{AB}(\bF^n)$
    defined by \eqref{DAB},
    in the setting where $A$ (red curve) and $B$ (blue curve)
    are 1D manifolds in $\mathbb{R}^2$.
    (a) $\bF \in A \cap B$ is a transversal intersection
    (here the center manifold is the single point $\bF$).
    (b) Non-transversal case, with generic quadratic separation between
    the manifolds.
    Each plot shows the iterates $\bF^0,\bF^1,\dots$ (grey dots),
    and the construction of the update vector (green solid line)
    $\bF^1-\bF^0$ as the difference between a projected reflection
    (green star) and a projection (green square).
    Note that in each plot the green solid and dotted lines are equal as displacement vectors.
    \label{f:diffmap}
  }
 \end{figure}

In the phase retrieval problem, let us assume that the unknown $\bF$
and its support $S$ are adequate in the sense of definition
\eqref{adequatedef}, with the magnitude torus $\TA_{\ba}$ defined by
the DFT magnitude data $\ba=\ba_{\bF}$.  
As a torus, $\TA_{\ba}$ is obviously not a convex
set.  The alternating projection method \eqref{apiter}, which we
write here in the form
\begin{equation}
  \bF^{n+1} := P_{\TA_{\ba}}\circ P_{B}(\bF^{n})
\end{equation}
is well known to be prone to converge to points that are not in
$\TA_{\ba}\cap B.$ The stable fixed points of the alternating
projection map are points $\tilde{\bF}\in\TA_{\ba}$ such that
$\tilde{\bF}$ and $P_B(\tilde{\bF})$ jointly form a non-zero, local
minimum of the Euclidean distance between the two sets.  Empirically
these exist in great profusion, and the alternating projection method
rarely (if ever) converges to a point in $\TA_{\ba}\cap B.$ See
\cite{BEGM} for a more extensive discussion.

In an attempt to avoid such false local minima and improve
convergence, a variety of modifications to the alternating projection
map have been introduced that involve reflection operators as well as
projections.  Quite a few variants have appeared in the literature
\cite{bauschke2002,chapman2006,elser2003,elser2007,fienup1987,miao1999,miao2015},
and we will limit our attention to a special case of Fienup's hybrid
input-output (HIO) method \cite{fienup1982}, which is also a special
case of the ``difference map'' approach due to Elser {\em et al.}
\cite{elser2003}.  Letting $A$ and $B$ now denote general sets, with
closest point projections $P_A$ and $P_B$, the map that we iterate is
\begin{equation}
  D_{AB}(\bF)\;:=\;\bF+P_{A}\circ R_B(\bF)-P_B(\bF),
  \label{DAB}
\end{equation}
where $R_B$ is the ``reflection'' around $B$ defined by
$R_B(\bF):=2P_B(\bF)-\bF$; see Figure~\ref{f:diffmap}.  This is
Fienup's HIO method with $\beta = 1$ and a specific instance of
Elser's, difference map as well. If $A$ is a linear subspace, then
$D_{AB}$ is also the Douglas-Rachford map, which is defined to be
\begin{equation}
  T_{B,A}=\frac{1}{2}\left[R_A\circ R_B+\Id\right],
\end{equation}
see~\cite{borwein}.

Since we are not testing all possible HIO or difference map variants,
we will call the specific method we use here a ``hybrid iterative map."
If $\bF^*$ is a fixed point of $D_{AB}$ then
\begin{equation}
  P_{A}\circ R_B(\bF^*)=P_B(\bF^*),
\end{equation}
in other words, the point $\bF^{**}:=P_B(\bF^*)$ lies in $A\cap B$. The
iterates are defined by $\bF^{n+1}=D_{AB}(\bF^n)$, and approximate
reconstructions are given by
\begin{equation}
  \br^n := P_B(\bF^n)~.
  \label{r}
\end{equation}
If the iterates
converge, then, assuming that $P_B$ is continuous at the limit point, the
sequence $ \{ \br^n \}$ converges to a point in $A\cap B.$

The fixed point set of $D_{AB}$ can be much larger than the set of
intersections. Given a point $\bF\in A\cap B$, we let
\begin{equation}
  L_A:=P_A^{-1}(\bF) \qquad \text{ and } \quad L_B:=P_B^{-1}(\bF).
\end{equation}
The \emph{center manifold} defined by $\bF$ (see Figure~\ref{f:diffmap})
is then the set
\begin{equation}
  \cC_{AB}^{\bF}:=R_{B}^{-1}(L_A)\cap L_B.
\end{equation}
The center manifold for any $\bF\in A \cap B$ contains the part of the
fixed point set for the map $D_{AB},$ which ``points to'' $\bF.$ Given
an image $\bF$ with support in a small support set $S$ (see Definition
\ref{ssdef}), and choosing $B=B_S$ and $A=\TA_{\ba}$, then the part of
the center manifold near to $\bF$ is given by
$\cC_{\TA_{\ba}B_S}^{\bF}=(\bF+B_S^{\bot})\cap N_{\bF}\TA_{\ba}$, and with
$\dim\cC_{\TA_{\ba}B_S}^{\bF}>|J|/4$.  However the map
$D_{\TA_{\ba}B_S}\restrictedto_{\cC_{\TA_{\ba}B_S}^{\bF}}=\Id$, which
is only neutrally stable.

While it is true that the fixed point sets are contained in the center
manifolds defined by points in $\TA_{\ba}\cap B_S,$ there are other
subsets that are attracting. If the pair of points $(\bF_1,\bF_2)\in\TA_{\ba}\times B_S$
defines a critical point of the map $d_{\TA_{\ba}B_S}:\TA_{\ba}\times
B_s\to\bbR_+,$ 
\begin{equation}
  d_{\TA_{\ba}B_S}(\bF,\bg)=\|\bF-\bg\|_2,
\end{equation}
then the set
\begin{equation}
 \cC_{\TA_{\ba} B_S}^{\bF_1,\bF_2}\overset{d}{=}(\bF_2+B_S^{\bot})\cap N_{\bF_1}\TA_{\ba}
 \end{equation}
contains the line segment from $\bF_1$ to $\bF_2.$ Indeed this
intersection is again a subset with dimension about $|J|/4$ (if
$d=2$). From low dimensional examples it appears that these sets can
define attracting basins, even if the critical point at
$(\bF_1,\bF_2)$ is not a local minimum.  The existence of these
attracting sets seems to complicate the dynamics of hybrid map iterations.

The fiber of the tangent bundle
$T_{\bF}\TA_{\ba}$ at $\bF\in \TA_{\ba}$ is the best linear
approximation to $\TA_{\ba}$ near to $\bF$; hence a linearization of
the problem of locating points in $\TA_{\ba}\cap B_S$ is to locate
points in the intersection of the affine subspaces
$T_{\bF}\TA_{\ba}\cap B_S$. With this as motivation we first analyze
the behavior of the map $D_{AB}(\bF)$ when $A$ and $B$ are
linear subspaces.

\subsection{Linear Subspaces}\label{sec4.1}
For the case of a linear subspace, $W\subset\bbR^N,$ the map $P_W$ is
the orthogonal projection onto $W$ and $R_W$ is the orthogonal
reflection with fixed point set $W$. Let $A$ and $B$ denote linear
subspaces of $\bbR^N$.
Let us first consider the linear model for the
benign \emph{transversal} intersection case.
We have $A\cap B=\{\bzero\}$, i.e.\ a single
isolated point,
and $\dim A+\dim B < N$ as befits the phase retrieval application.
(For example, for $d=2$, $\dim A = |J|/2$ and $\dim B \le |J|/4$ when the constraint is adequate.)
To analyze the iteration defined by $D_{AB}$ we split $\bbR^N$ into
the following subspaces $A$, $B$, and $C:=(A+B)^{\bot}=A^{\bot}\cap
B^{\bot}.$ The subspace $C$ is the center manifold defined by
$\{\bzero\}=A\cap B$ for this case. Let $U$, $V$, $Y$ denote matrices whose
columns are orthonormal bases for $A$, $B$ and $C$ respectively. If
$\bF=U\bx_1+V\bx_2+Y\bx_3,$ then, in this representation, the map
$D_{AB}$ takes the form
\begin{equation}
  D_{AB}(\bx_1,\bx_2,\bx_3)
= \left(\begin{matrix} 2H^tH& H^t&0\\
-H& 0 & 0\\0&0&\Id\end{matrix}\right)\left(\begin{matrix}\bx_1\\\bx_2\\\bx_3\end{matrix}\right),
\end{equation}
where $H=V^tU$.
In~\cite{BEGM}, it is shown that the upper $2\times 2$ block matrix 
$\left(\begin{matrix} 2H^tH& H^t\\
-H& 0 \end{matrix}\right)$
is a contraction, and therefore the
map is contracting in directions normal to the center manifold $C$, and
$\lim_{n\to\infty}D^n_{AB}(\bx_1,\bx_2,\bx_3)=(\bzero,\bzero,\bx_3)$.
This contraction is visible as the convergent spiral in
Figure~\ref{f:diffmap}(a), where $C=\{\bF\}$.
Its rate of contraction is determined by largest singular value of $H$.
The limit point then yields the desired solution
under the projection $P_B$.

The correct linear model for a \emph{non-transversal} intersection is
similar, but $A\cap B=F$ is now a subspace of positive dimension. We now
split $\bbR^N$ as $A_0+B_0+F+C,$ where $A_0=A\cap F^{\bot},$ 
$B_0=B\cap F^{\bot},$ and $C=(A+B)^{\bot}.$ If $U$, $V$, $X$, $Y$
denote matrices whose columns are
orthonormal bases for $A_0$, $B_0$, $F$, $C$ respectively, then, with
$\bF=U\bx_1+V\bx_2+X\bx_3+Y\bx_4,$ we have:
\begin{equation}
  D_{AB}(\bx_1,\bx_2,\bx_3,\bx_4)=\left(\begin{matrix} 2H^tH&H^t&0&0\\
-H & 0& 0 & 0 \\
0 & 0& \Id & 0 \\0&0&0&\Id\end{matrix}\right)\left(\begin{matrix}\bx_1\\ \bx_2
\\\bx_3\\\bx_4\end{matrix}\right).
\end{equation}
As before, the leading $2\times 2$ block is a contraction, and therefore
\[ \lim_{n\to\infty}D^n_{AB}(\bx_1,\bx_2,\bx_3,\bx_4)=
(\bzero,\bzero,\bx_3,\bx_4). \] 
Crucially, in this linear model $D_{AB}$ is the identity operator in both
the $F$- and $C$-directions.

For the non-linear phase retrieval problem, the intersections of
interest are, as shown above, generally non-transversal. In this case,
linearization at the intersection point tells one nothing about the
map's behavior, even very near to the center manifold.  More
precisely, because the intersections of $\TA_{\ba}\cap B_S$ are
isolated points $\bF$, the affine subspace $C=N_{\bF}\TA_{\ba}\cap
B_S^{\bot}$ remains the linear model for the center manifold in the
non-linear case. The subspace $F=T^0_{\bF}\TA_{\ba}\cap B_S$
constitutes a positive-dimensional set normal to the center manifold
where the map $D_{AB}$ is not known to be contracting.

In the 2-dimensional, quadratic non-transversal case shown in
Figure~\ref{f:diffmap}(b), one observes geometric contraction towards
$C$.  In simple, low dimensional examples of this sort convergence,
even in the non-transversal case, is often observed. In fact the
problem becomes, in some sense, easier as non-transversality causes
the dimension of the target center manifold to increase.  On the other
hand, very small, but non-zero angles lead to very slow convergence.
The spiral trajectory in Figure~\ref{f:diffmap}(a) is also
note-worthy, as it indicates that the linearization of $D_{AB}$ at the
limit point has complex eigenvalues, which is a phenomenon that
persists in the phase retrieval problem. These questions are discussed
in detail in~\cite{BEGM}.

In the phase retrieval problem, where the dimension $|J|$ is large,
and the geometry is much more complicated, non-transversality seems to
preclude convergence. Indeed, it is common to find that the hybrid
map iteration \emph{stagnates} at a substantial distance from the
center manifold, whenever $\dim T^0_{\bF}\TA_{\ba}\cap B_S>0.$

\begin{definition}
  An algorithm has stagnated if the distances from subsequent iterates to the
  nearest exact intersection point remain almost constant, and much larger than
  machine precision; moreover the distances between successive iterates are
  also essentially constant, and much larger than machine precision.
\end{definition}

This behavior is almost always observed when using hybrid iterative 
map-based algorithms on noise-free data coming from
images that are not tightly constrained by the support mask,
as we show next.
The
failure of transversality not only renders the problem of finding points in
$\TA_{\ba}\cap B$ ill-posed, but also prevents standard algorithms for finding
these points from converging.

In applications, a variety of possible support information is possible,
such as a bounding rectangle, bounding disc, etc.
Here we use quite an optimistic estimate for our
knowledge of the support, namely that the true support is known
up to a ``padding'' of $p$ pixels. The notion of a $p$-pixel
neighborhood is defined in~\eqref{eqn2.10.002}.
In applications this might possibly derive from
knowledge of a lower-resolution version of the
target image; note that it includes much more information than
merely a reasonably accurate bounding rectangle.

\section{The Performance of the Hybrid Iterative Maps} \label{sec5}

The theory presented in the previous sections makes rather specific predictions
as to how the hybrid iterative map will behave on various sorts
of images, and different sorts of auxiliary information. 
For the support constraint, one expects to see that the
iterates $\{\bF^n\}$ of such a map stagnate, and that
the differences between the
approximate reconstructions $\{\br^n\}$ and the nearest exact intersection
point, $\bF$, should lie mostly in directions belonging to
$T_{\bF}\TA_{\ba}\cap B_S$. 
In Section~\ref{sec5.1} we show that 
these predictions are largely verified in
practice. 

In Section~\ref{sec5.2}, we instead consider the non-negativity
constraint. While we still assume that the image has small support,
that information is not explicitly used.
When $B=B_+$, $\pa B_+$ is not a smooth space,
but is rather stratified by the number of vanishing coordinates.
The strata are orthants in Euclidean spaces of various dimensions and
the intersections with
$\TA_{\ba}$ lie on the boundary of the orthant.
In light of this, the intersection $T_{\bF}\TA_{\ba}\cap\pa B_+$ is a 
reasonable measure of the transversality of the intersection 
$\TA_{\ba}\cap B_+$ at $\bF$. 
The more coordinates that vanish at a point, the more
directions in which $\pa B_+$ is strictly convex near to that point.
This suggests that the intersections
between $\TA_{\ba}$ and $B_+$ have a better chance to be transversal,
and therefore hybrid map-based algorithms should work better with
this auxiliary information. We will see, in Section~\ref{sec5.2}, that both
expectations are indeed true. 

If $\bF\in\TA_{\ba}$ is a non-negative image, then it is obvious that
the zeroth DFT coefficient $\hat f_{\bzero}=\|\bF\|_1,$ the
$\ell_1$-norm of $\bF$. As follows from the triangle inequality, the $\ell_1$-norm is strictly
minimized on the magnitude torus $\TA_{\ba}$ exactly at such
single-signed images.  Hence
for non-negative images one can use the $\ell_1$-norm to define a
different constraint, and therefore different algorithms.  Let
$B^1_{r}$ denote the $\ell_1$ ball of radius $r=|\hat f_{\bzero}|$. The
analysis of the intersection $T_{\bF}{\TA_{\ba}}\cap {\pa B_+},$ where
$\bF$ is a non-negative image in $\TA_{\ba},$ has the somewhat
unexpected consequence that
\begin{equation}\label{eqn42}
  T_{\bF}{\TA_{\ba}}\cap {\pa B_+}=T_{\bF}{\TA_{\ba}}\cap {\pa B^1_r}.
\end{equation}
That is, the failure of transversality of these two intersections agree
exactly, and therefore algorithms based on using $B=B_+$ can be expected
to behave similarly to those using $B=B^1_r.$
We find that this is true, on average,
though, as the maps involved are non-linear, individual runs of these
algorithms can behave quite differently. This is also briefly explored
in Section~\ref{sec5.2}.

\begin{figure}[t]  
    \centering
    \begin{subfigure}[H]{.8\textwidth}
                \centerline{\includegraphics[height= 3.75cm]{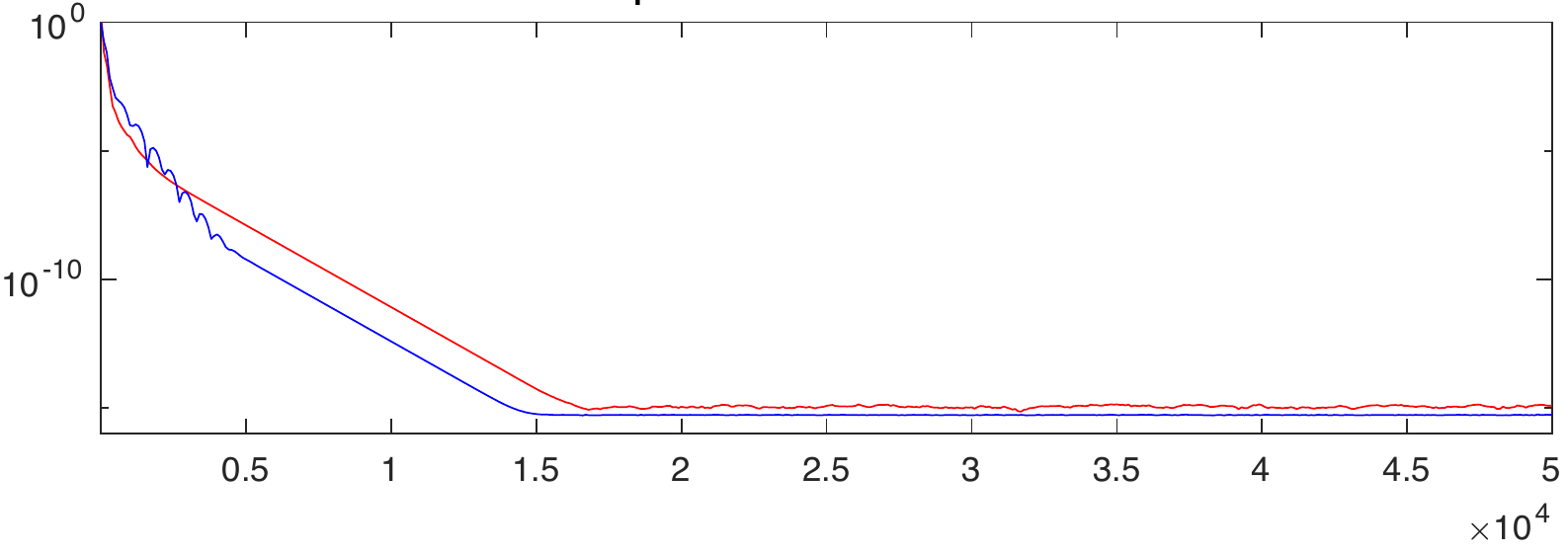}}
        \caption{ The iterates lie in an attracting basin with
          $\dim T_{\bF^{(2,2)}}\TA_{\ba}\cap B_{S_2}=0.$ }
    \end{subfigure}
\begin{subfigure}[H]{.8\textwidth}
                \centerline{\includegraphics[height= 3.75cm]{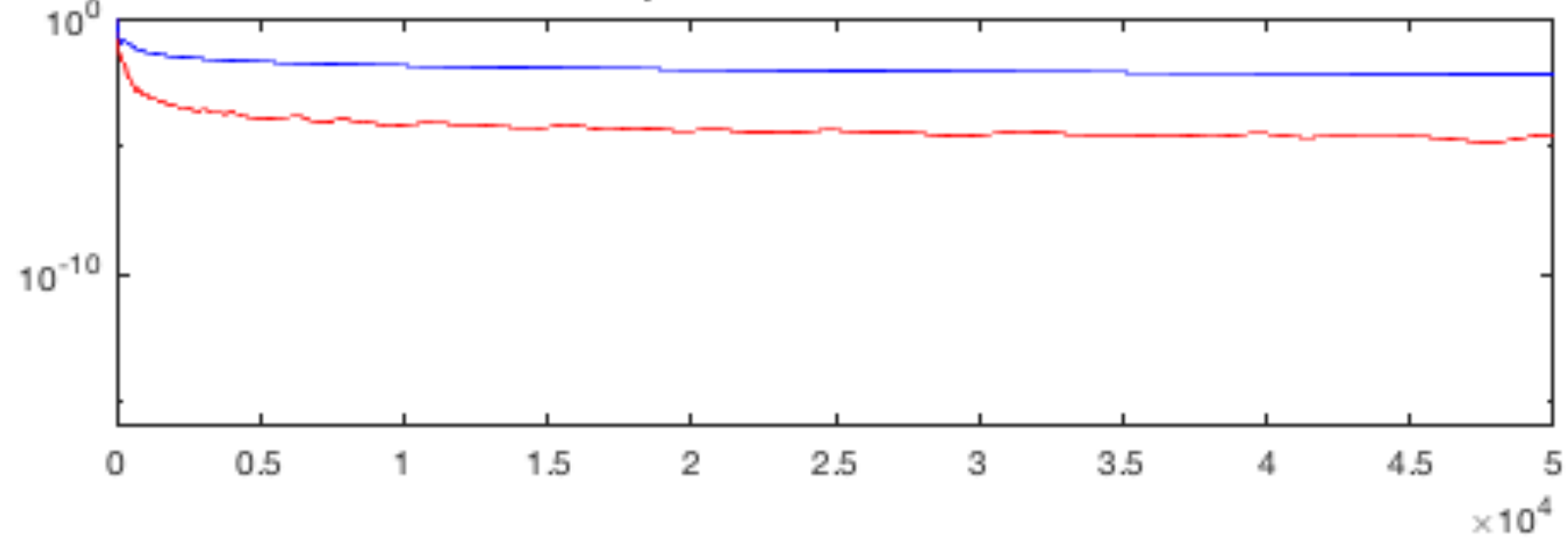}}
        \caption{The iterates lie in an attracting basin with
          $\dim T_{\bF^{(1,2)}}\TA_{\ba}\cap B_{S_2}=1.$ }
    \end{subfigure}
\caption{An illustration of how  the convergence properties of the
  hybrid iterative map using $D_{\TA_{\ba} B_{S_2}}$ 
  depend on the dimension of the 
  $\dim T_{\bF^{(\bv)}}\TA_{\ba}\cap B_{S_2}$. See Example~\ref{exam1}.
  The true errors are shown in
  blue, residual \eqref{E} in red.}
\label{fig5}
\end{figure}   

\subsection{The Support Constraint}\label{sec5.1}

In this section we examine the dependence of hybrid iterative maps
on smoothness ($k$) and padding of the support ($p$) for images of the
types used in Section~\ref{sec2.4.004}. When it is clear which image
is intended, we use $S_{p}$ to refer to $S_{\bF,p}.$

\begin{example}\label{exam1}
This example gives compelling evidence for the central importance of
the failure of transversality. We employ a piecewise
constant $256\times 256$ (i.e.\ $2N=256$) image $\bF$ with support
condition $S_{2}$, which is the exact support padded by $p=2$ pixels as
defined in~\eqref{eqn2.10.002}.  The intersection $\TA_{\ba}\cap
B_{S_2}$ contains 25 points, which are the trivial associates
$\{\bF^{(\bv)}: \; \|\bv\|_{\infty}\leq2\}$.
The dimension of each intersection
$T_{\bF^{(\bv)}}\TA_{\ba}\cap B_{S_2}$ depends on $\bv$. At
$\bv=\bzero$ this dimension attains the maximum of 12.
Each of the center manifolds
$\{\cC^{\bF^{(\bv)}}_{\TA_{\ba}B_{S_2}}:\|\bv\|_{\infty}\leq 2\}$
defines a basin of attraction for the hybrid map $D_{\TA_{\ba}B_{S_2}}.$
Starting at a random point on $\TA_{\ba}$ the iterates
seem to eventually fall into one of these basins of attraction.

Letting $\bF^{0}$ denote the starting point,
and writing $B = B_{S_2}$,
we have $\bF^{n}=D_{\TA_{\ba}B}(\bF^{n-1})$ the $n$th iterate. These
points are eventually close to points on a center manifold, but not
very close to the point in $\TA_{\ba}\cap B$ that defines it.
The sequence of approximate reconstructions is defined by \eqref{r}.
The plots in Figure~\ref{fig5} show the \emph{true error}
$\|\br^n-\bF^{(\bv)}\|_2$ (in blue), where
$\bF^{(\bv)}$ is the trivial associate of the true image
closest to $\br^n$, and the
\emph{residual} (in red), which is defined to be
\begin{equation}
  E(\bF^n) := \|P_{B}(\bF^n)-P_{\TA_{\ba}}\circ R_{B}(\bF^n)\|_2~.
  \label{E}
\end{equation}
Throughout this paper, the true errors are plotted in blue, and the
residuals in red. In a real experiment, only the residual is observable.

Recalling that $\ba$ is the data vector, the Lipschitz bound \eqref{lip}
implies
$$
\|\cM(\br^n) - \ba\|_2 \le C_{\cM} \|\br^n - P_{\TA_{\ba}}(\br^n)\|_2
 \leq C_{\cM} \|\br^n - P_{\TA_{\ba}}\circ R_{B}(\bF^n)\|_2 =C_{\cM} E(\bF^n),
$$
 where the middle inequality follows from the definition of $P_{\TA_{\ba}}$. This
 inequality shows that the data residual norm (left side), a measure
 of the extent to which the approximate reconstructions satisfy the DFT-magnitude
 constraints, is controlled by our plotted quantity \eqref{E}.
 For the hybrid map we also have $E(\bF^n)=\|\bF^{n+1}-\bF^{n}\|$,
 so \eqref{E} also provides an indicator as to whether the iterates
 are converging. Our plots are semilog-plots with a logarithmic
 $y$-axis; a linear decrease  therefore indicates  exponential
 (geometric) decay.

In Figure~\ref{fig5}(a) the iterates have settled into the attracting basin
defined by the associate
$\bF^{(2,2)}$. At this point the intersection with $B_{S_2}$ is
transversal; it is quite apparent that, by the 17,000th iterate, the
approximate reconstructions have converged to this intersection point to
machine precision. In Figure~\ref{fig5}(b) the iterates have settled into the
attracting basin defined by $\bF^{(1,2)},$ for which $\dim
T_{\bF^{(1,2)}}\TA_{\ba}\cap B_{S_2}=1$.
The iterates appear to have largely
stagnated after about the 35,000th iterate, with an true error of about $10^{-2}$
and a residual of about $10^{-4}.$ It follows from~\eqref{eqn21}, and the
observation
that the residual is about the square of the true error, that the 
differences
$\br^n-\bF^{(1,2)}$ are likely to lie largely along a common tangent
direction. Indeed, a more careful analysis of these differences, given
in~\cite{BEGM}, verifies this expectation.
\end{example}

\begin{figure}[t]   
    \centering
    \begin{subfigure}[H]{.45\textwidth}
                \includegraphics[height= 4.5cm]{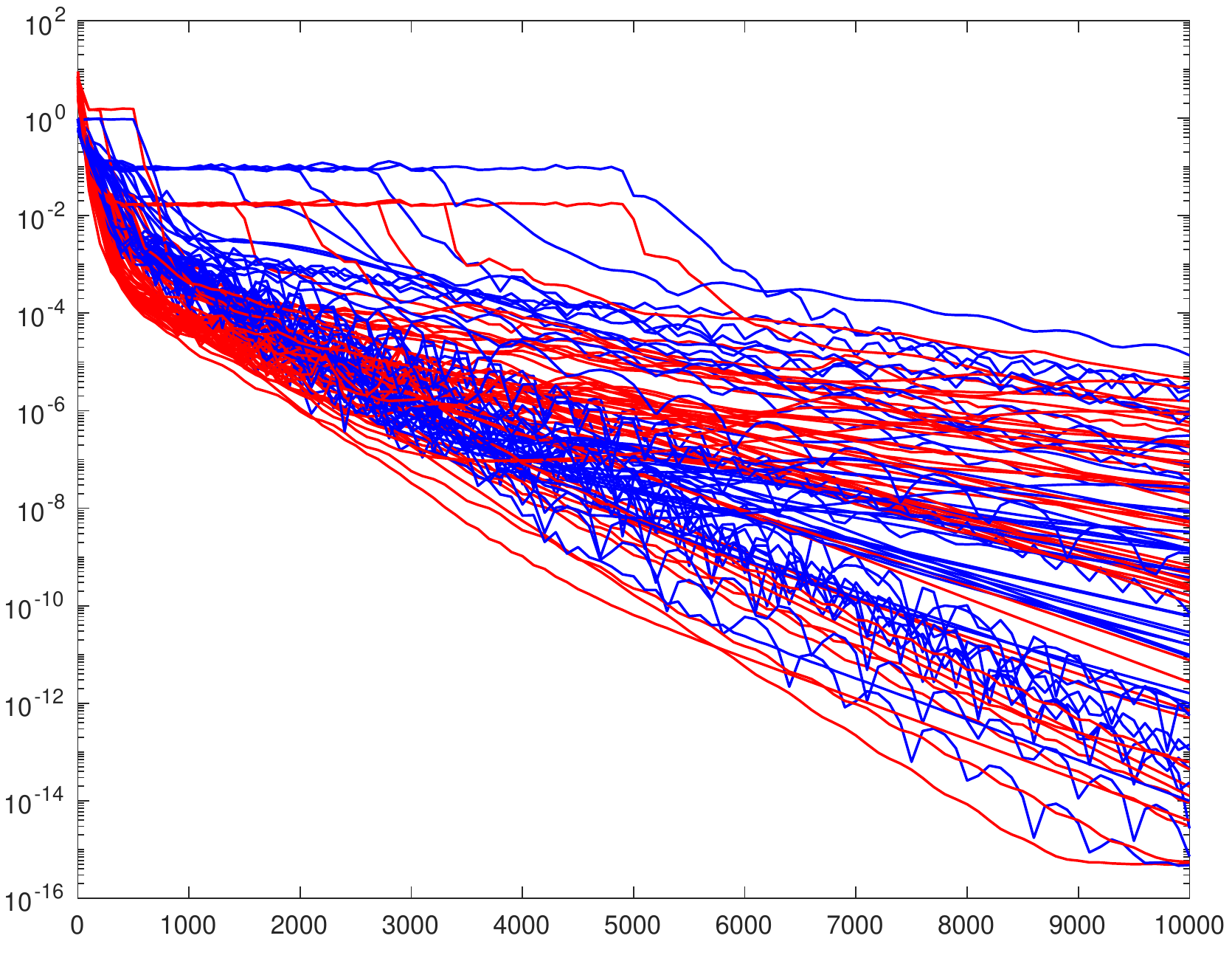}
        \caption{True errors and residuals for $D_{\TA_{\ba} B_{S_p}}$ with $k=0, p=1.$}
    \end{subfigure}\qquad
\begin{subfigure}[H]{.45\textwidth}
                \includegraphics[height= 4.5cm]{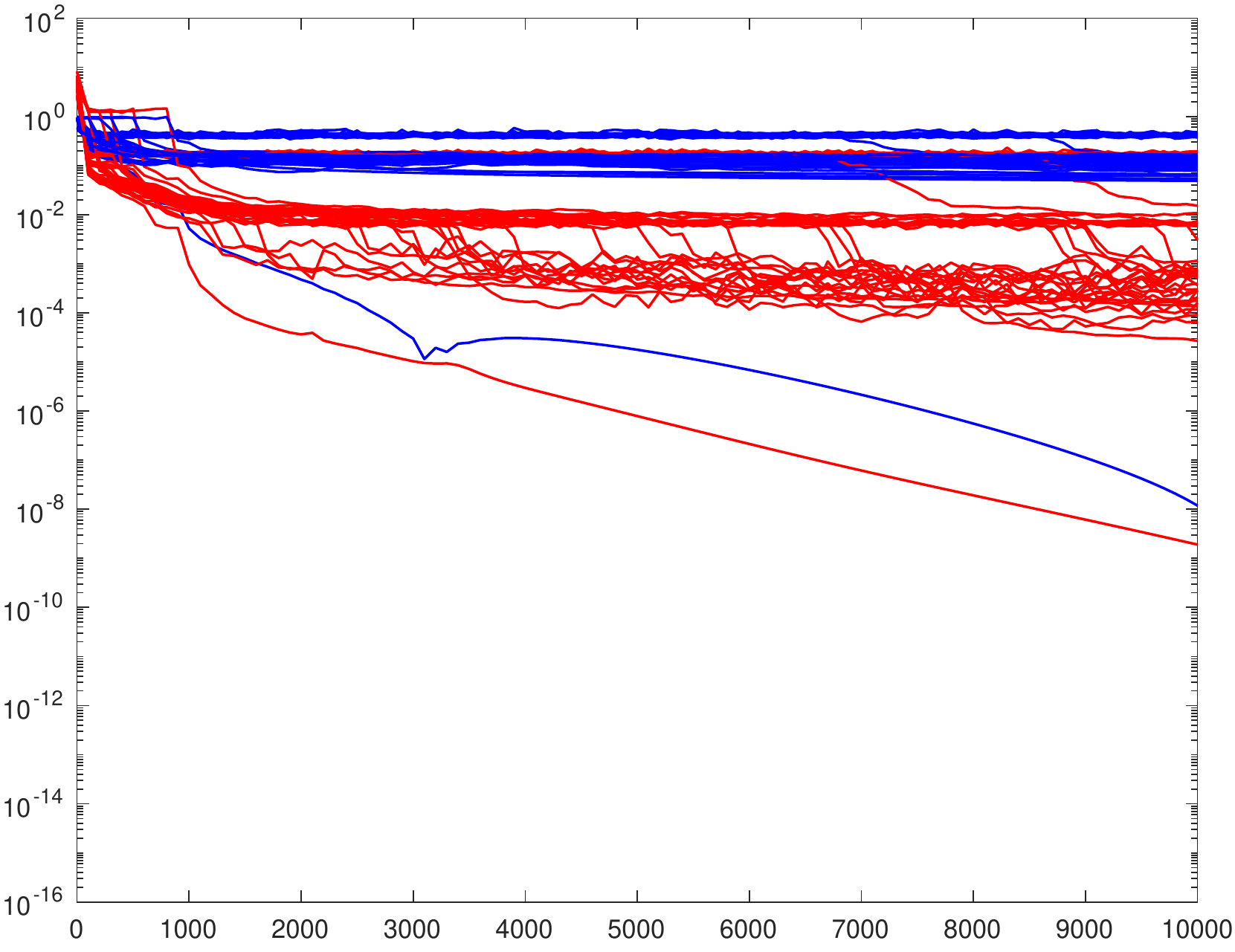}
        \caption{True errors and residuals for $D_{\TA_{\ba} B_{S_p}}$ with $k=0, p=3.$}
\end{subfigure}
\caption{The convergence properties of the hybrid map algorithm
  $D_{\TA_{\ba} B_{S_p}},$ with $k=0$; $p=1,3,$ 50 random restarts on
  a $256\times 256$-image. See Example~\ref{e:startpt}.}
\label{Fig6}
  \end{figure}  

    \begin{figure}[t] 
    \centering
    \begin{subfigure}[H]{.45\textwidth}
                \includegraphics[height= 4.5cm]{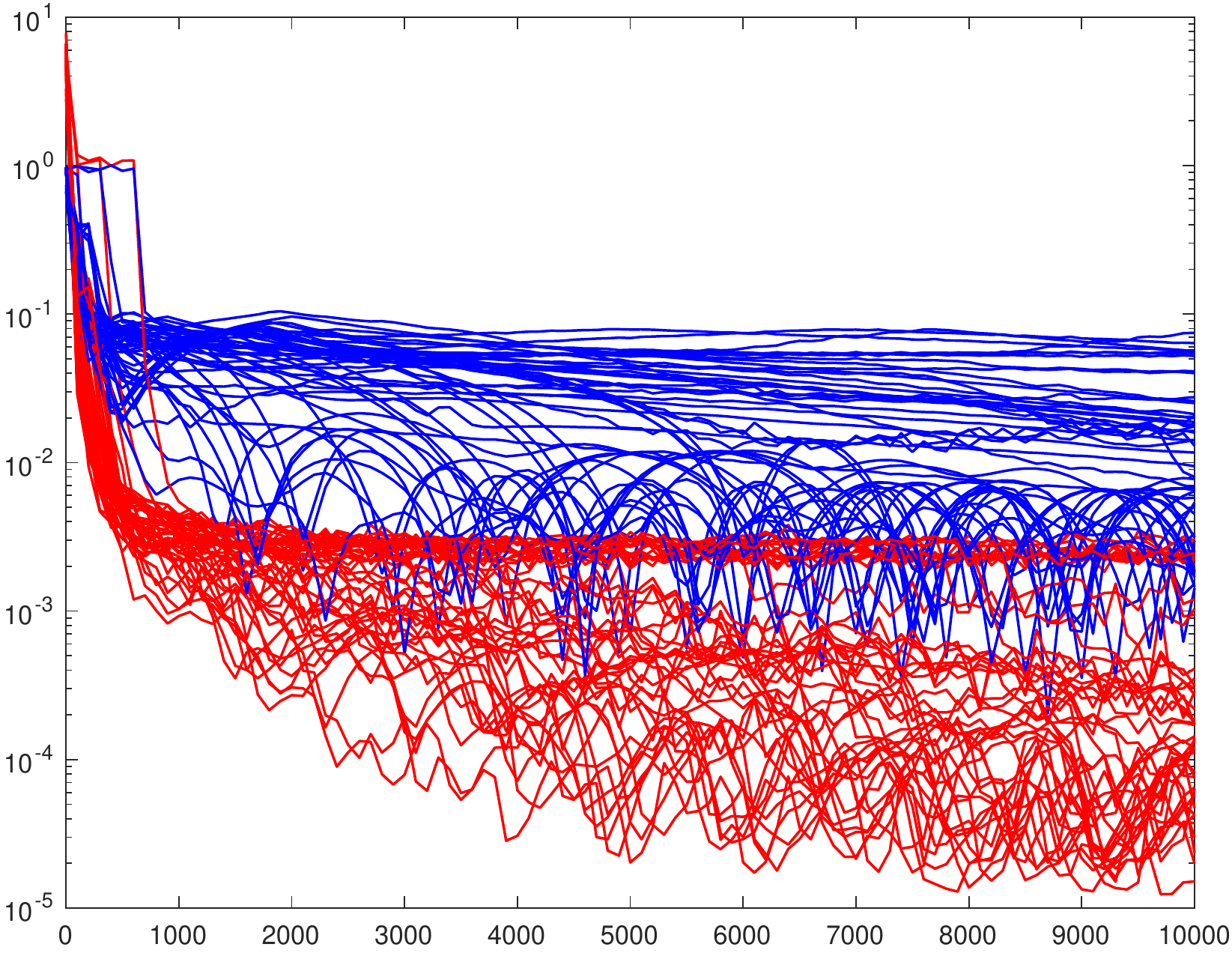}
        \caption{True errors and residuals for $D_{\TA_{\ba} B_{S_1}}$ with
          $k=2, p=1,$ without convolution.}
    \end{subfigure}\qquad
\begin{subfigure}[H]{.45\textwidth}
                \includegraphics[height= 4.5cm]{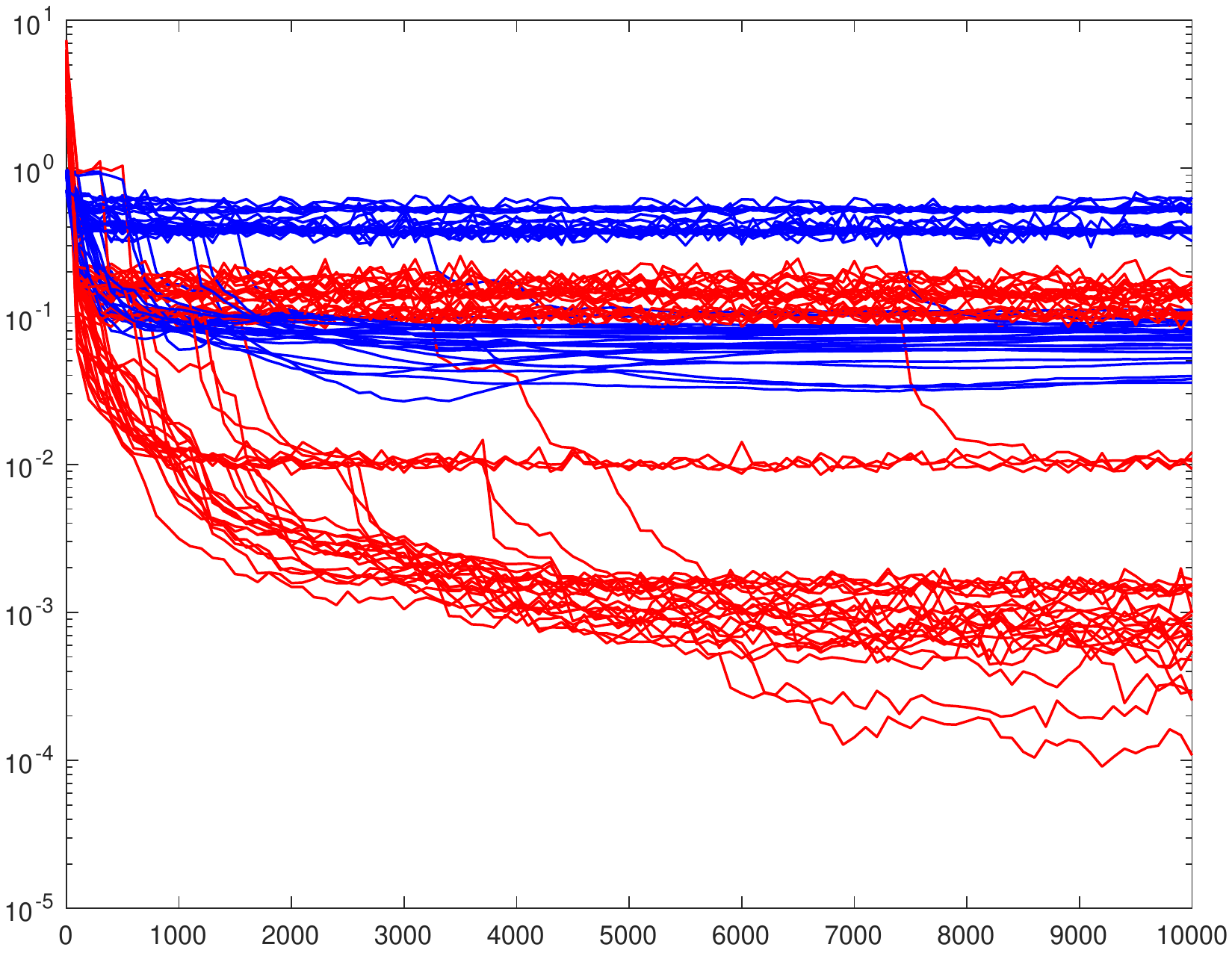}
        \caption{True errors and residuals for $D_{\TA_{\ba} B_{S_3}}$ with
          $k=2, p=3,$ without convolution.}
\end{subfigure}
\begin{subfigure}[H]{.45\textwidth}
                \includegraphics[height= 4.5cm]{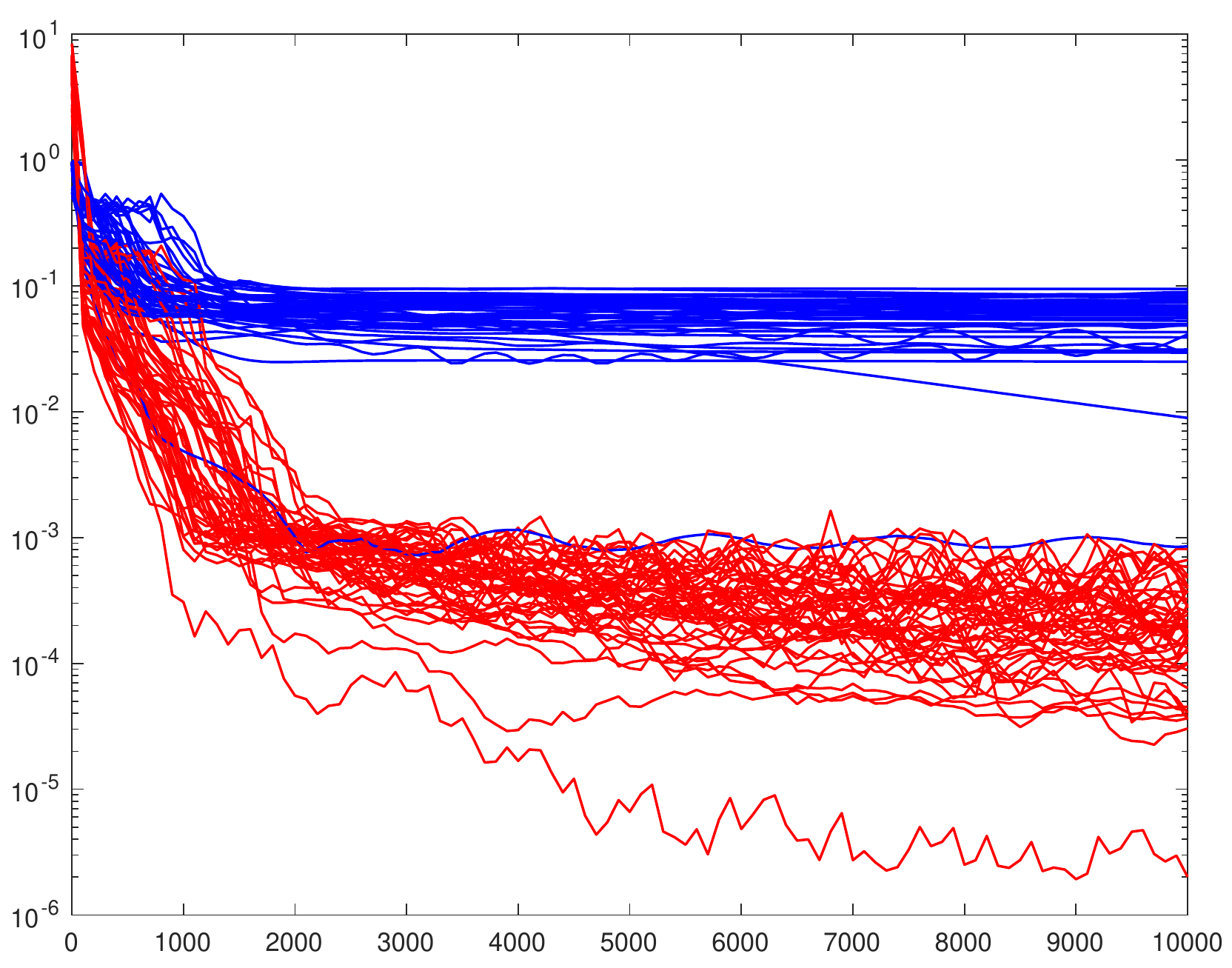}
        \caption{True errors and residuals for $D_{\TA_{\ba} B_{S_1}}$ with
          $k=2, p=1,$ with convolution.}
    \end{subfigure}\qquad
\begin{subfigure}[H]{.45\textwidth}
                \includegraphics[height= 4.5cm]{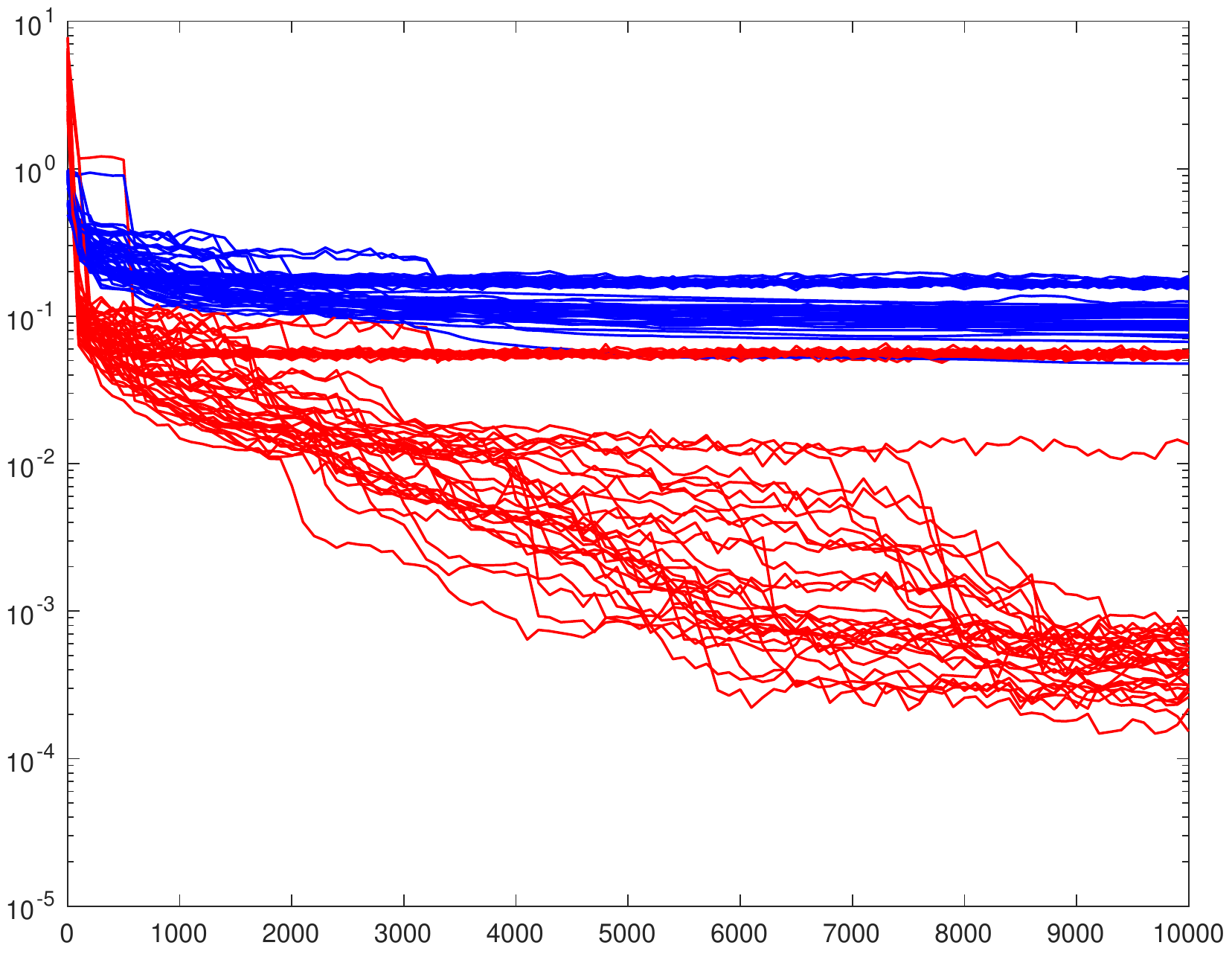}
        \caption{True errors and residuals for $D_{\TA_{\ba} B_{S_3}}$ with
          $k=2, p=3,$ with convolution.}
    \end{subfigure}
\caption{In (a,b) we examine the convergence properties of the
  hybrid map $D_{\TA_{\ba} B_{S_p}},$ on a $k=2$ image
  created without convolution; $p=1,3,$ and, in (c,d) on an image
  created by convolving the $k=0$ image with a Gaussian, chosen so
  that the power spectra matches the $k=2$ case as well as possible.
  See Example~\ref{e:startpt}.
  As
  before, $p=1,3.$ For each plot we used 50 random restarts on a
  $256\times 256$-image; the axes on all four plots are the same.}
\label{Fig6.1}
    \end{figure}  

\begin{example}\label{e:startpt}
  In this example we explore the effects of choosing different
  starting points for a variety of images with different levels of
  smoothness, using an algorithm based on the maps $D_{\TA_{\ba}
    B_{S_p}},$ with $p=1,3.$ Figures~\ref{Fig6}(a,b) shows the
  behavior of 10,000 iterates of these algorithms for a piecewise
  constant image, $k=0,$ and $p=1,3;$ Figures~\ref{Fig6.1}(a,b) are
  similar, but with a smoother image, for which $k=2,$ and $p=1,3.$ In
  this experiment the image is defined as a sum of functions,
  as in~\eqref{eqn2.35.004}--\eqref{eqn2.36.004}, with $k=2.$ Finally
  in Figures~\ref{Fig6.1}(c,d) we show the results of a similar
  experiment where the image is smoothed by convolving with a
  Gaussian.  The width of the Gaussian is selected so that the power
  spectra of the images used in $(a,b)$ and $(c,d)$ are as similar as
  possible. For each image we show the true errors (blue) and residuals
  (red) for 50 random initial conditions.

  For a piecewise constant image ($k=0$) an algorithm based on
  iterating $D_{\TA_{\ba} B_{S_1}}$ seems to converge in essentially
  every trial, albeit with a wide range of rates.  Some of the true error
  curves exhibit scalloping behavior. The linearizations of
  $D_{\TA_{\ba} B_{S_1}}$ at such limiting fixed points have the
  surprising property that they are \emph{not} contractions. Instead
  they are highly non-normal maps, with complex eigenvalues of modulus
  less than $1$.  A very simple example of this phenomenon appears in
  Figure~\ref{f:diffmap}(a).  Other trajectories seems to be
  contracting uniformly toward a fixed point. With a looser support
  constraint, an algorithm based on $D_{\TA_{\ba} B_{S_3}}$ seems to have
  stagnated for all trials, except one. Each of the images $\{\bF^{(\bv)}:\:\|\bv\|_{\infty}\leq
  3\}$ defines an attracting basin for this map. Only at the
  ``corners'' $\{\bF^{(\pm3,\pm3)}\}$ is the intersection with
  $B_{S_3}$ transversal. For the single trial that shows convergence
  $\bv=(-3,3);$ this is the only trial that found an attracting basin
  defined by a transversal intersection.

  It is apparent that a smoother image and/or looser support
  constraint makes it much harder for these algorithms to converge. In
  the experiment whose results are plotted in Figure~\ref{Fig6.1}(a),
  the image has $k=2$ and we use a 1-pixel neighborhood of the true
  support for the support constraint. The scalloping curves strongly
  indicate that the iterates have fallen into an attracting basin
  defined by a transversal intersection, and that these iterates are
  very slowly converging to a fixed point.  The slow convergence
  (relative to the $k=0$ case) is a result of the much smaller
  non-zero angles between $T_{\bF^{(\bv)}}\TA_{\ba}$ and $B_{S_1},$
  caused by the smoothness, even for translates, $\bv,$ where the
  intersection is transversal. The scalloping of the (unobservable)
  true errors is reflected in a similar scalloping in the (observable)
  residuals. In other experiments the iterates seem to be very slowly
  convergent, or perhaps have stagnated. In these cases the residual
  is roughly the square of the true error, indicating approach along a
  direction lying in $T_{\bF^{(\bv)}}\TA_{\ba}\cap B_{S_1}.$
  
  The plots in Figure~\ref{Fig6.1}(b) indicate that all trials have
  stagnated, though they do seem to fall into two distinct groups.  In
  the first group, the true error is close to $1,$ suggesting that the
  iterates have not found an attracting basin defined by a true
  intersection.  In the second group the true errors are below $10^{-1}$
  and the residuals are often even smaller than the squares of the
  true errors.  Empirically, this seems to occur when the iterates find an
  attracting basin, not defined by a true intersection, but close to
  one that is. This sort of behavior can persist for a very large
  number of iterates (millions, at least).

  For the plots in Figure~\ref{Fig6.1}(c,d), we use images that are
  defined by convolving with a Gaussian. As shown in
  Sections~\ref{sec2.3.004}--\ref{sec2.4.004}, this leads to a large
  $\dim T_{\bF^{(\bv)}}\TA_{\ba}\cap B_{S_p},$ even for $\bv$ with
  $S_{\bF^{(\bv)}}\subset S_p.$ The most striking comparison is
  between Figure~\ref{Fig6.1}(a) and Figure~\ref{Fig6.1}(c): in (a)
  essentially all experiments terminate with an true error less than
  $10^{-1}$ and many appear to be slowly converging, with true errors often
  less than $10^{-3},$ whereas in (c), all but 2 cases 
  have stagnated with a true error very close to $10^{-1},$ and a residual
  close to $10^{-4}.$ In one case the true error appears to have stagnated
  at about $10^{-3}$ and in another, the true error is $10^{-2},$ and
  appears to be decreasing geometrically.  In (b) and (d) the
  looseness of the support constraint appears to be the dominant
  source of difficulty.  It is notable how different (a) and (b) are, but
  how similar (c) and (d) are. With a non-transversal intersection,
  more precise support information does little to improve the behavior
  of the algorithm.
\end{example}

From these examples we see that algorithms based on hybrid iterative 
maps often stagnate at a very substantial distance from any
true intersection point. Even for a piecewise constant image, the
iterates stagnate, most of the time, once the support constraint becomes
somewhat imprecise. The quantitative relationship between the true
errors and the residuals often indicates approach along common
tangent directions.

\subsection{The Positivity and $\ell_1$ Constraints}\label{sec5.2}

We turn now to the usage of non-negativity as auxiliary information,
and begin by recalling that non-negativity alone does not suffice for
generic uniqueness up to trivial associates. However, if we also assume
that the autocorrelation image $(\bF\star\bF$) has sufficiently small
support, then this does indeed define an adequate constraint for
$\TA_{\ba}\cap B_+$ to consist of finitely many points, which are
generically trivial associates.  A special case of the uniqueness
result proved in~\cite{BEGM} is
\begin{theorem}
  Let $M$ be a positive integer, let $J=\{-2M,-2M+1,\dots,2M\}^d$, and
  let $\TA_{\ba}$ be the
  magnitude torus defined by a non-negative image $\bF\in\bbR^J$ for which
  $S_{\bF\star\bF}\subset\{-M_0,-M_0+1,\dots,M_0\}^d$,
  where $M_0$ is the largest integer not exceeding $4M/3$.
  Then the intersection $\TA_{\ba}\cap B_+$ consists of finitely many
  points, which, generically, are trivial associates of $\bF$.
\end{theorem}

\noindent
For a non-negative image, $\bF,$ the usual containment
$S_{\bF\star\bF}\subset S_{\bF}-S_{\bF}$ is an equality. This fact
allows one to deduce an upper bound on $S_{\bF}$ from the bound on
$S_{\bF\star\bF}.$ The theorem then follows from Hayes' uniqueness
theorem.

\begin{remark}It should be noted that the autocorrelation image is determined
by the measured data $\{|f_{\bj}|^2:\:\bj\in J\}$ and therefore the
support condition on $\bF\star\bF$ is, in principle, verifiable.  The
theorem is stated for images of size $(4M+1)^d$; there is an analogous
result for images of any size, whose precise statement depends on the
dimensions of the image mod 4; see~\cite{BEGM}.
\end{remark}

The analysis in the case of the support constraint suggests that the
``transversality'' of the intersection at $\bF\in\TA_{\ba}\cap B_+$
will strongly influence the behavior of algorithms based on the map
$D_{\TA_{\ba}B_+}.$ When $\TA_{\ba}\cap B_+$ is finite, this
intersection actually lies in $\pa B_+,$ which is not smooth, but
is a piecewise affine space. Therefore a
reasonable measure of the failure of transversality is
$T_{\bF}\TA_{\ba}\cap \pa B_+.$ It turns out that to study these
intersections it is very helpful to consider the $\ell_1$-norm as a
function on $T_{\bF}\TA_{\ba}.$ In~\cite{BEGM}, it is shown that the
$\ell_1$-norm on this affine subspace assumes its minimum value at
$\bF.$ The intersection at $\bF$ is transversal, i.e. locally
$T_{\bF}\TA_{\ba}\cap B_+=\{\bF\},$ if and only if this is a strict
minimum. This analysis also establishes the equality in
equation~\eqref{eqn42}, and shows that the intersection is a proper
convex cone lying in an orthant of a Euclidean space. The analysis
leads to a practical method for computing these intersections in
concrete examples.

Using this approach, we have considered many examples of the type
defined in~\eqref{eqn2.35.004}--\eqref{eqn2.36.004} with various
values of $k>0,$ and have never found an example with a
non-transversal intersection. We have also carried out these
computations for a collection of $128\times 128$ images, defined by
convolution of a piecewise constant image with $G_k,$ with values of
$k$ ranging from $0$ to $6.$ The results are shown in
Table~\ref{tab2}; the $\dim T_{\bF}\TA_{\ba}\cap B_+$ increases slowly
with $k.$ As $\pa B_+$ is strictly convex in many directions near to
points in $\TA_{\ba}\cap B_+,$ one might expect algorithms based on
$D_{\TA_{\ba}B_+}$ to work better than those based on
$D_{\TA_{\ba}B_{S_p}},$ In the following example we show that, even
for images defined by convolution, this is, indeed, the case.
\begin{table}[h]
\begin{center}
    \begin{tabular}{| c || c | c | c |c|c|c|c|}
    \hline
     $k$ & $0$ & $1$& $2$ & $3$& $4$&$5$&$6$ \\ \hline\hline
     $\dim T_{\bF}\TA_{\ba}\cap\pa B_+$&0 &0 &4 &10 &18 &22 &34 \\ \hline
 \end{tabular}
\end{center}
\caption{Table showing the dimensions of $T_{\bF}{\TA_{\ba}}\cap \pa B_+$ for
  $k=0,1,2,3,4,5,6.$ }\label{tab2}
\end{table}

\begin{example}
{\rm In Figure~\ref{fig6} we show the results of 10,000 iterates of
  $D_{\TA_{\ba}B_+}$ with 25 random starting points for each of four
  $256\times 256$ images (i.e.\ $N=128$).  The images, constructed
  using~\eqref{eqn2.35.004}--\eqref{eqn2.36.004}, have varying degrees
  of smoothness with $k=0,2,4$ in (a), (b), (c), respectively. For
  comparison, in (d) we show the results with an image defined by
  convolution with a Gaussian, where the width is selected so the
  power spectrum is similar to the $k=2$ case.  The plots in (a) show
  geometric convergence, with a wide range of rates. In (b) and (c)
  most examples quickly achieve errors in the $10^{-2}-10^{-3}$ range,  and then the error plots display the characteristic scalloping
  behavior seen in Figure~\ref{Fig6.1}(a). These trajectories are, in
  fact, very slowly converging.  The plots shown in (d) indicate that,
  for images defined by convolution, the algorithm again stagnates,
  though the ultimate true error is a little smaller with the
  positivity constraint than with the support constraint. This
  reflects the fact that $\pa B_+$ is considerably more convex, near
  to a point in $\TA_{\ba}\cap B_+$ than a linear subspace like
  $B_{S_p}.$ Once again the quadratic relationship between the true
  error and the residual indicates that the trajectory ultimately lies
  along a common tangent direction.
  
\begin{figure}[h]
    \centering
    \begin{subfigure}[H]{.4\textwidth}
        \centering
        \includegraphics[height= 4cm]{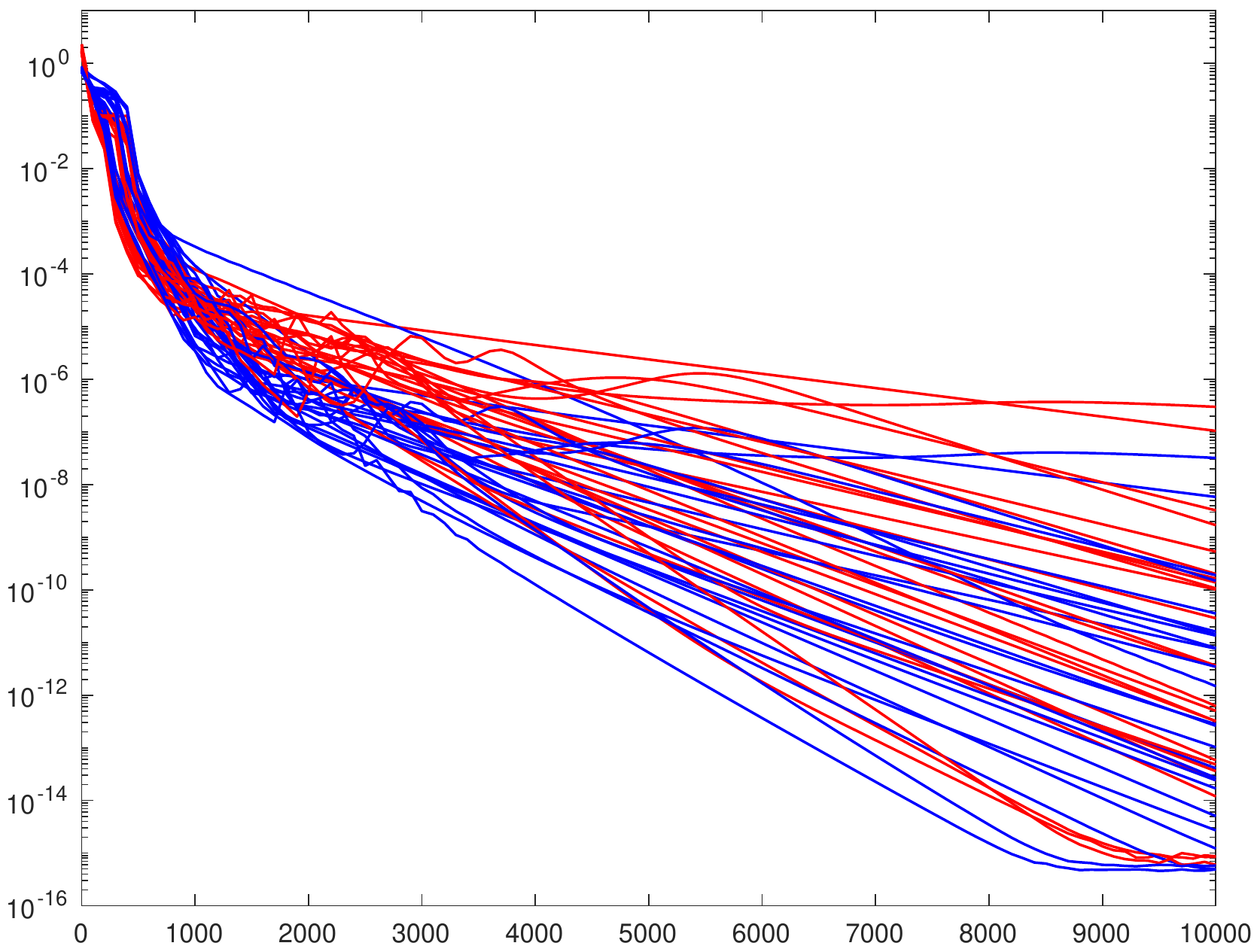}
        \caption{$k=0$}
    \end{subfigure}\qquad
\begin{subfigure}[H]{.4\textwidth}
        \centering
        \includegraphics[height= 4cm]{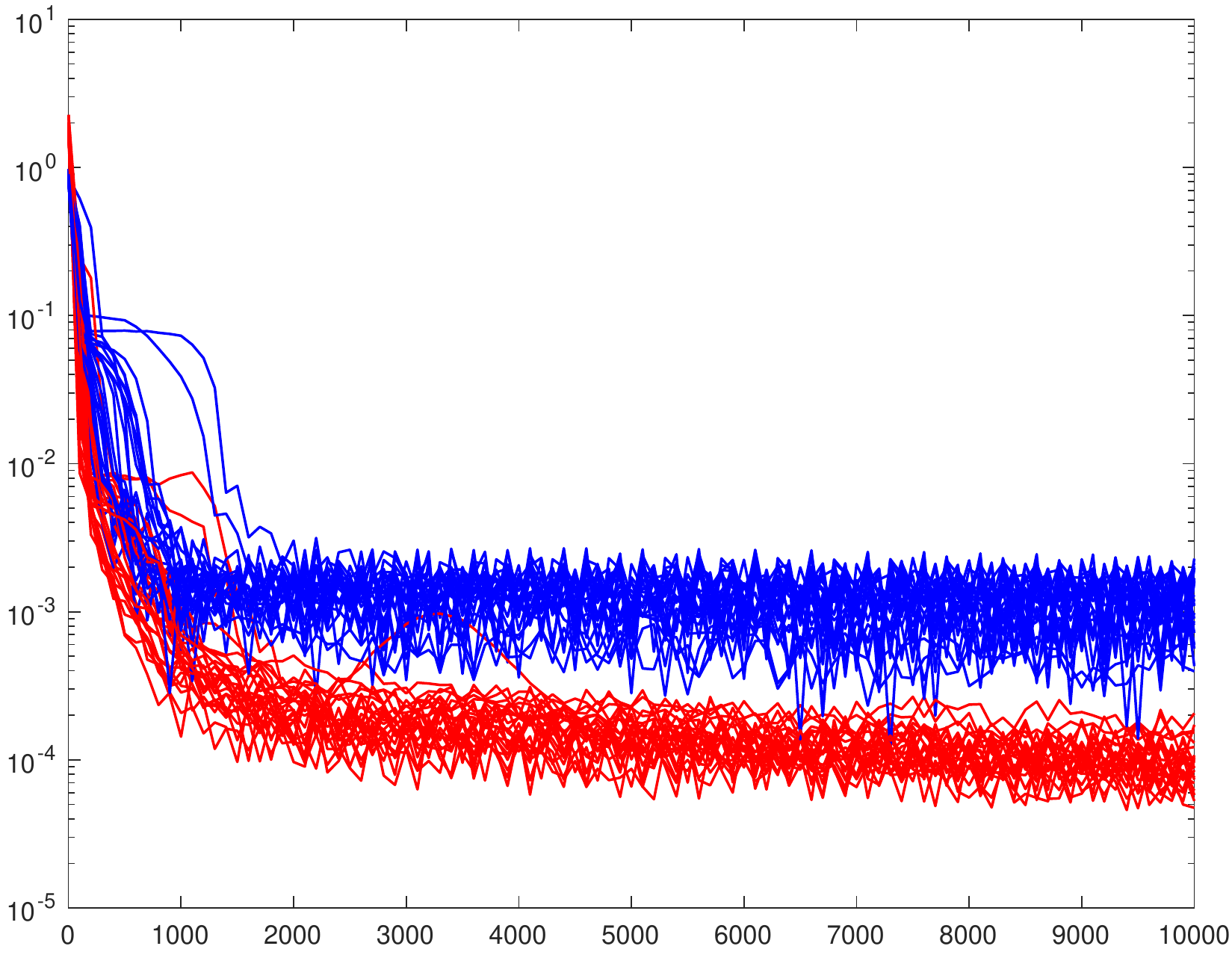}
        \caption{$k=2$}
    \end{subfigure}
\begin{subfigure}[H]{.4\textwidth}
        \centering
        \includegraphics[height= 4cm]{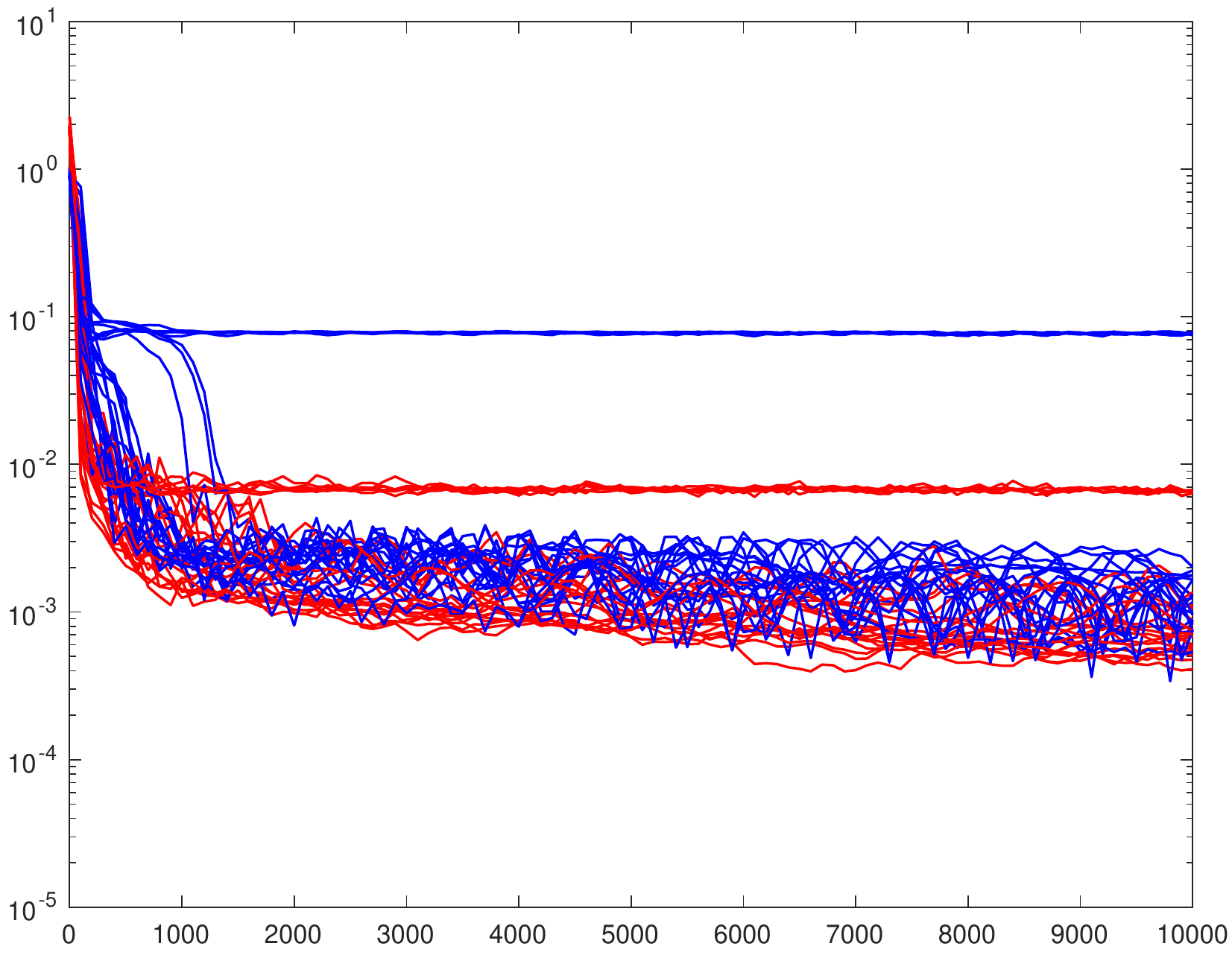}
        \caption{$k=4$}
    \end{subfigure}\qquad
\begin{subfigure}[H]{.4\textwidth}
        \centering
        \includegraphics[height= 4cm]{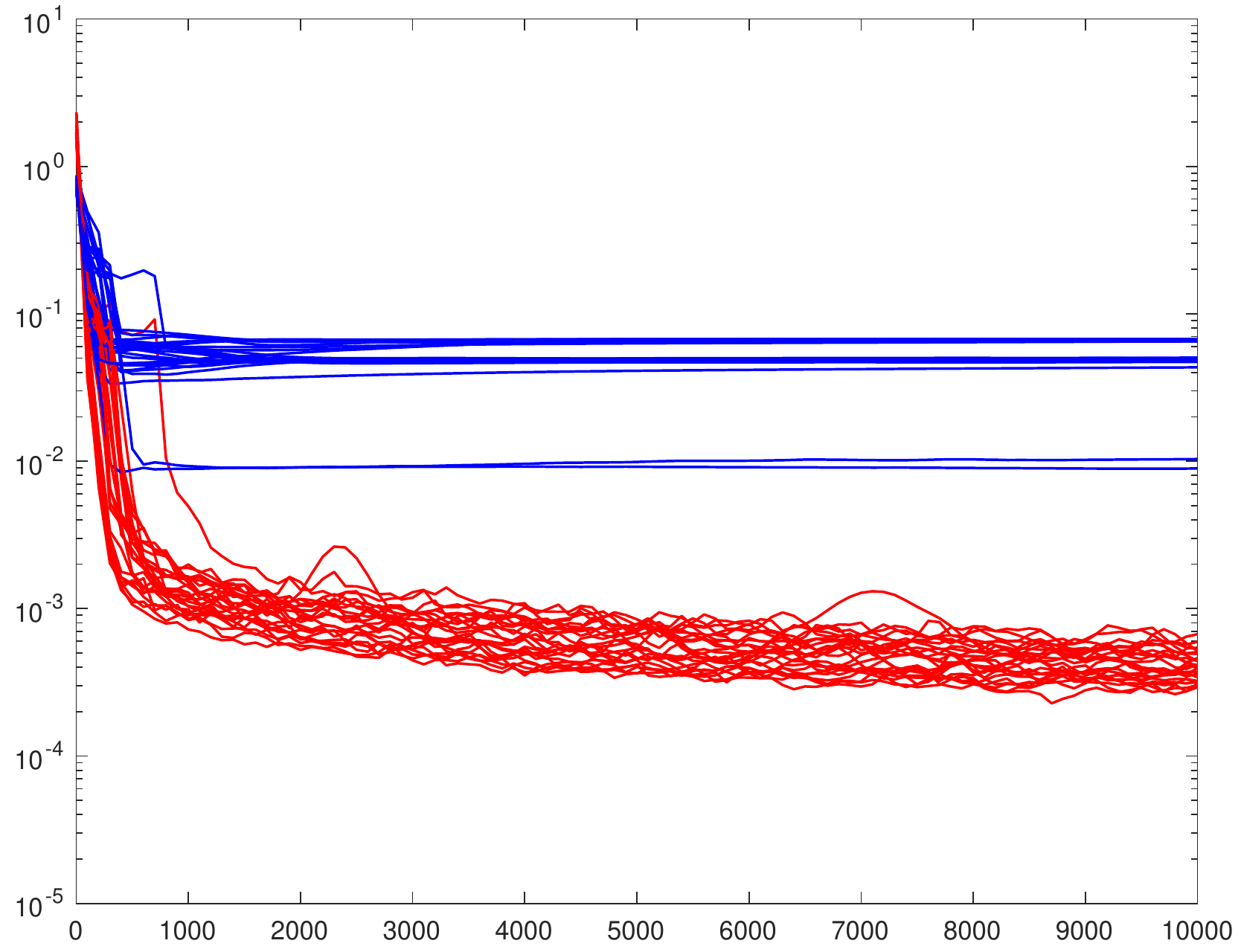}
        \caption{Convolution image}
    \end{subfigure}
    \caption{25 random restarts of 10,000 iterates of $D_{\TA_{\ba}B_+}$ on
      non-negative images with varying degrees of smoothness.}\label{fig6}
\end{figure}
}
\end{example}

\section{Overcoming the difficulty of classical phase retrieval} \label{sec7}

The foregoing sections provide compelling evidence that the intrinsic
difficulty in recovering the phase lies in the local geometry of the
intersections of a magnitude torus $\TA_{\ba}$ with a subset, $B,$
defined by the choice of auxiliary information. To improve the
situation one needs to break what is essentially an infinitesimal
symmetry in order to render these intersections more transversal. 
In practice, this can be achieved by collecting different
experimental data:
pthychography has become an important tool for this, consisting essentially of
rastering across the unknown
image with a mask, making a scattering measurement for each location.
This provides a much larger and richer data set to work with at the cost
of a longer, more involved experiment. Hybrid maps and other 
iterative phase retrieval methods work well with such data sets and converge
quite rapidly. 
Another way to obviate the classical phase retrieval problem is to record
in the near-field of the sample (the Fresnel regime).
For further discussion of pthychography, we refer
the reader to \cite{Dierolf2010,Marchesini2013,Rodenburg2007,Thibault2008} 
and the references therein. For a discussion of the mathematical issues in 
near-field imaging, see \cite{Hohage2017}. 

We limit our attention here to the 
coherent diffraction imaging setting (CDI), since it retains some 
advantages (including speed/timescale of acquisition), and would 
become an even more powerful technology if its associated
phase retrieval problem could be 
addressed robustly.
We propose two experimental modifications which could attain that end.

\subsection{Sharp Cut-off Mask} 
In biological applications one is
often seeking to image a sample of soft tissue. If one could 
cut the sample along a sharp edge, the object would be non-smooth. Moreover,
knowledge of the precise shape would provide for an 
accurate support constraint, which would in
turn break the infinitesimal translational symmetry. In practice, it
is better to use a mask that is not invariant under the
inversion symmetry; see~\eqref{eqn3}.

\begin{figure}[H]
    \centering
    \begin{subfigure}[t]{.85\textwidth}
      \centering
                \includegraphics[width= 10cm]{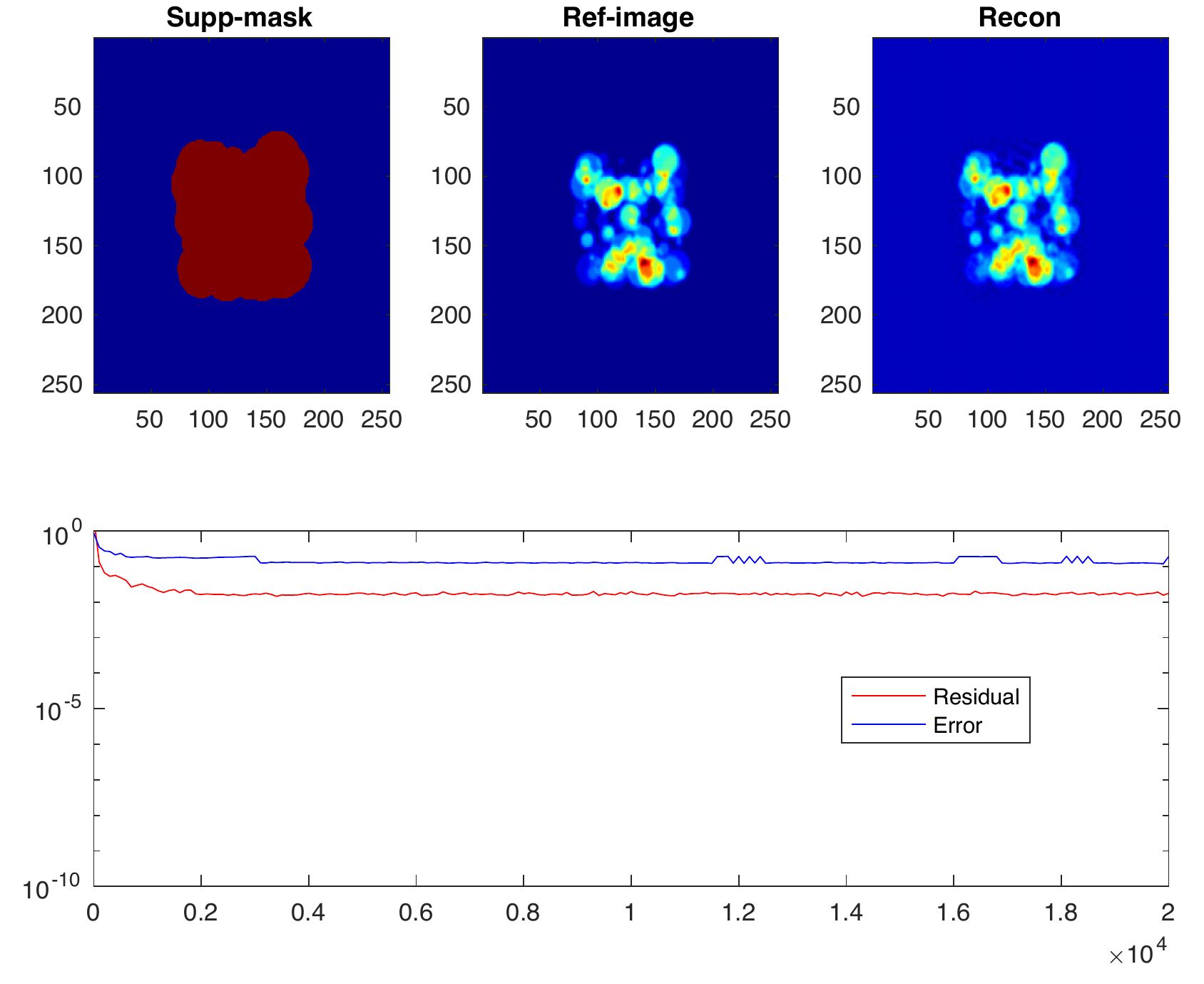}
        \caption{ 20,000 iterates of $D_{\TA_{\ba}B_{S_1}}$ for a soft
          object ($k=6$). }
    \end{subfigure}
    \begin{subfigure}[t]{.85\textwidth}
      \centering
                \includegraphics[width = 10cm]{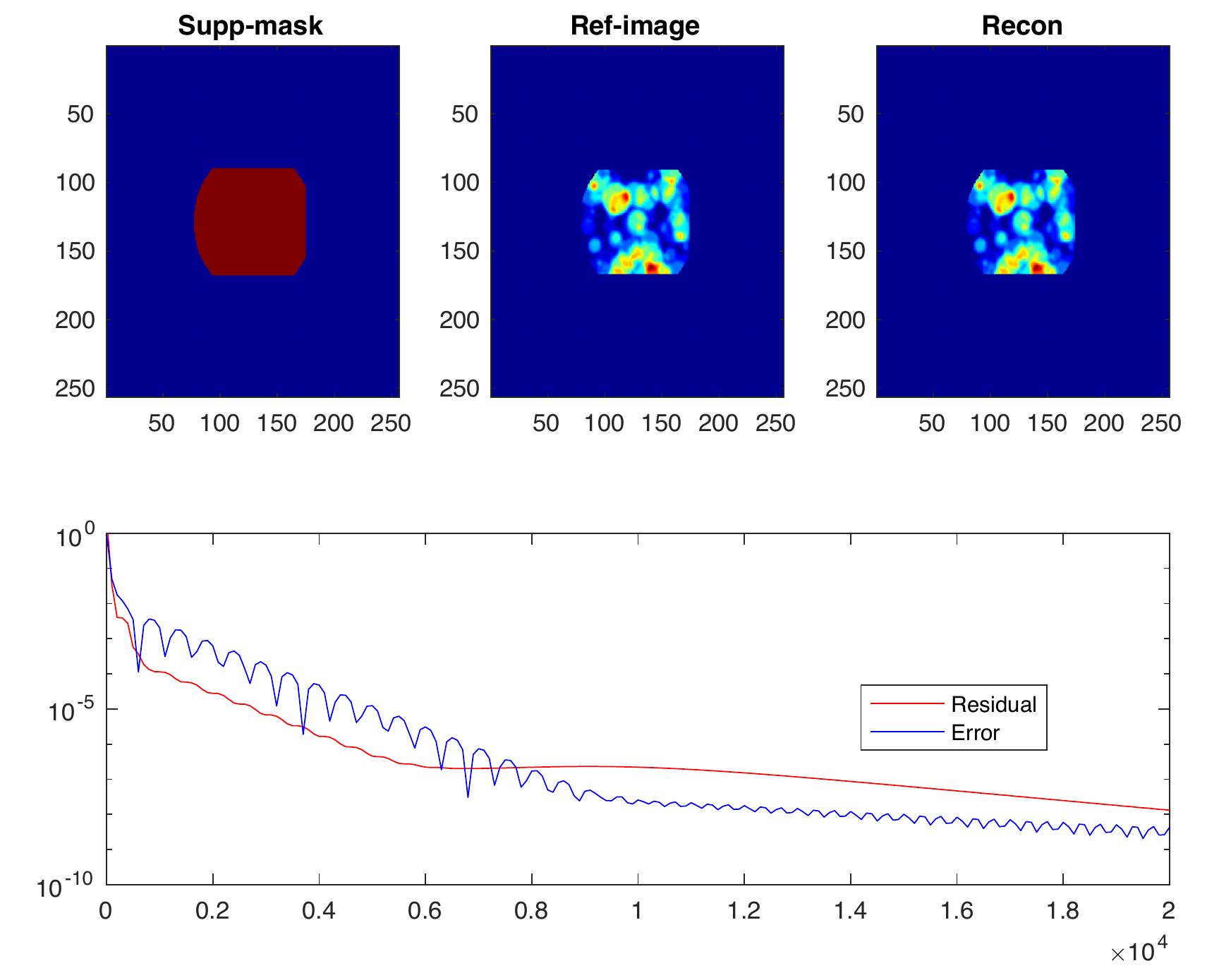}
        \caption{20,000 iterates of $D_{\TA_{\ba}B_{S_1}}$ for a soft
          object ($k=6$) cut-off with a sharp mask.}
    \end{subfigure}
    \caption{An illustration of how the convergence properties of a
      hybrid iterative map based on $D_{\TA_{\ba} B_{S_1}}$ are affected 
      by a sharp
      cut-off mask. The images in the top row, from left to right are the
      support-mask, the reference image, the reconstructed image. The bottom
      plots show the true error (blue) and the residual \eqref{E} (red).}
\label{fig8}
\end{figure}

In spite of the fact that the material being imaged may be soft,
it is possible to
obtain very high resolution.  The examples shown in Figure~\ref{fig8} were
selected as the best results from 20 random initial conditions, for each of the
two experimental set-ups.  Figure~\ref{fig8}(a) shows the result obtained
when running an algorithm defined by the map $D_{\TA_{\ba}B_{S_1}}$ on the
smooth image without a sharp cut-off. The set $S_1$ is the 1-pixel neighborhood
of the set $\{\bj:\:|f_{\bj}|\geq 10^{-10}\}.$ As expected the iterates
quickly stagnate, whereas, in Figure~\ref{fig8}(b) we see that a sharp cut-off
allows for geometric convergence, where we use as a support constraint
the 1-pixel neighborhood of the region bounded by the sharp cut-off. With a
2-pixel support neighborhood the performance degrades markedly. For both of
these experiments the 19 other runs yielded results that were only slightly worse.

\subsection{External Holography}
A second experimental modification (and perhaps one that is easier to
carry out), consists of what we will refer to as {\em external
  holography}.  For this, we imagine placing a \emph{known} hard
object in the exterior of the (perhaps soft) object one would like to
image.  Several related ideas appear in the literature. One is called
double blind Fourier holography, and was recently considered
in~\cite{leshem2016,raz2014}, with a reconstruction method based on a
mixture of Fourier and linear algebraic ideas. Another approach, using
more complex reference objects, is found in the recent work of
Barmherzig, Cand\`es, et al.,
see~\cite{Candesetal2019a,Candesetal2019b}. The reconstruction method
in this approach is largely algebraic.

\begin{figure}[H]
    \centering
                \includegraphics[width = 10cm]{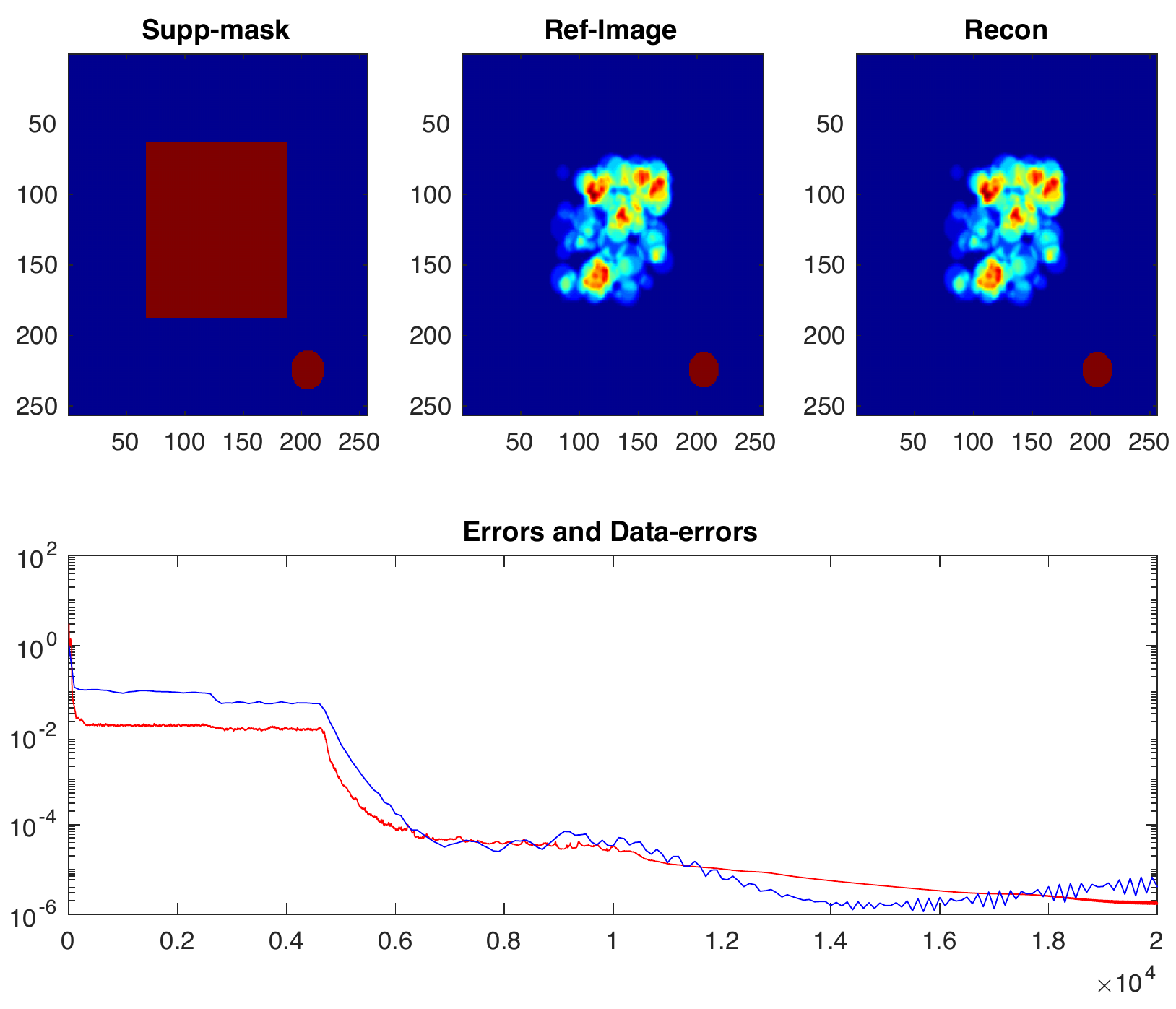}
\caption{An illustration of how  the convergence properties of the
  hybrid map $D_{\TA_{\ba} B_{S_1}}$ are improved by the
  addition of a hard external object. See the caption of the previous
  figure for explanation of the plots and curves.}
\label{fig9}
\end{figure}

Here, we simply make use of the hybrid map based on 
$D_{\TA_{\ba} B_{S_1}},$ where the support mask is the 1-pixel neighborhood
of the smallest rectangle that encloses the object, along with the
1-pixel neighborhood of the exterior hard object.
The {\em shape} of the external object must be precisely
known; its location is less important and can be determined as part of
an iteration step. 
As we see in Figure~\ref{fig9}, which uses the same soft
object as that employed in Figure~\ref{fig8}(a), inclusion of the external object
leads to geometric convergence.
Using the 2-pixel neighborhood of the external
object leads to results similar to those in Figure~\ref{fig8}(a). Once
again, we have shown the best outcome obtained from 20 independent
trials. Some of the other trials gave markedly worse results than the
one shown here. If we use the exact support of the external object,
then the results consistently resemble those in Figure~\ref{fig8}(b).

It is worth noting that 
external holography overcomes the microlocal
non-uniqueness described in Section~\ref{sec3.2}. For simplicity, suppose the
$\bF_0$ has a decomposition as $\bF_0=\bF_1+\bF_2,$ where $\bF_1$ and $\bF_2$
are compact objects for which the supports of $\hbf_1$ and $\hbf_2$ are
disjoint to a high degree of accuracy. Let $\bg$ denote the external object;
for translation directions $\bv_1,\bv_2,$ we observe that,
for all $\bk$,
\begin{multline}\label{eqn44}
  |f_{1\bk}+f_{2\bk}+g_{\bk}|^2-|f^{(\bv_1)}_{1\bk}+f^{(\bv_2)}_{2\bk}+g_{\bk}|^2
=\\
2\Re\left[f_{1\bk}\overline{g_{\bk}}(1-e^{\frac{\pi i\bk\cdot\bv_1}{N}})
+f_{2\bk}\overline{g_{\bk}}(1-e^{\frac{\pi i\bk\cdot\bv_2}{N}})\right]
+O(|f_{1\bk}f_{2\bk}|).
\end{multline}
Assuming that DFT coefficients $g_{\bk}$ decay slowly, the first term on the right
hand side of~\eqref{eqn44} is typically many orders of magnitude
larger than the error term $O(|f_{1\bk}f_{2\bk}|).$ This shows that
a hard external object effectively breaks the microlocal translational
symmetry that leads to $\epsilon$-non-uniqueness.

\section{Conclusions}

In this paper, we have described a framework for analyzing the
classical phase retrieval problem, where only the magnitude of the
Fourier transform of an unknown object is measured, typically in
combination with some information about its support.  Perhaps most
alarmingly, we have shown that the problem is classically
ill-posed---that is, with typical support information, the locally
defined inverse map is only H\"older continuous. Moreover, one can
easily construct objects that are quite distinct, but have 
supports, and magnitude DFT data that are indistinguishable to any
precision $\epsilon>0$.  While some such counterexamples are clearly
pathological, others are not (as shown in section
\ref{nonuniquesec}). This leads to two open mathematical questions:
how dense is the set of $\epsilon$-non-unique objects in the space of all
objects, and can one determine, for a given data set, whether phase
retrieval is even possible, at a given precision?

Assuming that, in the generic case, the problem is solvable, we have
also shown that phase retrieval typically involves finding the
intersection of two sets that do not meet transversally. It is
precisely this failure of transversality that prevents the local
inverse from being Lipschitz continuous. The formal linearization at a
non-transversal intersection has infinite condition number. Beyond
this, small angles between the set $B,$ defined by the auxiliary data,
and the fiber of the tangent bundle to $\TA_{\ba},$ at the
intersection point, prevent standard iterative methods from
converging. These effects are mitigated by having an object with a
sharp boundary and accurate support information. As the external
holography example shows, it suffices to have these properties for a
component of the object being imaged.

The mathematical foundations of this paper are presented in detail in \cite{BEGM},
and we are currently working on modifications of the experimental protocol
(other than pthychography or near-field imaging,  see \cite{Hohage2017})
that will lead to better-posed inverse problems. The results of that work will
be reported at a later date.

\bibliography{phaseret,alla-k}
\bibliographystyle{siam}

\end{document}